\documentclass[final,5p,times,twocolumn,authoryear]{elsarticle}

\usepackage{amssymb}
\usepackage{amsfonts}
\usepackage{latexsym, amssymb,amsmath, amsbsy,amsopn, amstext}
\usepackage{amsmath, amsbsy}
\usepackage{hyperref}
\usepackage[none]{hyphenat}
\usepackage{float}
\usepackage{ifpdf}
\usepackage{xcolor}
\usepackage{tikz}
\usetikzlibrary{arrows,shapes,snakes,shadows,positioning,automata,patterns}
\usetikzlibrary{trees,decorations.pathmorphing,decorations.markings}
\def\diag{diag}

\def\0{{\bf 0}}

\newtheorem{thm}{Theorem}[section]
\newtheorem{lem}[thm]{Lemma}

\newtheorem{exa}[thm]{Example}
\newtheorem{rem}[thm]{Remark}
\newtheorem{alg}[thm]{Algorithm}
\newtheorem{assum}{Assumption}

\journal{Elsevier}

\begin{document}
\begin{frontmatter}

\title{A distributed primal-dual algorithm for computation of generalized Nash equilibria  with shared affine coupling constraints via operator splitting methods}
\tnotetext[t1]{This work was supported by NSERC Discovery Grant (261764).}

\author[label1]{Peng Yi\corref{cor1}, Lacra Pavel}
\address[label1]{Department of Electrical and Computer Engineering, University of Toronto}
\ead{peng.yi@utoronto.ca, pavel@control.toronto.edu}

\begin{abstract}
In this paper, we propose a distributed primal-dual algorithm for computation of  a generalized Nash equilibrium (GNE) in noncooperative games over network systems. In the considered game, not only each player's  local objective function depends on other players' decisions, but also the feasible decision sets of all the players are coupled together with a globally shared affine inequality constraint. Adopting the variational  GNE, that is the solution of a variational inequality, as a refinement of GNE, we introduce a  primal-dual algorithm that  players can use to seek it in a distributed manner. Each player only needs to know its local objective function, local feasible set, and a local block of the affine constraint. Meanwhile, each player only  needs to observe the decisions on which its local objective function explicitly depends through the interference graph and share information related to multipliers with its neighbors through a multiplier graph. Through a primal-dual analysis and an augmentation of variables, we reformulate the problem as finding the zeros of a sum of monotone operators.   Our distributed primal-dual algorithm is based on forward-backward operator splitting methods. We prove its convergence to the variational  GNE for fixed step-sizes under some mild assumptions. Then a distributed algorithm with inertia is also introduced and analyzed for variational  GNE seeking. Finally, numerical simulations for network Cournot competition are given to illustrate the algorithm efficiency and  performance.
\end{abstract}


\end{frontmatter}



\section{Introduction}

%decision problems in network systems. solution concept, optimization-- centralized optimization and distributed optimization
Engineering network systems, like power grids, communication networks, transportation networks and sensor networks,
play a foundation role in modern society.  The efficient and secure operation of various network systems relies on efficiently  solving
 decision and control problems arising in those large scale network systems. In many decision problems, the nodes can be regarded as agents that  need to make  local decisions possibly limited by the shared network resources within  local feasible sets. Meanwhile, each agent has a local cost/utility function to be optimized, which depends on the decisions of other agents. The traditional manner for solving such decision problems over networks is the {\it centralized optimization approach}, which relies on a control center to gather the data of the problem and  to optimize the social welfare (usually taking the form of the sum of local objective functions) within the local and global constraints.  The centralized optimization approach may not be suitable for decision problems over large scale networks, since it needs bidirectional communication between all the network nodes  and the control center,  it is not robust to the failure of the center node, and the computational burden for the center is unbearable. It is also not preferable because the privacy of each agent might be compromised when the data is transferred to the center. Recently,  a {\it distributed optimization approach} is proposed as an alternative methodology for solving decision problems in network systems (\cite{yipeng}, \cite{shi} and \cite{zeng2016}). In the distributed optimization approach, the data is distributed throughout the network nodes and there is no control center, and each agent in the network can just utilize  its local data and share information with its neighbour agents to compute  its local decision that corresponds to the optimal solution of the social welfare optimization problem. Therefore, the distributed optimization approach  overcomes the drawbacks  of the centralized optimization approach by decomposing the data, computation and communication to each agent.
Moreover,  each agent has the authority and autonomy to formulate its own objective function without worrying about privacy leaking out.  However, both approaches adopt the same solution concept, that is the optimal social welfare solution with local and global constraints, as the solution criterion of  decision problems in network systems.

%why NE is a more desirable solution
However, optimal solutions of social welfare may not be proper solution concepts  in many applications.
In fact,  with the deregulation and liberalization of markets, there is no guarantee that the agents will not deliberately deviate from their local optimal solutions to increase (decrease) own local utility (cost), possibly by deceiving to utilize more network resources. In this paper, we consider the game theoretic approach where each agent in the network  has its own local autonomy and rationality.  In such a setup of  multiple interacting rational players making decisions in a noncooperative environment, Nash Equilibrium (NE) is a more reasonable solution. In an NE, no player can increase (decrease) its local utility (cost) by unilaterally changing its local decision, therefore, no agent has the incentive to deviate from  it. In other words,  NE is a self-enforceable solution in the sense that once NE is computed all the agents will execute that NE.
Recently, there has been  increasing interest in adopting  game theory  and NE as the modeling framework and solution concept for various network decision problems, like wireless communication systems  (\cite{pang3} and \cite{asuman}), network flow control (\cite{baser}), optical networks (\cite{pavel2}) and smart grids (\cite{hu}).

%In some cases, GNE should be considered
Moreover,  in engineering network systems, not only the  local objective function of each agent depends on other players' decisions, but also the feasible set of each local decision could depend on other agents' decisions, because the agents may compete for the utilization of some shared or common network resources like bandwidth, spectrum or power.
This type of network decision problems can be modeled as {\it noncooperative games with coupling constraints},  and {\it generalized Nash equilibrium}  (GNE) can be adopted as its solution. The study of GNE dates back to the social equilibrium concept proposed by \cite{debreu},
and flourished  in the last two decades with  applications in practical problems like environment pollution games (\cite{krawczyk2}), power market design (\cite{krawczyk1}), optical networks (\cite{pavel1}), wireless communication (\cite{pang}).
Interested readers can refer to
\cite{faccinei1} for a historical review of GNE, and refer to  \cite{faccinei2} for recent developments, and  to \cite{pang} for a technical treatment.

%how to compute NE and GNE, traditional method, and distributed method
Even though NE or GNE  is a reasonable and expectable solution as  a result of  multiple  rational agents  making decisions in a noncooperative manner, how to arrive at an  NE (GNE) is by no way a self-evident task.
Each player needs to know the complete game information, including objective functions and feasible decision sets of all the other agents, in order to compute NE in an introspective manner.
It gets much more complicated for computing GNE, because the agents also have to consider the coupling in the feasible decision sets. In fact, as of yet there is no universal manner to efficiently compute GNE in games with coupling constraints (\cite{harker} and \cite{faccinei2}),  except for games with shared coupling constraints (\cite{facchinei3}).
Moreover,  for games in large scale network systems, it is quite unrealistic and undesirable to  assume that each  agent could have the complete information of the whole network,  because this  implies  prohibitive communication and computation burden and no privacy protection.
{\em Therefore, each player (agent) should compute  its local decision corresponding with  an NE or GNE in a distributed manner, somehow resembling the distributed optimization approach.
In other words, each agent should only utilize its local objective function, local feasible set and possible local data related to coupling constraints,  and should only share information with its neighbouring agents to compute its local decision in the NE (GNE).}  This turns out to be an emerging research topic and gets studied in  \cite{pavel3}, \cite{shanbhag1}, \cite{lygeros1}, \cite{hu} and \cite{kar1}, etc.

%The concern of this paper, the problem, the solution concept
Motivated by the above,  in this work we consider a distributed algorithm for iterative computation of  GNE in  noncooperative games with shared  affine coupling constraints over network systems.
The considered noncooperative game has each agent's local objective function depending on  other agents' decisions as specified by an {\it interference graph},
and has also an affine constraint shared by all agents,  coupling all  players'  feasible decision sets.
The considered game model covers many practical problems, like the power market model in  \cite{krawczyk1}, environment pollution game in \cite{krawczyk2}, power allocation game in communication systems in \cite{shanbhag4}. A (centralized) numerical algorithm was recently studied in \cite{shanbhag3} for quadratic objective functions and in \cite{dreves} for linear objective functions.
% solution concept
Generally speaking, the GNE of the considered   game  may not be unique.
In this work, we adopt the variational GNE, that corresponds with the solution of a variational inequality  proposed in \cite{facchinei3}, to be  a refinement GNE solution. The variational GNE is a particular type of the normalized equilibrium proposed in \cite{rosen}, and enjoys  a nice economical interpretation that all the agents have the same shadow price for shared network resources without any discrimination as pointed in \cite{shanbhag2}. Furthermore, the variational GNE enjoys a sensitivity and stability property  (\cite{pang} and \cite{faccinei1}),  hence we adopt it as the desirable solution.

%The algorithm, the contributions....
We propose a new type of distributed algorithm that agents can use to compute the variational GNE by only manipulating their local data and communicating with neighbouring agents. Observing that
the KKT condition of the corresponding variational inequality requires all agents to reach consensus on
the multiplier of the shared affine constraint, we introduce a local copy of  the multiplier and an auxiliary variable for each player. To enforce the  consensus of local  multipliers, we  use a reformulation that incorporates the {\it Laplacian matrix} of a connected  graph. Motivated by the {\it forward-backward operator splitting method} for finding zeros of  a sum of monotone operators (refer to \cite{combettes1}) and the recent primal-dual algorithm proposed in \cite{condat} for optimization problems with linear composition terms,  we propose a novel distributed algorithm for iterative computation of GNE.
The main idea is to introduce a suitable metric matrix  and to split the equivalent reformulation into two monotone operators. An operator splitting method has been adopted for NE computation in a centralized manner in \cite{combettes2}.
A different splitting idea is adopted here appropriate for distributed GNE computation.  Moreover,  a distributed algorithm with inertia is also proposed and investigated,  motivated by the acceleration algorithms in \cite{attouch1}, \cite{attouch2}, \cite{hendrickx} and \cite{intertia1}, most of which only focused on optimization problems.
The  convergence of the proposed algorithms is  verified under suitable fixed step-size choice and some mild assumptions on the objective functions and communication graphs.

% review of some related works
The recent works of \cite{zhuminghui},  \cite{sayed},  \cite{liangshu} and  \cite{lygeros2} are closely related with this work
since all of them are concerned with  the distributed algorithm for  seeking GNE of noncooperative games with coupling constraints.
\cite{zhuminghui} address the GNE seeking for  the case where each player has non-shared local coupling constraints. Assuming that each player can observe other players' decisions on which its local objective function and local constraint functions depend through the interference graph, \cite{zhuminghui} propose a distributed primal-dual GNE seeking algorithm based on variational inequality methods, and show algorithm convergence  under diminishing step-sizes.
\cite{sayed} investigate the distributed GNE seeking under stochastic data observations.
The authors assume that the coupling constraints  have a locally shared property that if one player
has its one of local constraints dependent on the decision of another player, then this constraint must be shared between those two players. Their algorithm design is based on a penalty-type gradient method.
Under the assumption that each player can observe the decisions on which its local objective function and constraint functions depend, \cite{sayed} utilize a gradient type algorithm to seek the pure NE of the game derived by penalizing the coupling and local constraints.
They show that their algorithm can reach a  region near the pure penalized NE with a constant step-size, which will approach a GNE  if the penalizing parameter goes to infinity.
Both \cite{liangshu} and  \cite{lygeros2} consider the  distributed algorithm  for seeking a variational GNE of the {\em aggregative game} with globally shared affine coupling constraints. This represents a particular type of game where the players' local objective functions  depend on  some aggregative variables of all agents' decisions.
\cite{liangshu} assume  that each player  has  local copies  of both the aggregative variables  and  the multipliers,
and combine a finite-time convergent {\it continuous-time} consensus dynamics and a projected {\it gradient flow} to derive their distributed GNE seeking dynamics.
Meanwhile, \cite{lygeros2} adopt the asymmetric projection algorithm for variational inequalities to design their variational GNE
seeking algorithm. However \cite{lygeros2} assume that there is an additional central node for the update of the common multiplier, and only address quadratic objective functions.

Compared with these works, our paper has following contributions,

(i): The considered noncooperative game model is completely general, thus a  generalization of the {\it aggregative game} in \cite{liangshu}  and \cite{lygeros2}.   We further assume that the shared affine coupling constraint is also decomposed such that each player only knows its local contribution to the global constraint, that is only a sub-block matrix of the whole constraint matrix. In this sense, no player knows exactly the shared constraints, hence,  our problem model is also different from  the ones in \cite{sayed} and \cite{zhuminghui}. The decomposition of the coupling constraints, together with the localization of player's local objective function and local feasibility set, is quite appealing for iterative computation of GNE in large-scale network systems because this reduces the data transmission and computation burden, and protects the players' privacies.

(ii):The proposed distributed algorithms can compute the variational GNE iteratively under a more localized data structure and information observing structure compared to previous ones.  Firstly, each player only utilizes the local objective function and local feasible set, and its local sub-block matrix of the affine constraints.
  %  The decomposition of affine constraints has not been fulfilled with discrete-time algorithms,  noticing that  \cite{liangshu} achieved this with continuous-time dynamics.
     Secondly, we assume the players have two (different) information observing graphs, i.e.,  {\em interference graph} and {\em multiplier graph}. Each player only needs to observe the decisions that its local objective function {\it directly} depends on through the interference graph.  This type of information observation assumption has also been adopted in  \cite{sayed} and \cite{zhuminghui}. Meanwhile, each player only needs to share information related to multipliers with its neighbouring agents through another multiplier graph. Here it is not required that each player should know the decisions that coupling constraints depend on, as assumed in \cite{sayed} and \cite{zhuminghui}. Therefore, our information sharing (observing) structure is more localized and sparse.

(iii): The algorithm development and convergence analysis is motivated by the operator splitting method (\cite{combettes1}), different from the penalized method adopted in \cite{sayed} and  the variational inequality approach in \cite{liangshu}, \cite{zhuminghui} and \cite{lygeros2}.
    Based on this  operator splitting approach, we prove the algorithm converges to the variational GNE under fixed step-sizes.
    Note that neither convergence nor non-bias  estimation is achieved in \cite{sayed} under fixed step-sizes, while \cite{zhuminghui} achieve convergence with diminishing step-sizes.
    On the other hand, compared with \cite{combettes2}, this paper addresses the GNE seeking under coupling constraints,   adopts a different splitting technique, and achieves fully distributed computations. The operator splitting method is  powerful and provides additional insights. Moreover,  a distributed algorithm with inertia is proposed and analyzed, resembling the acceleration algorithms in optimization (\cite{nesterov} and \cite{hendrickx}). The algorithm performance is illustrated via  numerical experiments of network Cournot  competitions with bounded market capacities.

The paper is organized as follows.
Section \ref{sec_notation} gives the notations and preliminary background.
Section \ref{sec_game_and_algorithm} formulates the noncooperative game and gives the distributed algorithm for iterative computation of a GNE.
Section \ref{sec_algorithm_develpo} shows how the operator splitting  method motivates the  algorithm development, and Section \ref{sec_alg_converge} presents the algorithm convergence analysis. Then a distributed GNE seeking algorithm with inertia is proposed and analyzed in Section \ref{sec_alg_inertial}.
Finally, a network Cournot competition with bounded market capacities is formulated with numerical studies in Section \ref{sec_cournot_simulaton}, while concluding remarks are given in  Section \ref{sec_concluding}.

\section{Notations and preliminary background}\label{sec_notation}

In this section, we review the notations and  some preliminary notions in  monotone operators and graph theory.

{\it Notations}:  In the following,
$\mathbf{R}^m$ ($\mathbf{R}^m_{+}$) denotes the $m-$dimesional (nonnegative) Euclidean space.
For a column vector $x \in \mathbf{R}^m$ (matrix $A\in \mathbf{R}^{m\times n}$),
$x^T$ ($A^T$)  denotes its transpose.
$x^Ty=\langle x, y\rangle$ denotes the inner product of $x,y$, and $||x||_2= \sqrt{x^Tx}$ denotes the norm induced by inner product $\langle\cdot,\cdot\rangle$.
Given a symmetric positive definite matrix $G$, denote the $G$-induced inner product $\langle x,  y\rangle_{G}=\langle Gx,y\rangle$. The $G$-matrix induced norm, $|| \cdot||_{G}$,  is defined as $||x ||_G=\sqrt{\langle Gx,x\rangle}$.
Denote  by $||\cdot||$  any matrix induced norm in the Euclidean space.
Denote $\mathbf{1}_m=(1,...,1)^T \in \mathbf{R}^m$ and
$\mathbf{0}_m=(0,...,0)^T \in \mathbf{R}^m$. For column vectors $x, y$, $x \geq (>) y$ is understood  componentwise.
$diag \{A_1, . . . ,A_N\}$ represents
the block diagonal matrix with matrices $A_1, . . . ,A_N$ on its
main diagonal. $Null(A)$ and $Range(A)$ denote the null space and range space of matrix $A$, respectively.
Denote $col(x_1,....,x_N) $ as the column vector stacked with column vectors $x_1,...,x_N$.
$I_n$ denotes the identity matrix in $\mathbf{R}^{n\times n}$.
For a matrix $A=[a_{ij}]$, $a_{ij}$ or $[A]_{ij}$
stands for the matrix entry in the $i$th row and $j$th column of $A$. We also use
$[x]_{k}$ to denote the $k-$th element in column vector $x$.
Denote $\times_{i=1,...,n}\Omega_i$ or $\prod_{i=1}^n \Omega_i$ as the Cartesian product of the sets $\Omega_i,i=1,...,n$.
Denote $int(\Omega)$ as the interior  of $\Omega$, and $bd(\Omega)$ as the boundary set of  $\Omega$.
Define the projection of $x$ onto a set $\Omega$ by $P_{\Omega}(x)=\arg\min_{y\in \Omega} ||x-y ||_2$.
A set $\Omega$ is a convex set if $\lambda x +(1-\lambda)y\in \Omega,\; \forall \lambda \in [0,1], \forall x,y \in \Omega$.
An extended value proper  function $f: \mathbf{R}^{m}\rightarrow \mathbf{R}$ is  a convex function if
$f(\lambda x +(1-\lambda)y ) \leq \lambda f(x) + (1-\lambda) f(y), \forall \lambda \in [0,1], \forall x,y \in dom f $.

\subsection{Monotone operators}
The following concepts are reviewed from \cite{combettes1}.
Let $\mathfrak{A}:\mathbf{R}^m \rightarrow 2^{\mathbf{R}^m}$ be a set-valued operator. Denote ${\rm Id}$ as the identity operator, i.e, ${\rm Id}(x)=x$.
The domain of $\mathfrak{A}$ is $dom\mathfrak{A}= \{x\in \mathbf{R}^m| \mathfrak{A}x \neq \emptyset\}$ where $\emptyset$ stands for the empty set, and the range of $\mathfrak{A}$ is $ran\mathfrak{A}=\{y \in \mathbf{R}^m|  \exists x\in \mathbf{R}^m, y\in \mathfrak{A}x\}$. The graph of $\mathfrak{A}$ is $gra\mathfrak{A}=\{(x,u) \in \mathbf{R}^m\times \mathbf{R}^m| u\in \mathfrak{A}x\}$, then the inverse of $\mathfrak{A}$ is defined through its graph as $gra\mathfrak{A}^{-1}=\{(u, x)  \in \mathbf{R}^m\times \mathbf{R}^m| (x, u)\in gra \mathfrak{A}\}$.
The zero set of operator $\mathfrak{A}$ is $zer\mathfrak{A}=\{x\in \mathbf{R}^m | \mathbf{0} \in \mathfrak{A}x\}$.
Define the resolvent of operator $\mathfrak{A}$ as $R_{\mathfrak{A}}=({\rm Id}+\mathfrak{A})^{-1}$.
An operator $\mathfrak{A}$ is called monotone if
$\forall (x,u), \forall(y,v)\in gra\mathfrak{A}$, we have
$\langle x-y, u-v\rangle \geq 0.$
Moreover, it is maximally monotone if $gra\mathfrak{A}$ is not strictly contained in the graph of any other monotone operator. $R_\mathfrak{A}$ is single-valued and $domR_{\mathfrak{A}}=\mathbf{R}^m$ if $\mathfrak{A}$ is maximally monotone \footnote{ Proposition 23.7 in \cite{combettes1}}.
For a proper {\it lower semi-continuous convex} (l.s.c.) function $f$, its sub-differential is a set-valued operator $\partial f: domf\rightarrow 2^{\mathbf{R}^m}$ and
\begin{equation}\label{subdifferential}
\partial f:\; x \mapsto \{g\in \mathbf{R}^m | f(y)\geq f(x)+ \langle g, y-x \rangle, \forall y\in domf\}.
\end{equation}
$\partial f$ is a maximally monotone operator \footnote{   Theorem 20.40  in \cite{combettes1} }.
Then $Prox_{f}= R_{\partial f}:\mathbf{R}^m\rightarrow dom f$ is called the proximal operator of $f$ \footnote{ Proposition 16.34 in \cite{combettes1}},
i.e.\begin{equation}\label{proximal}
Prox_{f} :  x \mapsto  \arg\min_{u\in dom f } f(u)+\frac{1}{2} ||u-x ||_2^2.
\end{equation}
Define the indicator function of set $\Omega$ as
$\iota_{\Omega}(x)=\begin{cases}
                   0,         & x\in \Omega; \\
                   \infty,    & x\notin \Omega.
                   \end{cases} $
For a closed convex set $\Omega$, $\iota_{\Omega}$ is a proper l.s.c. function.
$\partial \iota_{\Omega}$ is just the normal cone operator of set $\Omega$, that is $\partial \iota_{\Omega}(x) = N_{\Omega}(x)$ \footnote{Example 16.12 of \cite{combettes1}} and
 \begin{equation}\label{normalcone}
   N_{\Omega}(x)
    =\begin{cases}
     \emptyset  &  x\notin \Omega\\
       \{v| \langle v, y-x\rangle\leq 0, \forall y\in \Omega\} & x\in bd(\Omega) \\
         \mathbf{0}    & x\in int(\Omega)
      \end{cases}
 \end{equation}
In this case, we also have \footnote{Example 23.4 of \cite{combettes1}}
\begin{equation}\label{projection_and_resolvent}
Prox_{\iota_{\Omega}}(x)=R_{N_{\Omega}}(x)=P_{\Omega}(x).
\end{equation}

For a single-valued operator $T:\mathbf{R}^m\rightarrow \mathbf{R}^m$, a point $x\in \mathbf{R}^m$ is a fixed point of $T$ if $Tx=x$, and the set of fixed points of $T$ is denoted as $FixT$. The composition of operators $\mathfrak{A}$ and $\mathfrak{B}$, denoted by $ \mathfrak{A}\circ \mathfrak{B}$, is defined via its graph $gra \mathfrak{A}\circ \mathfrak{B}=\{(x,z)| \exists y \in ran \mathfrak{B}, (x,y)\in gra \mathfrak{B}, (y,z)\in gra \mathfrak{A} \}$. We also use $\mathfrak{A}\mathfrak{B}$ to denote the composition $ \mathfrak{A}\circ \mathfrak{B}$ when they are single-valued. Similarly, their sum $\mathfrak{A}+\mathfrak{B}$ is defined as $gra (\mathfrak{A}+\mathfrak{B})=\{(x,y+z)| (x,y)\in gra \mathfrak{A}, (x,z)\in gra \mathfrak{B}\}$.
Suppose operators $\mathfrak{A}$ and $\mathfrak{B}$ are maximally monotone and
$0\in int(dom \mathfrak{A}-dom \mathfrak{B})$, then $\mathfrak{A}+\mathfrak{B}$ is also maximally monotone\footnote{Corollary 24.4 in \cite{combettes1}}. Further suppose that $\mathfrak{A}$ is single-valued, then
$zer(\mathfrak{A}+\mathfrak{B})=Fix R_{\mathfrak{B}}\circ({\rm Id}-\mathfrak{A})$\footnote{Proposition 25.1 in \cite{combettes1}}, which helps to formulate the basic forward-backward operator splitting algorithm  for finding zeros of a sum of monotone operators.

\subsection{Graph theory}

The following concepts are reviewed from \cite{god}.
The information sharing or exchanging among the agents is described by  graph
$\mathcal{G}=(\mathcal{N},\mathcal{E})$.
$\mathcal{N}=
\{1,\cdots,N\}$ is the set of agents, and the edge set $\mathcal{E} \subset \mathcal{N}\times \mathcal{N} $ contains all the  information interactions.
If agent $i$ can get  information from agent $j$, then $(j,i) \in \mathcal{E}$ and
agent $j$ belongs to agent $i$'s  neighbor set $\mathcal{N}_i=\{ j | (j,i) \in
\mathcal{E}\}$, and  $i \notin \mathcal{N}_i$.
$\mathcal{G}$ is said to be undirected when
$(i,j)\in \mathcal{E}$ if and only if $(j,i)\in \mathcal{E}$.
A path of graph $\mathcal{G}$ is a sequence of distinct agents in
$\mathcal{N}$ such that any consecutive agents in the sequence
correspond to an edge of graph $\mathcal{G}$. Agent $j$ is said to be
connected to agent $i$ if there is a path from $j$ to $i$.
$\mathcal{G}$ is said to be  connected if any two agents are
connected.

Define the weighted adjacency matrix $W=[w_{ij}]\in \mathbf{R}^{N\times N}$ of $\mathcal{G}$ with $w_{ij} >0$ if
$j\in \mathcal{N}_i$ and $w_{ij}=0$ otherwise. Assume
$W=W^T$ for undirected graphs.
Define the weighted degree matrix
$Deg= diag\{d_1,\cdots,d_N\}=diag\{ \sum_{j=1}^N w_{1j},..., \sum_{j=1}^N w_{Nj}\}.$
Then the weighted Laplacian of graph $\mathcal{G}$ is
$L=Deg-W.$
When graph $\mathcal{G}$ is a  connected and undirected graph, 0 is a simple
eigenvalue of Laplacian $L$ with the eigenspace $\{ \alpha \mathbf{1}_N| \alpha\in \mathbf{R}\}$,  and
$L \mathbf{1}_N=\mathbf{0}_N$, $\mathbf{1}^T_{N} L=\mathbf{0}^T_N$, while all other eigenvalues are positive.
Denote the eigenvalues of $L$   in an ascending order as $0<s_2\leq \cdots \leq s_N$, then by Courant-Fischer Theorem,
\begin{equation}\label{lap}
\min_{x\neq 0,\atop{\mathbf{1}^Tx=0}} x^T L x =s_2||x||_2^2, \quad \max_{x\neq 0} x^T L x = s_N||x ||^2_2.
\end{equation}
Denote $d^*=\max \{d_1,\cdots,d_N\}$ as the maximal weighted degree of  graph $\mathcal{G}$, then we have the following estimation,
\begin{equation}\label{degree_of_graph_lap}
 d^* \leq s_N \leq 2d^*.
\end{equation}

\section{Problem formulation and distributed algorithm}\label{sec_game_and_algorithm}

\subsection{Game formulation}
Consider a group of agents (players) $\mathcal{N}=\{1,\cdots, N \}$  that seek the generalized Nash equilibrium (GNE) of a  noncooperative game with coupling constraints defined as follows.
Each player $i\in \mathcal{N}$ controls its local decision (strategy or action) $x_i\in \mathbf{R}^{n_i}$.
 Denote $\mathbf{x} =col(x_1,\cdots,x_N) \in \mathbf{R}^n$ as the
decision profile,  i.e., the stacked vector of all the agents' decisions where $\sum_{i=1}^N n_i=n$.
 Denote $\mathbf{x}_{-i}=col(x_1,\cdots,x_{i-1},x_{i+1},\cdots,x_{N})$
as the decision profile stacked vector of the agents' decisions except player $i$.
Agent $i$ aims to optimize its local objective function within its feasible decision set.
The local objective function for agent $i$ is $f_i(x_i,\mathbf{x}_{-i}): \mathbf{R}^n \rightarrow \mathbf{R}$.
Notice that the local objective function of agent $i$ is coupled with  other players' decisions
(however, may not be explicitly  coupled with all other players' decisions).
Moreover, the feasible decision set of  player $i$ also depends on the decisions of the other players with $X_i(\mathbf{x}_{-i}): \mathbf{R}^{n-n_i} \rightarrow 2^{\mathbf{R}^{n_i}}$ denoting a set-valued map that maps  $\mathbf{x}_{-i}$ to the feasible decision set of agent $i$. The aim of agent $i$ is to find the best-response strategy set  given the other players' decision $\mathbf{x}_{-i}$,
\begin{equation}\label{GM}
\min_{x_i\in \mathbf{R}^{n_i}} \;  f_i(x_i,\mathbf{x}_{-i})   \;  \quad s.t., \;\; x_i\in  X_i(\mathbf{x}_{-i}).
\end{equation}
The GNE $\mathbf{x}^*=col(x_1^*,\cdots,x_N^*)$ of the game  in \eqref{GM} is obtained  at the intersection of all the players' best-response sets, and  is defined as:
\begin{equation}
x_i^* \in \arg\min f_i(x_i,\mathbf{x}^*_{-i}),  \; s.t. \; x_i\in  X_i(\mathbf{x}^*_{-i}), \qquad \forall i\in \mathcal{N}.
\end{equation}

Here we consider the  GNE seeking in  noncooperative games where the couplings between  players' feasible sets are specified by  globally shared affine constraints.
Denote $$ X \subset \mathbf{R}^m := \prod_{i=1}^N \Omega_i \bigcap \{ \mathbf{x} \in \mathbf{R}^n | A\mathbf{x} \geq b\},$$
where  $\Omega_i \subset \mathbf{R}^{n_i}$ is a private feasible decision set of player $i$, and
 $A=[A_1,\cdots,A_N]\in \mathbf{R}^{m\times n}$ with $ A_i \in \mathbf{R}^{m\times n_i}$, and $b\in \mathbf{R}^m$.
 Denote $\Omega=\prod_{i=1}^N \Omega_i$.
Given the globally shared set $X$ (which may not be known by any agents),  the following set-valued map gives  the feasible decision set map of agent $i$:
$X_i(\mathbf{x}_{-i}):=\{ x_i\in \mathbf{R}^{n_i}: (x_i,\mathbf{x}_{-i}) \in X \}$, or in other words:
\begin{equation}\nonumber
X_i(\mathbf{x}_{-i}):=\{ x_i\in \Omega_i | A_i x_i \geq b- \sum_{j \neq i,j\in \mathcal{N}} A_j x_j  \}.
\end{equation}

Hence, each agent has a local feasible constraint $x_i \in \Omega_i$, and there exists a coupling constraint shared by all  agents with
sub-matrix $A_i$ characterizing how  agent $i$ is involved in the coupling constraint (shares the global resource).
Notice that  agent $i$ may only know its local $A_i$, in which case the globally shared affine constraint couples the agents' feasible decision sets, but is not known by any agents.

\begin{rem}
We consider  affine coupling constraints for various reasons. Even though not as general as the
nonlinear constraints considered in \cite{pavel1} and \cite{zhuminghui}, this setup does enjoy quite strong modeling
 flexibility. As pointed out in page 191 of \cite{faccinei1}, `` However, it should be noted that the jointly convex assumption on the constraints $\cdots$
practically is likely to be satisfied only when the joint constraints $g_{\mu}=g,\mu=1,...,N$ are
linear, i.e. of the form $Ax \leq b$ for some suitable matrix $A$ and vector $b$."  In fact,
many existing generalized Nash game models adopt affine coupling constraints, as well documented in
\cite{shanbhag3} and \cite{dreves}.
\end{rem}
\begin{assum}\label{assum1}
For player $i$, $f_i(x_i,\mathbf{x}_{-i})$ is a differentiable convex function with respect to $x_i$ given any fixed $\mathbf{x}_{-i}$, and $\Omega_i$ is  a closed  convex set.  $X$ has  nonempty interior point (Slater's condition).
\end{assum}

Suppose $\mathbf{x}^*$ is a GNE of game \eqref{GM}, then for agent $i$, $x_i^*$ is the optimal solution to
the following convex optimization problem:
\begin{equation}\label{op_1}
\min_{x_i} f_i(x_i,\mathbf{x}^*_{-i}), \quad s.t.\; x_i\in \Omega_i, A_ix_i \geq b-\sum_{j \neq i,j\in \mathcal{N}} A_j x^*_j.
\end{equation}
Define a local Lagrangian function for agent $i$ with multiplier $\lambda_i\in \mathbf{R}^m_{+}$ as
\begin{equation}\label{local_lagrangian}
L_i(x_i,\lambda_i;\mathbf{x}_{-i}) =  f_i(x_i,\mathbf{x}_{-i})+\lambda^T_{i}(b-A\mathbf{x}).
\end{equation}
When  $x_i^*$ is an optimal solution to \eqref{op_1}, there exists $\lambda_i^*\in \mathbf{R}_{+}^m$ such that
the following  optimality conditions (KKT) are satisfied:
\begin{equation}
\begin{array}{lll}
\mathbf{0} \in  \nabla_{x_i} L_i(x^*_i,\lambda^*_i;\mathbf{x}^*_{-i}) + N_{\Omega_i}(x_i^*), x^*_i \in \Omega_i \\
\langle \lambda^*_i, b-A\mathbf{x}^* \rangle=0, b-A\mathbf{x}^*\leq \mathbf{0}, \lambda_i^*\geq \mathbf{0},
\end{array}
\end{equation}
These can be equivalently written in the following form by using \eqref{local_lagrangian} and the definition of the normal cone operator in \eqref{normalcone}
\begin{equation}
\begin{array}{lll}\label{kkt_1}
\mathbf{0} \in \nabla_{x_i} f_i(x^*_i,\mathbf{x}^*_{-i}) - A_i^T \lambda_i^* +N_{\Omega_i}(x^*_i)  \\
\mathbf{0} \in  (A\mathbf{x}^*-b) + N_{\mathbf{R}^m_{+}}(\lambda_i^*)
\end{array}
\end{equation}
In fact, since $N_{\Omega}(x)=\emptyset$ when $x\notin \Omega$, it must hold that $x^*_i\in \Omega_i$ and $\lambda^*_i\in \mathbf{R}^m_{+}$ when \eqref{kkt_1} is satisfied.
Furthermore,  $N_{\mathbf{R}^m_{+}} = \prod_{j=1}^m N_{\mathbf{R}_{+}}$. If $[\lambda^*_i]_{k}=0$, then   $N_{\mathbf{R}_{+}}([\lambda_i^*]_{k})= -\mathbf{R}_{+},$ and hence $[b-A\mathbf{x}^*]_{k}\leq 0$;
 and if $[\lambda^*_i]_{k}>0$, we have $N_{\mathbf{R}_{+}}([\lambda_i^*]_{k})= 0,$ hence $[b-A\mathbf{x}^*]_{k}= 0$. Therefore, $b-A\mathbf{x}^*\leq \mathbf{0}$ and $\langle \lambda^*_i, b-A\mathbf{x}^*\rangle=0$.
Denote $\bar{\lambda}=col(\lambda_1,\cdots,\lambda_N)$. Therefore, by Theorem 4.6 in \cite{faccinei1} when $(\mathbf{x}^*,\bar{\lambda}^*)$ satisfies KKT \eqref{kkt_1}  for $i=1,...,N$,  $\mathbf{x}^*$ is a GNE of the game in \eqref{GM}

According to the above discussions, given $\mathbf{x}^*$ as a GNE of game in \eqref{GM}, its corresponding Lagrangian multipliers for the globally shared affine coupling constraint may be different for the agents, i.e.,
$\lambda_1^*\neq \lambda_2^*\neq ,\cdots,\neq\lambda_N^*$.
In this work, we aim to seek a GNE with the same Lagrangian multiplier for all the agents, which has a
nice interpretation from the viewpoint of  variational inequality.

Define
\begin{equation}\label{pseudogradient}
F(\mathbf{x})= col(\nabla_{x_1} f_1(x_1,\mathbf{x}_{-1}), \cdots, \nabla_{x_N} f_N(x_N,\mathbf{x}_{-N})),
\end{equation}
which is usually called the {\em pseudo-gradient}.
The variational inequality (VI) approach to find a GNE  of  game \eqref{GM} is to find the solution of the following $VI(F,X)$:
\begin{equation}\label{vi}
Find \; \mathbf{x}^* \in X, \langle F(\mathbf{x}^*), \mathbf{x}-\mathbf{x}^*\rangle \geq 0, \forall \mathbf{x}\in X.
\end{equation}
Let us check the KKT condition for $VI(F,X)$ in \eqref{vi}. In fact,  $\mathbf{x}^*$ is a solution to $VI(F,X)$ in \eqref{vi} if and only if $\mathbf{x}^*$ is the optimal solution to the following optimization problem:
\begin{equation}\label{vi_op}
\min_{\mathbf{x} \in \mathbf{R}^n} \langle F(\mathbf{x}^*), \mathbf{x}\rangle, \quad s.t., \quad \mathbf{x}\in \Omega,  A\mathbf{x} \geq b
\end{equation}
According to the optimization formulation of $VI(F,X)$ in \eqref{vi_op},  if $\mathbf{x}^*$ solves $VI(F,X)$, there  exists $\lambda^*\in \mathbf{R}^m$ such that the following optimality conditions (KKT) are satisfied:
\begin{equation}
\begin{array}{lll}\label{kkt_2}
\mathbf{0} \in \nabla_{x_i} f_i(x^*_i,\mathbf{x}^*_{-i}) - A_i^T \lambda^* +N_{\Omega_i}(x^*_i),i=1,\cdots,N  \\
\mathbf{0} \in  (A\mathbf{x}^*-b) + N_{\mathbf{R}^m_{+}}(\lambda^*)
\end{array}
\end{equation}

By comparing the two sets of KKT conditions in \eqref{kkt_1} and \eqref{kkt_2} we obtain,
\begin{thm}\footnote{Theorem 2.1 of \cite{facchinei3}}\label{thm_ve}
Suppose Assumption \ref{assum1} holds. Every solution $\mathbf{x}^*$ of  $VI(F,X)$ in \eqref{vi} is a GNE of game in \eqref{GM}.
Furthermore, if $\mathbf{x}^*$ together with $\lambda^*$ satisfies
the KKT conditions for the $VI(F,X)$, i.e.,  \eqref{kkt_2},  then $\mathbf{x}^*$ together with $\lambda_1^*=,\cdots,=\lambda^*_N=\lambda^*$ satisfies
the KKT conditions for the GNE, i.e., \eqref{kkt_1}.
\end{thm}

The solution $\mathbf{x}^*$ of $VI(F,X)$ in \eqref{vi} is termed as a {\it variational GNE} or normalized equilibrium of the  game with coupling constraints in \eqref{GM}. A variational GNE enjoys an economical interpretations of no price discrimination and  has better stability and sensitivity properties,
therefore, can be regarded as a refinement of GNE  (refer to \cite{shanbhag2} for more discussions).
This paper will propose a novel distributed  algorithm  for the agents to
 find a solution of $VI(F,X)$ in \eqref{vi}, thus provides a distributed coordination mechanism  such that  the  coupling constraint is met and  a variational GNE of the game is found.

Define two operators $\mathfrak{A}$ and $\mathfrak{B}$, both from $\Omega\times \mathbf{R}^{m}_{+}$ to $\mathbf{R}^n \times \mathbf{R}^m$ as follows,
\begin{equation}\label{operator_tilde_A_B}
\begin{array}{l}
\mathfrak{A}:\; \left(
                  \begin{array}{c}
                    \mathbf{x} \\
                    \lambda\\
                  \end{array}
                \right)
 \mapsto \left(
           \begin{array}{c}
             F(\mathbf{x}) \\
             -b \\
           \end{array}
         \right)\\
\mathfrak{B}:\;\left(
                  \begin{array}{c}
                    \mathbf{x} \\
                    \lambda\\
                  \end{array}
                \right)
 \mapsto
              \left(
                \begin{array}{c}
                  N_{\Omega}(\mathbf{x}) \\
                  N_{\mathbf{R}^m_{+}}(\lambda)) \\
                \end{array}
              \right)+
              \left(
                \begin{array}{cc}
                  \mathbf{0} & -A^T \\
                  A & \mathbf{0} \\
                \end{array}
              \right)\left(
                  \begin{array}{c}
                    \mathbf{x} \\
                    \lambda\\
                  \end{array}
                \right)

\end{array}
\end{equation}

 By the definition of $F(\mathbf{x})$ in \eqref{pseudogradient}, the KKT conditions \eqref{kkt_2} can be equivalently written as $\mathbf{0}\in (\mathfrak{A} + \mathfrak{B})col(\mathbf{x}^*,\lambda^*)$. Notice that $\mathfrak{B}$ is a maximally monotone operator (similar arguments for this can be found in Lemma \ref{lem_monotone}). Hence,  if $F(\mathbf{x})$ has some additional properties, then solving $VI(F,X)$ can be converted to the problem of finding zeros of a sum of monotone operators.

\begin{assum}\label{assum2}
$F(\mathbf{x})$ defined in \eqref{pseudogradient} is strongly monotone with parameter $\eta$ over $\Omega$: $\langle F(\mathbf{x})-F(\mathbf{y}), \mathbf{x}-\mathbf{y}\rangle \geq \eta ||\mathbf{x}-\mathbf{y}  ||_2^2, \forall \mathbf{x}, \mathbf{y}\in \Omega$, and
$\theta-$ Lipschitz continuous over $\Omega$: $||  F(\mathbf{x})-F(\mathbf{y})||_2 \leq \theta ||\mathbf{x}-\mathbf{y}  ||_2,\forall \mathbf{x}, \mathbf{y}\in \Omega$.
 \end{assum}

\begin{rem}
Assumption \ref{assum2} has also been adopted in \cite{sayed}, \cite{lygeros2} and \cite{zhuminghui}.
Assumption \ref{assum2} guarantees that there exists a unique solution to $VI(F,X)$ in \eqref{vi} \footnote{Theorem 2.3.3 of \cite{pang2}}, thus
guarantees the existence of a GNE for game in \eqref{GM}. However,  the GNE of \eqref{GM} may not be unique. The algorithm for computing {\rm all} the GNE of noncooperative games with coupling constraints is still an opening research topic, and interested readers can refer to \cite{pseng}.
This work aims to provide a distributed algorithm for iterative computation of  a variational GNE of the considered game, which  enjoys nice economical interpretations and stability properties.
\end{rem}

\subsection{Distributed algorithm}

% distributed data structure
 In practice,  each player only knows its private information in game \eqref{GM}, especially when the players interact over large scale networks. It is quite natural that each player can only know its local objective function $f_i(x_i, \mathbf{x}_{-i})$ and local feasible set  $\Omega_i$, which  cannot be shared with other players, because those data contain its local private information, such as cost function, preference and action ability.  Moreover, matrix $A_i$ specifies how  player $i$ participates in the resource allocation or market behavior, hence also contains private information, and $b$ can be decomposed as $b=\sum_{i=1}^N b_i$.
Thus matrix $A_i$ can be regarded as a map from decision space $\mathbf{R}^{n_i}$ to resource space $\mathbf{R}^m$, while vector $b_i$ can be regarded as a local contribution or observation for the global resource. For example, if there are total $m$ markets, and each player $i$ produces a kind of product with amount of $x_i\in \mathbf{R}_{+}$, then $A_i \in \mathbf{R}^m_{+} $ satisfying $\mathbf{1}_m^TA_i=1 $ and $\mathbf{0}\leq A_i \leq \mathbf{1}$ just specifies how each player allocates its production to each market. In this case, if  $\tilde{b}_i \in \mathbf{R}^m_{+}$ is a local observation of market capacity vector, then the true market capacities can be taken as $b=\sum_{i=1}^N b_i=\sum_{i=1}^N\frac{1}{N}\tilde{b}_i$.

Therefore, we assume that  player $i$ only {\it knows} its local $f_i(x_i, \mathbf{x}_{-i})$, $\Omega_i$, and matrix $A_i$ and $b_i$ with $A=[A_1,\cdots,A_N]$ and $b=\sum_{i=1}^N b_i$.  In other words,  player $i$ has a local {\it first-order oracle} of $f_i(x_i, \mathbf{x}_{-i})$ which returns  $\nabla_{x_i} f_i(x_i,\mathbf{x}_{-i})$ given $(x_i,\mathbf{x}_{-i})$,
meanwhile player $i$ can manipulate $\Omega_i$, $A_i$ and $b_i$ for its local computation.

\vskip 2mm
% information observation structure
The players need to find the solution to $VI(F,X)$ in \eqref{vi} in a distributed manner by  local information observation and
sharing, hence find a variational GNE of the game in \eqref{GM} without any coordinator.
To facilitate the local coordination between agents, here we specify two graphs,  $\mathcal{G}_f$ and $\mathcal{G}_{\lambda}$, related to the local information observations or exchanging between players.

Graph $\mathcal{G}_f$, termed as
{\it interference graph}, is defined according to the dependence relationships between the agents' objective functions and the other players' decisions, which is also called {\em graphical model} for games in computer science (refer to \cite{graphmodel}).
For graph $\mathcal{G}_f=(\mathcal{N}, \mathcal{E}_f)$,  $(j,i)\in \mathcal{E}_f$ if the objective function of agent $i$, $f_{i}(x_i,\mathbf{x}_{-i})$ {\bf explicitly} depends on the decision of player $j$.
We define $\mathcal{N}^f_i=\{j | (j,i)\in \mathcal{E}_f\}$  as the set of {\it interference neighbors} whose decisions directly influence the objective function of player $i$. Therefore, the objective function of player $i$ can also be written as $f_i(x_i,\{x_j\}_{j\in \mathcal{N}^f_i})$, and the local oracle of player $i$  returns $\nabla_{x_i} f_i(x_i,\{x_j\}_{j\in \mathcal{N}^f_i})$, i.e., $\nabla_{x_i} f_i(x_i,\mathbf{x}_{-i})$,  given $\{x_j\}_{j\in \mathcal{N}^f_i}$. The local oracle might compute $\nabla_{x_i} f_i(x_i,\{x_j\}_{j\in \mathcal{N}^f_i})$ by approximating  $\nabla_{x_i} f_i(x_i,\mathbf{x}_{-i})$ with local objective function value observations (taking the simultaneous perturbation stochastic approximation in \cite{spsa} as an example), or by utilizing the estimation techniques developed in \cite{pavel4}.

On the other hand, for the coordination of the feasibility of action sets and the consensus of local  multipliers $\lambda_1^*=,\cdots,=\lambda^*_N=\lambda^*$ in Theorem \ref{thm_ve}, we also assume that the agents can exchange certain local information through a {\it multiplier graph} $\mathcal{G}_{\lambda}=(\mathcal{N},\mathcal{E}_{\lambda})$. $(j,i)\in \mathcal{E}_{\lambda}$ if player $i$ can receive certain information from player $j$, while the information to be shared through $\mathcal{G}_{\lambda}$ will be specified later.
Thereby, player $i$ has its {\it multiplier neighbors} $\mathcal{N}^{\lambda}_i=\{j| (j,i)\in \mathcal{E}_{\lambda}\}$. $W=[w_{ij}]$ is the weighted adjacency matrix associated with multiplier graph $\mathcal{G}_{\lambda}$, and $L$ is the corresponding weighted Laplacian matrix.

\begin{assum}\label{assum3}
$\mathcal{G}_{\lambda}=\{\mathcal{N},\mathcal{E}_{\lambda}\}$ is undirected and connected. $W=W^T$.
\end{assum}
\begin{rem}
 We assume that each agent can observe the decisions  which its local objective function directly depends on through interference graph $\mathcal{G}_{f}$. Therefore, player $i$
 can get its local gradient $\nabla_{x_i} f_i(x_i,\{x_j\}_{j\in \mathcal{N}^f_i})$. This type of local information
 observation model has also been adopted in \cite{sayed} and \cite{zhuminghui}.
 On the other hand,  player $i$'s feasible decision set may depend on any other player $k$'s decision even if $f_{i}(x_i,\mathbf{x}_{-i})$ does not explicitly depend on the decision of player $k$, i.e. $k\notin \mathcal{N}_{i}^f$.
  In fact, player $i$'s feasible decision set {\bf implicitly} depends on  all other players' decisions  through the globally shared affine coupling constraint: $A\mathbf{x}\geq b$.
To ensure that the globally shared coupling constraint is satisfied and all the agents have the same local  multipliers, all players must coordinate which necessarily requires that multiplier graph $\mathcal{G}_{\lambda}$ must be connected. Therefore,  $\mathcal{G}_f$ and $\mathcal{G}_{\lambda}$ could be two different information observation or information sharing graphs because they serve for different purposes.
\end{rem}

We are ready to present the main distributed algorithm after the introduction  of algorithm notations.
Agent (player) $i$ controls its local decision variable $x_i\in \mathbf{R}^{n_i}$ and local Lagragian multiplier $\lambda_i \in \mathbf{R}^m$. Meanwhile, we also assume that  player $i$ has a local auxiliary variable $z_i\in \mathbf{R}^m$ for the  coordinations needed to  satisfy the  affine coupling constraint and to reach  consensus of the  local multipliers $\lambda_i$.
As indicated before, player $i$ can compute $\nabla_{x_i}f_i(x_i,\mathbf{x}_{-i})$ by observing the adversary players' decisions that its local objective function  $f_i(x_i,\mathbf{x}_{-i})$ directly depends on, that is the decisions of its interference neighbors in $\mathcal{N}^f_i$.
On the other hand,  player $i$ can also share information related to the local multiplier $\lambda_i$ and local auxiliary variable $z_i$ with its multiplier neighbours in $\mathcal{N}^{\lambda}_i$ through multiplier network $\mathcal{G}_{\lambda}$.

Next we show the distributed algorithm for agent $i$.
\begin{alg}\label{dal_1}
\quad \\
\noindent\rule{0.49\textwidth}{0.7mm}
\begin{align}
&x_{i,k+1}   =  P_{\Omega_i}\big[x_{i,k}-\tau_i ( \nabla_{x_i} f_i(x_{i,k},\{x_{j,k}\}_{j\in \mathcal{N}^f_i}) -A_i^T \lambda_{i,k})\big] \label{al_1}\\
&z_{i,k+1}   =  z_{i,k}                + \nu_i  \sum_{j\in \mathcal{N}^{\lambda}_i} w_{ij}(\lambda_{i,k}-\lambda_{j,k}) \label{al_2}\\
&\lambda_{i,k+1} = P_{\mathbf{R}^{m}_{+}}\Big\{\lambda_{i,k} - \sigma_i\big[ 2 A_i x_{i,k+1}-A_ix_{i,k}-b_i + \sum_{j\in\mathcal{N}^{\lambda}_i} w_{ij}(\lambda_{i,k}-\lambda_{j,k})  \nonumber\\
& \quad \qquad +2\sum_{j\in \mathcal{N}^{\lambda}_i} w_{ij}(z_{i,k+1}-z_{j,k+1})-\sum_{j\in \mathcal{N}^{\lambda}_i} w_{ij}(z_{i,k}-z_{j,k})  \big]\Big\}  \label{al_3}
\end{align}
\noindent\rule{0.49\textwidth}{0.7mm}
$\tau_i,\nu_i, \sigma_i >0$ are fixed constant step-sizes of player $i$, and $W=[w_{ij}]$ is the weighted adjacency matrix of $\mathcal{G}_{\lambda}$.
\end{alg}

Algorithm \ref{dal_1} runs sequentially as follows. At the iteration time $k$,
player $i$  gets $ \nabla_{x_i}f_i(x_{i,k}, \mathbf{x}_{-i,k})$ by observing $x_{j,k},j\in \mathcal{N}^f_i$ through interference graph $\mathcal{G}_f$,  and  updates $x_{i,k}$  with \eqref{al_1}; meanwhile, player $i$ gets $\lambda_{j,k}, j\in \mathcal{N}^{\lambda}_i $ through multiplier graph $\mathcal{G}_{\lambda}$, and updates $z_{i,k}$ by \eqref{al_2}. Then player $i$ gets $z_{j,k+1}, j\in \mathcal{N}^{\lambda}_i $ through multiplier graph $\mathcal{G}_{\lambda}$ and updates $\lambda_{i,k}$ with  \eqref{al_3} that also employs the most recent local information $x_{i,k+1}$.

Intuitively speaking, \eqref{al_1} employs the projected gradient descent of the local Lagrangian function in \eqref{local_lagrangian}. \eqref{al_2} can be regarded as the discrete-time intergration for the consensual errors of local copies of multipliers, which will ensure the consensus of $\lambda_i$ eventually. In fact, a similar dynamics has been employed in distributed optimization in \cite{leijinlong}.  Finally, \eqref{al_3} updates local multiplier by a  combination of the projected gradient ascent of local Lagrangian function \eqref{local_lagrangian} and a proportional-integral dynamics for consensual errors of multipliers. Section \ref{sec_algorithm_develpo} will give a detailed algorithm development from the viewpoint of operator splitting methods.

 Algorithm \ref{dal_1} is  a totally distributed algorithm and has following key features:

 i). The full data decomposition and privacy protection is achieved   since each player only needs to know its local objective function $f_i(x_i,\mathbf{x}_{-i})$ and local feasible set $\Omega_i$.

 ii). The matrix $A$ is decomposed and each block $A_i$ is kept by player $i$,  hence the privacy of each player is protected because $A_i$ describes how player $i$ is involved in the market or competition.

 iii). Each player only needs to observe the decisions  which its local objective function directly depends on, and only needs to share information with its multiplier neighbours, through $\mathcal{G}_{f}$ and $\mathcal{G}_{\lambda}$, respectively. Both graphs usually have sparse edge connections,  therefore, the observation or communication burden is relieved. Furthermore, the information observation related with decisions and  the information sharing related with multipliers is decoupled and accomplished with different graphs $\mathcal{G}_f$ and $\mathcal{G}_{\lambda}$, respectively. Therefore, those two information sharing processes can work in a parallel manner, and can be designed independently.

 iv). The algorithm converges with fixed step-sizes under some mild conditions, and  works in a  Gauss-Seidel manner that utilizes the most recent information when updating $\lambda_i$. Moreover, \eqref{al_1} and \eqref{al_2}  can even be computed  in parallel for  player $i$.

\section{Algorithm development}\label{sec_algorithm_develpo}

In this section, we first show how Algorithm \ref{dal_1} is developed and provide the motivations behind the algorithm's convergence analysis. Then we  verify that the limiting  point of Algorithm \ref{dal_1} solves the $VI(F,X)$ in \eqref{vi}, and thus finds a variational GNE of the game in \eqref{GM}.

Algorithm \ref{dal_1} is inspired by the forward-backward splitting methods for finding zeros
of the sum of monotone operators (\cite{combettes1}) and the primal-dual algorithm for optimization with linear composition terms by \cite{condat}. The key difference are the specific operator splitting form and the augmentation of variables to achieve distributed computations. Next, we  systematically show how to write Algorithm \ref{dal_1} in the form of a  forward-backward operator splitting algorithm.

Let us define some notations  to write Algorithm \ref{dal_1} in a compact form.
 Denote $\mathbf{x}_{k}=col(x_{1,k},\cdots,x_{N,k})\in \mathbf{R}^{n}$,
 $\bar{\lambda}_k=col(\lambda_{1,k},\cdots,\lambda_{N,k})\in \mathbf{R}^{mN}$,
$\bar{z}_k=col(z_{1,k},\cdots,z_{N,k})\in \mathbf{R}^{mN}$,  $\bar{b}=col(b_1,\cdots,b_N)\in \mathbf{R}^{mN}$,  $\bar{\tau}=\diag\{\tau_1 I_{n_1},\cdots,\tau_N I_{n_N}\}\in \mathbf{R}^{n\times n} $, $\bar{\nu}=\diag\{\nu_1 I_{m},\cdots,\nu_N I_{m}\}\in \mathbf{R}^{mN\times mN}$, $\bar{\sigma}=\diag\{\sigma_1 I_m,\cdots,\sigma_N I_m\}\in \mathbf{R}^{mN\times mN}$, $\Lambda= \diag\{A_1,\cdots,A_N\}\in \mathbf{R}^{mN\times n}$ and  $\Lambda^T=\diag\{A_1^T,\cdots,A_N^T\}\in \mathbf{R}^{n\times mN}$.
$\bar{L}=L\otimes I_m$ where $L\in \mathbf{R}^{N\times N} $ is the weighted Laplacian matrix of multiplier graph $\mathcal{G}_{\lambda}$.

Using these notations, the definition of pseudo-gradient $F(\mathbf{x})$ in \eqref{pseudogradient} and
$P_{\prod_{i=1}^N \Omega_i}(col(x_1,\cdots,x_N))=  col(P_{\Omega_1}(x_1),\cdots,P_{\Omega_N}(x_N))$\footnote{Proposition 23.16 of \cite{combettes1}},
 Algorithm \ref{dal_1} can be written in a compact form as:
\begin{alg}\label{dal_2}
\quad \\
\noindent\rule{0.49\textwidth}{0.7mm}
\begin{align}
&\mathbf{x}_{k+1}   =  P_{\Omega}\big[\mathbf{x}_{k}-\bar{\tau} ( F(\mathbf{x}_{k}) -\Lambda^T \bar{\lambda}_{k})\big] \label{cal_1}\\
&\bar{z}_{k+1}   =  \bar{z}_{k}                + \bar{\nu}  \bar{L} \bar{\lambda}_{k} \label{cal_2}\\
&\bar{\lambda}_{k+1} = P_{\mathbf{R}^{mN}_{+}}\Big \{\bar{\lambda}_{k} - \bar{\sigma}\big[  \Lambda (2\mathbf{x}_{k+1}-\mathbf{x}_{k})-\bar{b} + \bar{L} \bar{\lambda}_{k}  +\bar{L} (2\bar{z}_{k+1}-\bar{z}_k)\big]\Big\} \label{cal_3}
\end{align}
\noindent\rule{0.49\textwidth}{0.7mm}
\end{alg}

Using the fact that $P_{\Omega}(x)=Prox_{\iota_{\Omega}}(x)=R_{N_{\Omega}}(x)$ in \eqref{projection_and_resolvent} and the definition of resolvent operator as
$R_{N_{\Omega}}(x)=({\rm Id}+N_{\Omega})^{-1}$,
equation \eqref{cal_1} be be written as
$ \mathbf{x}_{k+1}= ({\rm Id}+N_{\Omega})^{-1} [\mathbf{x}_{k}-\bar{\tau} ( F(\mathbf{x}_{k}) -\Lambda^T \bar{\lambda}_{k}) ]$, or equivalently,
\begin{equation}\label{eq_c_1}
\mathbf{x}_{k}-\bar{\tau} ( F(\mathbf{x}_{k}) -\Lambda^T \bar{\lambda}_{k}) \in \mathbf{x}_{k+1}+N_{\Omega}(\mathbf{x}_{k+1}).
\end{equation}
Notice that $\bar{\tau}^{-1} = \diag\{\frac{1}{\tau_1}I_{n_1},\cdots,\frac{1}{\tau_N}I_{n_N}\}$,
$\bar{\nu}^{-1} = \diag\{\frac{1}{\nu_1}I_{m},\cdots,\frac{1}{\nu_N}I_{m}\}$ and
$\bar{\sigma}^{-1} = \diag\{\frac{1}{\sigma_1}I_{m},\cdots,\frac{1}{\sigma_N}I_{m}\}$.
Furthermore, $N_{\Omega}$ is a  cone and $N_{\Omega}(\mathbf{x})=\prod_{i=1}^{N}N_{\Omega_i}(x_i)$, hence
$\bar{\tau}^{-1} N_{\Omega}(\mathbf{x})=N_{\Omega}(\mathbf{x})$.
Therefore, \eqref{eq_c_1} can be written as
\begin{equation}\label{eq_c_11}
- F(\mathbf{x}_{k}) \in  N_{\Omega}(\mathbf{x}_{k+1})-\Lambda^T \bar{\lambda}_{k+1}
+ \bar{\tau}^{-1}(\mathbf{x}_{k+1}-\mathbf{x}_{k})
    +\Lambda^T (\bar{\lambda}_{k+1}- \bar{\lambda}_{k}).
\end{equation}

 Moreover, $N_{\mathbf{R}^{mN}_{+}}$ is a cone,  $N_{\mathbf{R}^{mN}_{+}}(\bar{\lambda})=\prod_{i=1}^N N_{\mathbf{R}^m_{+}}(\lambda_i)$, and $\bar{\sigma}^{-1}N_{\mathbf{R}^{mN}_{+}}=N_{\mathbf{R}^{mN}_{+}}$.  Then with similar arguments, equation \eqref{cal_3} can be written as:
\begin{equation}\label{eq_c_2}
\bar{\lambda}_{k} - \bar{\sigma}\Big[  \Lambda (2\mathbf{x}_{k+1}-\mathbf{x}_{k})-\bar{b} + \bar{L} \bar{\lambda}_{k}+\bar{L}(2\bar{z}_{k+1}-\bar{z}_k)\Big]
\in \bar{\lambda}_{k+1}+N_{\mathbf{R}^{mN}_{+}}(\bar{\lambda}_{k+1})
\end{equation}
or equivalently,
\begin{equation}\label{eq_c_22}
\begin{array}{lll}
 -[\bar{L} \bar{\lambda}_{k} -\bar{b}]& \in & N_{\mathbf{R}^{mN}_{+}}(\bar{\lambda}_{k+1})+\Lambda \mathbf{x}_{k+1}+\bar{L}\bar{z}_{k+1} \\
 &+&\Lambda (\mathbf{x}_{k+1}-\mathbf{x}_{k})+\bar{L}(\bar{z}_{k+1}-\bar{z}_k)+\bar{\sigma}^{-1}(\bar{\lambda}_{k+1}-\bar{\lambda}_{k})
\end{array}
\end{equation}

Therefore, equation\eqref{cal_2} together with \eqref{eq_c_11} and \eqref{eq_c_22} can be written in a compact form as:
\begin{align}\label{expanded_al}
&-\left(
  \begin{array}{c}
     F(\mathbf{x}_{k}) \\
    \mathbf{0} \\
    \bar{L} \bar{\lambda_k}-\bar{b} \\
  \end{array}
\right)
\in
\left(
  \begin{array}{c}
    N_{\Omega}(\mathbf{x}_{k+1}) -\Lambda^T \bar{\lambda}_{k+1} \\
      -\bar{L}\bar{\lambda}_{k+1}  \\
     N_{\mathbf{R}^{mN}_{+}}(\bar{\lambda}_{k+1})+\Lambda\mathbf{x}_{k+1}+\bar{L}\bar{z}_{k+1} \\
  \end{array}
\right)
\nonumber\\
&\qquad+
\left(
  \begin{array}{ccc}
    \bar{\tau}^{-1} & \mathbf{0} & \Lambda^T \\
    \mathbf{0} & \bar{\nu}^{-1} & \bar{L} \\
    \Lambda& \bar{L} &  \bar{\sigma}^{-1} \\
  \end{array}
\right)
\left(
  \begin{array}{c}
    \mathbf{x}_{k+1}-\mathbf{x}_{k} \\
    \bar{z}_{k+1}-\bar{z}_{k} \\
    \bar{\lambda}_{k+1}-\bar{\lambda}_{k} \\
  \end{array}
\right)
\end{align}

Denote
\begin{equation}\label{metric_matrix}
\Phi=\left(
  \begin{array}{ccc}
    \bar{\tau}^{-1} & \mathbf{0} & \Lambda^T \\
    \mathbf{0} & \bar{\nu}^{-1} & \bar{L}  \\
    \Lambda& \bar{L}&  \bar{\sigma}^{-1} \\
  \end{array}
\right)
\end{equation}
 Notice that matrix $\Phi \in \mathbf{R}^{(n+2mN) \times (n+2mN)}$ is symmetric due to $L=L^T$ and $\bar{L}=L\otimes I_m$.

\vskip 3mm
Denote $\bar{\Omega}=\Omega \times \mathbf{R}^{mN} \times \mathbf{R}_{+}^{mN}$. Define the  operators
$\bar{\mathfrak{A}}: \bar{\Omega}\rightarrow \mathbf{R}^{n+2mN}$ and $\bar{\mathfrak{B}}: \bar{\Omega}\rightarrow 2^{\mathbf{R}^{n+2mN}}$ as follows,
\begin{align}\label{operator_A_B}
&\mathfrak{A}:=\left(\begin{array}{c}
\mathbf{x}\\ \bar{z} \\ \bar{\lambda} \\ \end{array} \right) \mapsto
\left(
  \begin{array}{c}
    F(\mathbf{x}) \\
    \mathbf{0} \\
    \bar{L}  \bar{\lambda}-\bar{b} \\
  \end{array}
\right),  \quad \nonumber \\
&
\mathfrak{B}:=\left(\begin{array}{c}
\mathbf{x}\\ \bar{z} \\ \bar{\lambda} \\ \end{array} \right)
\mapsto
\left(
  \begin{array}{c}
    N_{\Omega}(\mathbf{x}) \\
    \mathbf{0} \\
    N_{\mathbf{R}^{mN}_{+}}(\bar{\lambda}) \\
  \end{array}
\right)
+
\left(
  \begin{array}{ccc}
    \mathbf{0} & \mathbf{0} & -\Lambda^T \\
    \mathbf{0} & \mathbf{0} & -\bar{L} \\
    \Lambda & \bar{L} & \mathbf{0} \\
  \end{array}
\right)
\left(\begin{array}{c}
\mathbf{x}\\ \bar{z} \\ \bar{\lambda} \\ \end{array} \right)
\end{align}

\begin{rem}
Operators $\bar{\mathfrak{A}}$ and $\bar{\mathfrak{B}}$ in \eqref{operator_A_B}
can be regarded as an extension  of
operators $\mathfrak{A}$ and $\mathfrak{B}$ in \eqref{operator_tilde_A_B} by augmenting $\lambda \in \mathbf{R}^m$ to $\bar{\lambda}\in \mathbf{R}^{mN}$ and introducing auxiliary variables $\bar{z}\in \mathbf{R}^{mN}$.
 Moreover, operators in \eqref{operator_A_B} utilize  $L\mathbf{1}_N = \mathbf{0}_N$ of  Laplacian matrix $L$ to ensure the consensus of $\lambda_i$, and utilize  $\mathbf{1}^T_N L = \mathbf{0}_N^T$ to ensure the feasibility of affine coupling constraints.
\end{rem}

The next result shows that Algorithm \ref{dal_1}, or equivalently Algorithm \ref{dal_2}, can be regarded as a forward-backward operator splitting method for finding zeros of a sum of operators,  or an iterative computation of fixed points of a composition of operators.
\begin{lem}\label{lem_fix_approximation}
Suppose $\Phi$ in \eqref{metric_matrix} is positive definite, and operators $\bar{\mathfrak{A}}$ and $\bar{\mathfrak{B}}$ in \eqref{operator_A_B} are maximally monotone. Denote $T_1:={\rm Id}-\Phi^{-1}\bar{\mathfrak{A}} $ and $T_2:=({\rm Id}+\Phi^{-1}\bar{\mathfrak{B}})^{-1}$. Then any limiting point of Algorithm \ref{dal_1}, i.e., $col(\mathbf{x}^*, \bar{z}^*,\bar{\lambda}^*)$,  is a zero of $\bar{\mathfrak{A}}+\bar{\mathfrak{B}}$ and is a fixed point of $T_2\circ T_1$.
\end{lem}
{\bf Proof:}
Denote $\varpi= col(\mathbf{x}, \bar{z},\bar{\lambda})$,
then using \eqref{expanded_al}, \eqref{metric_matrix} and \eqref{operator_A_B},  Algorithm \ref{dal_1} can written in a compact form as follows:
\begin{equation}\label{compact_operator_1}
-\bar{\mathfrak{A}}(\varpi_{k}) \in \bar{\mathfrak{B}}(\varpi_{k+1})+\Phi(\varpi_{k+1}-\varpi_{k}).
\end{equation}
Since  $\Phi$ is symmetric and positive definite, we can write equation \eqref{compact_operator_1}
as $\varpi_{k}-\Phi^{-1}\bar{\mathfrak{A}}(\varpi_{k})\in \varpi_{k+1}+\Phi^{-1}\bar{\mathfrak{B}}(\varpi_{k+1})$, or equivalently,
\begin{equation}\label{compact_operator_2}
({\rm Id}-\Phi^{-1}\bar{\mathfrak{A}})(\varpi_{k}) \in ({\rm Id}+\Phi^{-1}\bar{\mathfrak{B}})(\varpi_{k+1}).
\end{equation}
Since  $\Phi^{-1}\bar{\mathfrak{B}}$ is maximally monotone (refer to Lemma \ref{lem_monotone_metric}), $({\rm Id}+\Phi^{-1}\bar{\mathfrak{B}})^{-1}$ is single-valued. Then by the definition of the inverse of a set-valued operator, Algorithm \ref{dal_1} is written as
\begin{equation}\label{com_fix}
\varpi_{k+1}=({\rm Id}+\Phi^{-1}  \bar{\mathfrak{B}})^{-1}({\rm Id}- \Phi^{-1} \bar{\mathfrak{A}}) \varpi_{k}: = T_2 \circ T_1 \varpi_{k}.
\end{equation}

Suppose that Algorithm \ref{dal_1}, or equivalently \eqref{com_fix},  converges to a   limiting point $\varpi^*$. Then by the  continuity of the right hand of Algorithm \ref{dal_1} (In fact, the right hand of Algorithm \ref{dal_1} is Lipschitz continuous due to Assumption \ref{assum2} and the nonexpansive property of projection operator), $\varpi^*=T_2T_1\varpi^*$. Therefore, any limiting point
of Algorithm \ref{dal_1} is a fixed point of the composition $T_2 \circ T_1$, and Algorithm \ref{dal_1} can be regarded as an iterative computation of fixed points of  $T_2 \circ T_1$.

By Theorem 25.8 of \cite{combettes1}, \eqref{com_fix} is  also the forward-backward splitting algorithm for finding zeros of a sum of  monotone operators, hence $\varpi^* \in zer(\Phi^{-1}\bar{\mathfrak{A}}+\Phi^{-1} \bar{\mathfrak{B}}) $ for any  limiting point $\varpi^*$.
Since $\Phi$ is positive definite, any limiting point $\varpi^*=col(\mathbf{x}^*,\bar{z}^*, \bar{\lambda}^*)$, i.e., $\varpi^*=T_2T_1 \varpi^*$,  also belongs to $zer(\bar{\mathfrak{A}}+\bar{\mathfrak{B}})$. In fact,
\begin{equation}
\begin{array}{ll}\nonumber
&\quad\varpi^*=({\rm Id}+\Phi^{-1}\bar{\mathfrak{B}})^{-1} ({\rm Id}- \Phi^{-1}\bar{ \mathfrak{A}}) \varpi^* \\
&\Leftrightarrow ({\rm Id}- \Phi^{-1} \bar{\mathfrak{A}}) \varpi^* \in ({\rm Id}+\Phi^{-1}  \bar{\mathfrak{B }}) \varpi^* \\
& \Leftrightarrow \mathbf{0}\in \Phi^{-1} (\bar{\mathfrak{A}}+\bar{\mathfrak{B }})(  \varpi^*)\Leftrightarrow
\mathbf{0}\in (\bar{\mathfrak{A}}+\bar{\mathfrak{B }})(  \varpi^*).
\end{array}
\end{equation}
\hfill $\Box$

\begin{rem}
The iteration $\varpi_{k+1}=T_2T_1\varpi_{k}$ is also known as {\rm Picard iteration} for iteratively approximating fixed points of $T_2T_1$ (refer to \cite{iterative_fix}).
Lemma \ref{lem_metric_matrix} will give a sufficient condition for $\Phi$ to be positive definite.
Lemma \ref{lem_monotone}  will give the condition that ensures $\bar{\mathfrak{A}}$ and $\bar{\mathfrak{B}}$ to be maximally monotone.
\end{rem}

\vskip 2mm
The next result shows that any limiting point of Algorithm \ref{dal_1},  that is,   any zero point of operator $\bar{\mathfrak{A}}+ \bar{\mathfrak{B}}$,  is a variational GNE of game \eqref{GM}.

\begin{thm}\label{thm_zero_is_correct}
Suppose that Assumptions \ref{assum1}-\ref{assum3} hold. Consider operators $\bar{\mathfrak{A}}$ and $\bar{\mathfrak{B}}$  defined in \eqref{operator_A_B}, and operators $\mathfrak{A}$ and $\mathfrak{B}$ defined in \eqref{operator_tilde_A_B}. Then the following statements hold:

(i): Given any $col(\mathbf{x}^*,\bar{z}^*,\bar{\lambda}^*) \in zer(\bar{\mathfrak{A}}+\bar{\mathfrak{B}})$,
then  $\mathbf{x}^*$ solves the  $VI(F,X)$ in \eqref{vi}, hence $\mathbf{x}^*$ is a variational GNE of game in \eqref{GM}. Moreover $\bar{\lambda}^*=\mathbf{1}_N\otimes \lambda^*$, and the  multiplier $\lambda^*$ together with $\mathbf{x}^*$  satisfy the  KKT condition in \eqref{kkt_2}, i.e.,  $col(\mathbf{x}^*,\lambda^*)\in zer(\mathfrak{A}+\mathfrak{B})$.

(ii): $zer(\mathfrak{A}+\mathfrak{B}) \neq \emptyset$ and  $zer(\bar{{\mathfrak{A}}}+\bar{{\mathfrak{B}}}) \neq \emptyset$.
\end{thm}

{\bf Proof:}
(i): By the definition of operators $\bar{\mathfrak{A}}$, $\bar{\mathfrak{B}}$ in \eqref{operator_A_B}, we have,
\begin{equation}\label{operator_sum}
\begin{array}{ll}
&\bar{\mathfrak{A}}+\bar{\mathfrak{B}}:  \;col(\mathbf{x},\bar{z},  \bar{\lambda})\mapsto\\
&\left(
                                 \begin{array}{ccc}
                                  \mathbf{0}                           & \mathbf{0}& -\Lambda^T  \\
                                    \mathbf{0}   &    \mathbf{0}          &-\bar{L}\\
                                  \Lambda        &    \bar{L}         &   \bar{L}   \\
                                 \end{array}
                               \right) \left(
                                         \begin{array}{c}
                                           \mathbf{x} \\
                                         \bar{z} \\
                                           \bar{\lambda}\\
                                         \end{array}
                                       \right)
                                 + \left(
                                                            \begin{array}{c}
                                                              \mathbf{0} \\
                                                              \mathbf{0}\\
                                                              -\bar{b}\\
                                                            \end{array}
                                                          \right)
                                                          +
                                        \left(
                                 \begin{array}{c}
                                   F(\mathbf{x})+ N_{\Omega}(\mathbf{x}) \\
                                   \mathbf{0}\\
                                   N_{\mathbf{R}^{mN}_{+}}(\bar{\lambda})  \\
                                 \end{array}
                               \right)
\end{array}
\end{equation}

Suppose that $col(\mathbf{x}^*,\bar{z}^*,\bar{\lambda}^*) \in zer(\bar{\mathfrak{A}}+\bar{\mathfrak{B}})$.
From the second line of \eqref{operator_sum}, $$-\bar{L}\bar{\lambda}^*=-L\otimes I_m \bar{\lambda}^* = \mathbf{0}_{mN}.$$
It follows that $\bar{\lambda}^*=\mathbf{1}_N\otimes \lambda^*, \lambda^*\in \mathbf{R}^m$ since $L$ is the weighted Laplacian of
multiplier graph $\mathcal{G}_{\lambda}$ and $\mathcal{G}_{\lambda}$ is connected due to Assumption \ref{assum3}.

Then by the first line of \eqref{operator_sum}, combined with  $\Lambda^T=\diag\{A_1^T,\cdots,A_N^T\}$ and $\lambda^*_1=\lambda^*_2=,\cdots,=\lambda_N^*=\lambda^*$, we have
\begin{equation}
\mathbf{0}_n  \in - \Lambda^T 1_N\otimes\lambda^* +  F(\mathbf{x}^*)+ N_{\Omega}(\mathbf{x}^*),
\end{equation}
or equivalently,
\begin{equation}\label{fact1}
\mathbf{0}_{n_i}  \in - A_i^T \lambda^* + \nabla_{x_i}f_i(x_i^*,\mathbf{x}^*_{-i}) + N_{\Omega_i}(x_i^*),  \; i=1,\cdots,N.
\end{equation}

By the third line of \eqref{operator_sum} and using $\bar{L}\bar{\lambda}^*=\mathbf{0}$, it follows that
$$ \mathbf{0}_{mN} \in \Lambda \mathbf{x}^* - \bar{b} + \bar{L} \bar{z}^*+ N_{\mathbf{R}^{mN}_{+}}(\bar{\lambda}^*).$$
This implies that  there exist $v_1,v_2,\cdots,v_N \in N_{R^{m}_{+}}(\lambda^*)$, such that
\begin{align}
 \mathbf{0}_{mN} &=  \Lambda \mathbf{x}^* - \bar{b} + L\otimes I_m \bar{z}^*+ col(v_1,\cdots,v_N). \nonumber
\end{align}
Multiplying both sides of above equation with $\mathbf{1}^T_N\otimes I_m$ and combining with  $\mathbf{1}^T L=\mathbf{0}^T$, we have
\begin{align}
\mathbf{0}_{m} &= (\mathbf{1}^T_N \otimes I_m) ( \Lambda \mathbf{x}^* - \bar{b} + L\otimes I_m \bar{z}^* + col(v_1,\cdots,v_N)) \nonumber  \\
   \;          &= \sum_{i=1}^N A_ix^*_i - \sum_{i=1}^N b_i+\sum_{i=1}^N v_i \nonumber
\end{align}
Due to the fact that $N_{\bigcap_{i=1}^N \Omega_i}=\sum_{i=1}^N N_{\Omega_i}$ if $\bigcap_{i=1}^N int(\Omega_i)\neq \emptyset$ {\footnote{Corollary 16.39 of \cite{combettes1}}, we have
\begin{equation}
\begin{array}{lll}\label{fact2}
\mathbf{0}_{m}  & \in  & \sum_{i=1}^N A_ix^*_i - \sum_{i=1}^N b_i+\sum_{i=1}^N N_{R^{m}_{+}}(\lambda^*)\\
\;              & \in  & \sum_{i=1}^N A_ix^*_i - \sum_{i=1}^N b_i+N_{ \bigcap_{i=1}^N R^{m}_{+}}(\lambda^*) \\
\;              & \in  & \sum_{i=1}^N A_ix^*_i - \sum_{i=1}^N b_i + N_{R^{m}_{+}}(\lambda^*).
\end{array}
\end{equation}

By \eqref{fact1} and \eqref{fact2},  for any $col(\mathbf{x}^*,\bar{\lambda}^*,\bar{z}^*)\in zer(\bar{\mathfrak{A}}+\bar{\mathfrak{B}})$, the KKT condition for $VI(F,X)$ in \eqref{vi}, i.e. \eqref{kkt_2},  is satisfied for $\mathbf{x}^*$,  $\lambda^*$.
 We conclude that $\mathbf{x}^*$ solves $VI(F,X)$ in \eqref{vi}, and is a variational GNE of game \eqref{GM} by Theorem \ref{thm_ve}. It also follows that $\lambda^*$ together with $\mathbf{x}^*$  satisfy the  KKT condition in \eqref{kkt_2}.
 This also implies $col(\mathbf{x}^*,\lambda^*) \in zer(\mathfrak{\mathfrak{A}}+\mathfrak{\mathfrak{B}})$.

\vskip 1.5mm
(ii) Under Assumptions \ref{assum1} and \ref{assum2}, the considered game in \eqref{GM} has a unique variational GNE $\mathbf{x}^*$, and there exists $\lambda^* \in \mathbf{R}^m$ such that the KKT condition \eqref{kkt_2} is satisfied, i.e. $col(\mathbf{x}^*,\lambda^*) \in zer(\mathfrak{A}+\mathfrak{B})$. Therefore $zer(\mathfrak{A}+\mathfrak{B}) \neq \emptyset$.

Then we need to show that there exists $col(\mathbf{x}^*, \bar{\lambda}^*, \bar{z}^*)$ such that $ col(\mathbf{x}^*, \bar{\lambda}^*, \bar{z}^*)\in zer(\bar{{\mathfrak{A}}}+\bar{{\mathfrak{B}}})$.

Take $\bar{\lambda}^*=\mathbf{1}_N \otimes \lambda^*$. Because $L\mathbf{1}_N=\mathbf{0}$, the second line of \eqref{operator_sum}
is satisfied.

Since $col(\mathbf{x}^*,\lambda^*) \in zer(\mathfrak{A}+\mathfrak{B})$,
$\mathbf{0} \in F(x^*)-A^T\lambda^* + N_{\Omega}(\mathbf{x}^*)$. Using $\lambda^*_1=\lambda^*_2=,\cdots,=\lambda_N^*=\lambda^*$, $\Lambda^T \bar{\lambda}^*= col(A_1^T \lambda^*, \cdots, A_N^T \lambda^*)=A^T \lambda^*$. Therefore, the first line of \eqref{operator_sum} is satisfied with $\mathbf{x}^*, \mathbf{1}_N \otimes \lambda^*$.

Moreover, with $col(\mathbf{x}^*,\lambda^*) \in zer(\mathfrak{A}+\mathfrak{B})$,
 $\mathbf{0}_{m} \in   A\mathbf{x}^* -b + N_{\mathbf{R}^m_{+}}(\lambda^*)$.
Then we need to show that there exists $\bar{z}^*=col(z_1^*,\cdots, z_N^*)\in \mathbf{R}^{mN}$, such that the third line of \eqref{operator_sum} is satisfied.
Take $v^*\in N_{\mathbf{R}^m_{+}}(\lambda^*)$ such that $ \mathbf{0} =  A\mathbf{x}^*  -b +v^*$.
Since $\bar{\lambda}^*=\mathbf{1}_N \otimes \lambda^*$ and
$N_{\mathbf{R}^{mN}_{+}}(\bar{\lambda}^*)=\prod_{i=1}^N N_{\mathbf{R}^m_{+}}(\lambda^*)$, take $v_1^*=v_2^*=,\cdots,=v_N^*=\frac{1}{N}v^* \in N_{\mathbf{R}^m_{+}}(\lambda^*)$.
Then $(\mathbf{1}^T_N\otimes I_m) ( \Lambda \mathbf{x}^* -\bar{b} +  \bar{L}\bar{\lambda}^* + col(v_1^*,\cdots,v^*_N)  )=\sum_{i=1}^N A_i x_i^* - \sum_{i=1}^N b_i +  v^*=A\mathbf{x}^* -b + v^*=\mathbf{0}_m$.
That is $\Lambda \mathbf{x}^* -\bar{b} +  \bar{L}\bar{\lambda}^* + col(v_1^*,\cdots,v^*_N)\in Null(\mathbf{1}^T_N\otimes I_m)$.
By the  fundamental theorem of linear algebra\footnote{Page 405 of  \cite{mayer_matrix_analsysis}}, $Null(\mathbf{1}^T_N\otimes I_m)=Range(\mathbf{1}_N\otimes I_m)^{\bot}$ and $Range(\bar{L})=Null(\bar{L})^{\bot}$ since $\bar{L}$ is  also symmetric.
Notice that $Null(\bar{L})=Range(\mathbf{1}_N\otimes I_m)$, hence $\Lambda \mathbf{x}^* -\bar{b} +  \bar{L}\bar{\lambda}^* + col(v_1^*,\cdots,v^*_N)\in  Range(\bar{L})$.
Noticing that $col(v_1^*,\cdots,v^*_N)\in N_{\mathbf{R}^{mN}_{+}}(\bar{\lambda}^*)$,  there exists $\bar{z}^*\in \mathbf{R}^{mN}$ such that the third line of \eqref{operator_sum} is satisfied with $\mathbf{x}^*, \bar{z}^*, \mathbf{1}_N \otimes \lambda^*$.
Therefore, $zer(\bar{{\mathfrak{A}}}+\bar{{\mathfrak{B}}}) \neq \emptyset$.
\hfill $\Box$

\section{Convergence Analysis}\label{sec_alg_converge}

In this section, we prove the convergence of Algorithm \ref{dal_1} by giving a sufficient step-size choice condition.  The analysis is based on the compact reformulation \eqref{com_fix}.
 We will first show that all the prerequisites in Lemma \ref{lem_fix_approximation} can be satisfied under suitable step-sizes.
Then \eqref{com_fix}, i.e., Algorithm \ref{dal_1},  can be regarded as a forward-backward splitting algorithm for finding zeros of
a sum of monotone operators,  or equivalently, an iterative computation of fixed points of a composition of operators.

In fact, some existing NE (GNE)  algorithms can also be regarded as a type of iterative computation of fixed points of operators, such as the  best-response learning dynamics (\cite{lygeros1}), relaxation algorithms based on Nikaido-Isoda function (\cite{krawczyk2} and \cite{krawczyk1}) and the proximal-best response algorithm in \cite{pang}. Most of above works built their algorithm convergence analysis on the {\em contractive property} of underlying operators.
However, the contractivity assumption on operators is usually quite restrictive.
 Herein we  resort to the theory of  {\em averaged operators} and {\em firmly nonexpansive operators} for convergence analysis.
Firstly we give some basic definitions and properties of averaged operators and firmly nonexpansive operators\footnote{Chapter 4 and Chapter 20 of \cite{combettes1}}.  All the following results are valid in Hilbert spaces, thus they hold in Euclidean spaces with any $G-$matrix induced norm $|| \cdot||_{G}$, given $G$ as a symmetric positive definite matrix.
Denote by $|| \cdot||$ an arbitrary  matrix induced norm in a finite dimensional Euclidean space.

An operator $T: \Omega \subset \mathbf{R}^m\rightarrow \mathbf{R}^m$  is nonexpansive if
it is $1-$Lipschitzian, i.e., $||T(x)-T(y) || \leq ||x-y||, \forall x,y \in \Omega.$
An operator $T$ is $\alpha-$averaged  if there exists a nonexpansive operator $T^{'}$ such that  $T=(1-\alpha){\rm Id}+\alpha T^{'}$.
Denote the class of $\alpha-$averaged operators as $\mathcal{A}(\alpha)$.
If $T\in \mathcal{A}(\frac{1}{2})$, then $T$ is also called firmly nonexpansive operator.
\begin{lem}\footnote{Proposition 4.25 of \cite{combettes1}}\label{lem_alpha_opertor}
Given an operator $T:\Omega \subset  \mathbf{R}^m \rightarrow \mathbf{R}^m $ and $\alpha\in (0,1)$, the following three statements are equivalent:

(i): $T\in \mathcal{A}(\alpha) $;

(ii): $||Tx-Ty ||^2 \leq ||x-y ||^2 -\frac{1-\alpha}{\alpha} || (x-y)-(Tx-Ty) ||^2, \forall x, y \in \Omega$

(iii): $||Tx-Ty ||^2+ (1-2\alpha)||x-y||^2 \leq 2(1-\alpha) \langle x-y, Tx-Ty \rangle, \forall x, y \in \Omega$.
\end{lem}
By (iii) of Lemma \eqref{lem_alpha_opertor}, $T\in \mathcal{A}(\frac{1}{2})$ if and only if
\begin{equation}\label{firmly_nonexpansively}
||T(x)-T(y) ||^2 \leq \langle x-y, T(x)-T(y) \rangle, \forall x, y \in \Omega.
\end{equation}
The operator $T$ is called $\beta-$cocoercive (or $\beta-$inverse strongly monotone) if  $\beta T$ is firmly nonexpansive, i.e.,
\begin{equation}\label{cocoercive}
\beta||T(x)-T(y) ||^2 \leq \langle x-y, T(x)-T(y) \rangle, \forall x, y \in \Omega.
\end{equation}

\begin{lem} \footnote{Theorem 18.15 of
\cite{combettes1}}\label{gradientcocoerceive}
For a convex differentiable function $f$ with $\vartheta-$ Lipschitzian gradient, we have $\nabla f $ to be $\frac{1}{\vartheta}-$cocoercive, i.e.,
\begin{equation}
\frac{1}{\vartheta} ||\nabla f(x)-\nabla f(y) ||_2^2 \leq \langle x-y, \nabla f(x)- \nabla f(y)\rangle.
\end{equation}
\end{lem}
Lemma \ref{gradientcocoerceive} is known as {\it Baillon-Haddad theorem}, and one elementary proof can  be found in  Theorem 2.1.5 of \cite{nesterov}.

\begin{lem}\footnote{Proposition 23.7 of \cite{combettes1}}\label{lem_maxiamlly_monotone}
If operator $\Delta$ is maximally monotone, then $T=R_\Delta=({\rm Id}+\Delta)^{-1}$ is firmly nonexpasive and $2R_{\Delta}-{\rm Id}$ is nonexpansive.
\end{lem}
Hence, the projection operator onto a closed convex set  is firmly nonexpansive since
$P_{\Omega}=Prox_{\iota_{\Omega}}=R_{\partial \iota_{\Omega}}=R_{N_{\Omega}}$  and $N_{\Omega}$ is maximally monotone\footnote{Example 20.41 and Proposition 4.8 of \cite{combettes1}}.

\vskip 2mm
In the following, we  analyze the maximal  monotonicity of operators $\bar{\mathfrak{A}}$, $\bar{\mathfrak{B}}$ in \eqref{operator_A_B}, the positive definite property of matrix $\Phi$, and the properties of operators  $T_1$ and $T_2$ defined in Lemma \ref{lem_fix_approximation} by giving sufficient step-sizes choice conditions,  which are shown in Lemma \ref{lem_monotone}, \ref{lem_metric_matrix} and \ref{lem_monotone_metric}, respectively.

\begin{lem}\label{lem_monotone}
Suppose Assumptions \ref{assum1}-\ref{assum3} hold. Given an Euclidean space with norm $|| \cdot||_2$, then

(i): Operator $\bar{\mathfrak{A}}$  in \eqref{operator_A_B} is $\beta-$cocoercive with  $0<\beta \leq \min\{\frac{1}{2d^*}, \frac{\eta}{ \theta^2}\}$ where $d^*$ is the maximal weighted degree of multiplier graph $\mathcal{G}_{\lambda}$, i.e.,
 $d^*=\max\{ \sum_{j=1}^N w_{1j},\cdots, \sum_{j=1}^Nw_{Nj}\}$, and $\eta,\theta$ are parameters in Assumption \ref{assum2};

(ii): Operator $\bar{\mathfrak{B}}$ in \eqref{operator_A_B} is maximally monotone.
\end{lem}
{\bf Proof:}
(i): According to the definition of $ \bar{\mathfrak{A}}$ in \eqref{operator_A_B} and the definition of $\beta-$cocercive in \eqref{cocoercive},
 we need to prove that
 \begin{equation}\label{equ_lem5_3_1}
 \begin{array}{ll}
 &\langle \left(
     \begin{array}{c}
        F(\mathbf{x}_1) \\
        0\\
        \bar{L} \bar{\lambda}_1 -\bar{b} \\
     \end{array}
   \right)-
   \left(
     \begin{array}{c}
        F(\mathbf{x}_2) \\
        0\\
        \bar{L} \bar{\lambda}_2 -\bar{b} \\
     \end{array}
   \right)
   ,   \left( \begin{array}{c}
                                           \mathbf{x}_1 \\
                                          \bar{z}_1\\
                                           \bar{\lambda}_1 \\
                                         \end{array}
                                       \right)
                                          -
                                \left( \begin{array}{c}
                                           \mathbf{x}_2 \\
                                          \bar{z}_2\\
                                           \bar{\lambda}_2 \\
                                         \end{array}
                                       \right)
                                          \rangle \geq \\
                                          &\beta ||\left(
     \begin{array}{c}
        F(\mathbf{x}_1) \\
        0\\
        \bar{L} \bar{\lambda}_1 -\bar{b} \\
     \end{array}
   \right)-
   \left(
     \begin{array}{c}
        F(\mathbf{x}_2) \\
        0\\
        \bar{L} \bar{\lambda}_2 -\bar{b} \\
     \end{array}
   \right)||_2^2, \forall \left( \begin{array}{c}
                                           \mathbf{x}_1 \\
                                          \bar{z}_1\\
                                           \bar{\lambda}_1 \\
                                         \end{array}
                                       \right),
                                \left( \begin{array}{c}
                                           \mathbf{x}_2 \\
                                          \bar{z}_2\\
                                           \bar{\lambda}_2 \\
                                         \end{array}
                                       \right) \in \bar{\Omega}
   \end{array}
 \end{equation}

Notice that $\bar{L} \bar{\lambda} -\bar{b}$ is the gradient of function $\tilde{f}(\bar{\lambda}):=\frac{1}{2} \bar{\lambda}^T \bar{L} \bar{\lambda} -\bar{b}^T\bar{\lambda} $. Moveover, $\tilde{f}(\bar{\lambda})$ is a convex function since $\nabla^2 \tilde{f}(\bar{\lambda})=\bar{L}$ is positive semi-definite due to Assumption \ref{assum3} \footnote{Proposition 17.10 of \cite{combettes1}}.
It can easily be verified that  $\bar{L} \bar{\lambda} -\bar{b}$ is $||L||_2-$Lipschitz continuous (notice that the eigenvalues of $\bar{L}$ are just the elements in  $col(0,s_2\cdots,s_N)\otimes \mathbf{1}_m$), therefore, by Lemma \ref{gradientcocoerceive}
\begin{equation}\label{equ_lem5_3_2}
\langle \bar{L} \bar{\lambda}_1 -\bar{b} - ( \bar{L} \bar{\lambda}_2 -\bar{b}), \bar{\lambda}_1-\bar{\lambda}_2    \rangle \geq \frac{1}{||L ||_2} || \bar{L} \bar{\lambda}_1  - \bar{L} \bar{\lambda}_2 ||_2^2.
\end{equation}
Since $||L||_2 \leq s_N$ and $ d^* \leq s_N \leq 2d^*$  by \eqref{degree_of_graph_lap} where $d^*=\max\{d_1,\cdots,d_N\}$ is the maximal weighted degree of multiplier graph $\mathcal{G}_{\lambda}$,  we have  $\frac{1}{||L ||_2} \geq \frac{1}{2d^*}$.

Meanwhile by Assumption \ref{assum2}, $F(\mathbf{x})$ is $\eta-$strongly monotone and $\theta-$ Lipschitz continuous over $\Omega$. By $ || F(\mathbf{x}_1)- F(\mathbf{x}_2)||_2^2 \leq \theta^2 ||\mathbf{x}_1- \mathbf{x}_2  ||_2^2$,
$\forall\; \mathbf{x}_1, \mathbf{x}_2\in \Omega$ we have
\begin{equation}\label{equ_lem5_3_3}
\langle  F(\mathbf{x}_1)- F(\mathbf{x}_2),  \mathbf{x}_1- \mathbf{x}_2 \rangle \geq
\eta ||\mathbf{x}_1- \mathbf{x}_2 ||_2^2  \geq  \frac{\eta}{ \theta^2}|| F(\mathbf{x}_1)- F(\mathbf{x}_2)||_2^2.
\end{equation}

Taking  $0<\beta \leq \min\{\frac{1}{2d^*}, \frac{\eta}{ \theta^2}\}$ and adding \eqref{equ_lem5_3_2} and \eqref{equ_lem5_3_3} yields \eqref{equ_lem5_3_1}. Thus  operator $\bar{\mathfrak{A}}$ is $\beta-$cocoercive.
By the definition of $\beta-$cocoercive in \eqref{cocoercive},  operator $\bar{\mathfrak{A}}$ is also monotone.
Since operator $\bar{\mathfrak{A}}$ is also single-valued, it is also maximally monotone.

(ii): The operator $\bar{\mathfrak{B}}$ in \eqref{operator_A_B} can be written as:
 　\begin{equation}\nonumber
   \left(
                                 \begin{array}{ccc}
                                  \mathbf{0}     &\mathbf{0} & -\Lambda^T \\
                                  \mathbf{0}   &\mathbf{0} &-\bar{L}\\
                                  \Lambda        &\bar{L}    &\mathbf{0}    \\
                                 \end{array}
                               \right) \left(
                                         \begin{array}{c}
                                           \mathbf{x} \\
                                         \bar{z}\\
                                           \bar{\lambda} \\
                                         \end{array}
                                       \right)
                                        + \left(
                                 \begin{array}{c}
                                   N_{\Omega}(\mathbf{x}) \\
                                   \mathbf{0}\\
                                   N_{\mathbf{R}^{mN}_{+}} \\
                                 \end{array}
                               \right):= \mathfrak{B}_1 + \mathfrak{B}_2.
   \end{equation}
Since $\bar{L}$ is symmetric,  $\mathfrak{B}_1$ is a skew-symmetric matrix, i.e., $\mathfrak{B}^T_1=-\mathfrak{B}_1$.
 Hence, $\mathfrak{B}_1$ is maximally monotone\footnote{Example 20.30 of \cite{combettes1}}.
 %we have
% \begin{equation}
% \begin{array}{ll}
% & \langle \left(
%              \begin{array}{c}
%                \mathbf{x}_1-\mathbf{x}_2 \\
%               \bar{z}_1-\bar{z}_2 \\
%               \bar{\lambda}_1-\bar{\lambda}_2 \\
%              \end{array}
%            \right),
%       \left(
%  \begin{array}{c}
%    -   diag\{A_i^T\}_{i=1}^N ( \bar{\lambda}_1-\bar{\lambda}_2) \\
%     -\bar{L} (\bar{\lambda}_1-\bar{\lambda}_2)  \\
%    \Lambda(\mathbf{x}_1-\mathbf{x}_2)  +  \bar{L }( \bar{z}_1-\bar{z}_2) \\
%  \end{array}
%\right) \rangle  \\ \nonumber
%& =   -(\mathbf{x}_1-\mathbf{x}_2)^T \Lambda^T ( \bar{\lambda}_1-\bar{\lambda}_2) \\
%&-( \bar{z}_1-\bar{z}_2)^T\bar{L} (\bar{\lambda}_1-\bar{\lambda}_2)+( \bar{\lambda}_1-\bar{\lambda}_2)^T \Lambda (\mathbf{x}_1-\mathbf{x}_2)\\
%&+ ( \bar{\lambda}_1-\bar{\lambda}_2)^T\bar{L }( \bar{z}_1-\bar{z}_2) \\
%&=0
%\end{array}
%\end{equation}

$\mathfrak{B}_2$ can be written as the direct sum of  $N_{\Omega} \times \mathbf{0}_{mN}\times N_{\mathbf{R}^{mN}_{+}} $.
Both $N_{\Omega}$ and $N_{\mathbf{R}^{mN}_{+}}$ are maximally monotone as normal cones of closed convex sets.
Obviously, $\mathbf{0}_{mN}$ is also maximally monotone as a single-valued operator. Furthermore, the direct sum of maximally monotone operators is also maximally monotone\footnote{Proposition 20.23 of \cite{combettes1}},  hence $\mathfrak{B}_2$ is also maximally monotone.

Obviously, $dom \mathfrak{B}_1 = \mathbf{R}^{n+2mN}$, hence $\bar{\mathfrak{B}}=\mathfrak{B}_1+\mathfrak{B}_2$ is also maximally monotone\footnote{ Corollary 24.4 of \cite{combettes1}}.
\hfill $\Box$

\begin{lem}\label{lem_metric_matrix}
Given  any $\delta >0$, if each player $i$ takes  $\tau_i>0$, $\nu_i>0$ and $\sigma_i>0$ as its local fixed step-sizes in Algorithm \ref{dal_1} that satisfy:
\begin{equation}
\begin{array}{lll}\label{step_size_choice}
\tau_i   & \leq  &   \frac{1}{\max_{j=1,...n_i} \{\sum^m_{k=1} | [A_i^T]_{jk} |\}+\delta},\\
\nu_i    & \leq  &   \frac{1}{2d_i+\delta} \\
\sigma_i & \leq  &   \frac{1}{\max_{j=1,...m} \{ \sum^{n_i}_{k=1} | [A_i]_{jk} | \}+2d_i+ \delta}
\end{array}
\end{equation}
then matrix  $\Phi$ defined in \eqref{metric_matrix} is positive definite, and $  \Phi - \delta I_{n+2mN} $ is  positive semi-definite.
\end{lem}
{\bf Proof:}
It is sufficient to show that  $\Phi-\delta I_{n+2mN}$ is positive semi-definite.
\begin{align}
\Phi-\delta I_{n+2mN}=
\left(
\begin{array}{ccc}
    \bar{\tau}^{-1}-\delta I_n & \mathbf{0} &   \Lambda^T \\
    \mathbf{0} & \bar{\nu}^{-1}-\delta I_{mN}  & \bar{L} \\
    \Lambda & \bar{L} &  \bar{\sigma}^{-1}-\delta I_{mN} \\
  \end{array}
\right),
\end{align}
 where
 $\bar{\tau}^{-1}-\delta I_n=\diag\{(\frac{1}{\tau_1}-\delta)I_{n_1},\cdots,(\frac{1}{\tau_N}-\delta)I_{n_N}\}$,
 $\bar{\nu}^{-1}-\delta I_{mN}=\diag\{(\frac{1}{\nu_1}-\delta)I_{m},\cdots,(\frac{1}{\nu_N}-\delta)I_{m}\}$,
 and
 $\bar{\sigma}^{-1}-\delta I_{mN}=\diag\{(\frac{1}{\sigma_1}-\delta) I_{m},\cdots,(\frac{1}{\sigma_N}-\delta )I_{m}\}$.
 One sufficient condition for  matrix  $\Phi-\delta I_{n+2mN}$ to be positive semi-definite is that  it is diagonally dominant with nonnegative diagonally elements, that is for every row of the matrix  the diagonal entry  is larger than or equal to the sum of the magnitudes of all the other (non-diagonal) entries in that row. This is equivalent to require that,
\begin{equation}
\begin{array}{lll}\label{equ_lem5_4_1}
\frac{1}{\tau_i}- \delta  & \geq  &   \max_{j=1,...n_i} \{\sum^m_{k=1} | [A_i^T]_{jk} |\}  \\
\frac{1}{\nu_i}- \delta   & \geq  &   \sum_{j=1}^m  | L_{ij} | = 2 d_i \\
\frac{1}{\sigma_i}-\delta &  \geq &   \max_{j=1,...m} \{ \sum^{n_i}_{k=1} | [A_i]_{jk} | \} + \sum_{j=1}^m  | L_{ij} | \\
\qquad                    & = &\max_{j=1,...m} \{ \sum^{n_i}_{k=1} | [A_i]_{jk} | \}+2d_i
\end{array}
\end{equation}
It can easily be verified  that if each agent  chooses its local step-sizes satisfying \eqref{step_size_choice}, then \eqref{equ_lem5_4_1} is satisfied.
\hfill $\Box$

Given a globally known parameter  $\delta$, each agent can {\it independently} choose its local step sizes $\tau_i$ ,$\nu_i$, and $\sigma_i$ with the rule given  in \eqref{step_size_choice}.

Suppose that the step-sizes $\tau_i, \mu_i, \sigma_i$ in Algorithm \ref{dal_1} are chosen such that $\Phi$ in \eqref{metric_matrix} is  positive definite.
  Thus we can define a norm induced by matrix $\Phi$ for a finite Euclidean space as $||x||_{\Phi}= \sqrt{\langle x, x\rangle_{\Phi}} = \sqrt{\langle \Phi x, x\rangle}$. The next result investigates the properties of operators $\Phi^{-1}\bar{\mathfrak{A}}, \Phi^{-1}\bar{\mathfrak{B}}$ and $T_1, T_2$ defined in Lemma \ref{lem_fix_approximation} under $\Phi-$induced norm $ || \cdot||_{\Phi}$.

\begin{lem}\label{lem_monotone_metric}
Suppose Assumptions \ref{assum1}-\ref{assum3} hold.
Take $0<\beta \leq \min\{\frac{1}{2d^*}, \frac{\eta}{ \theta^2}\}$ where $d^*$ is the maximal weighted degree of multiplier graph $\mathcal{G}_{\lambda}$, and $\eta,\theta$ are parameters in Assumption \ref{assum2}.
Take $ \delta > \frac{1}{2\beta}$.
Suppose that the step-sizes $\tau_i, \nu_i, \sigma_i$ in Algorithm \ref{dal_1} are chosen to satisfy \eqref{step_size_choice}.
Then the operators $\Phi^{-1}\bar{\mathfrak{A}}$ and $\Phi^{-1}\bar{\mathfrak{B}}$ with $\Phi$ in \eqref{metric_matrix} and $\bar{\mathfrak{A}},\bar{\mathfrak{B}}$ in \eqref{operator_A_B}, and operators $T_1={\rm Id}-\Phi^{-1}\bar{\mathfrak{A}},   T_2=R_{\Phi^{-1} \bar{\mathfrak{B}}}=({\rm Id}+\Phi^{-1}\bar{\mathfrak{B}})^{-1} $ defined as Lemma \ref{lem_fix_approximation}  satisfy the following properties under the $\Phi-$induced norm $||\cdot||_{\Phi}$:

(i).  $\Phi^{-1}\bar{\mathfrak{A}}$ is $\beta\delta-$cocoercive, and  $T_1 \in \mathcal{A}(\frac{1}{2\delta\beta})$.

(ii).  $\Phi^{-1}\bar{\mathfrak{B}}$ is  maximally monotone, and $T_2 \in \mathcal{A}(\frac{1}{2})$.

\end{lem}
{\bf Proof:}
(i): By the definition of cocoercivity in \eqref{cocoercive}, we need to prove $ \langle \Phi^{-1}\bar{\mathfrak{A}}(x)-\Phi^{-1}\bar{\mathfrak{A}}(y),  x-y\rangle_{\Phi} \geq  {\beta}{\delta}  || \Phi^{-1}\bar{\mathfrak{A}}(x)-\Phi^{-1}\bar{\mathfrak{A}}(y) ||^2_{\Phi}$, $\forall x,y \in \bar{\Omega}$.
Noticing that by the choice of parameters $\tau_i, \nu_i, \sigma_i$, we have that matrix $\Phi - \delta I_{n+2mN} $ is positive semi-definite from Lemma \ref{lem_metric_matrix}.
Denote $s_{max}(\Phi)$ and $s_{min}(\Phi)$ as the maximal and minimal eigenvalues of matrix $\Phi$. It must hold that $s_{\max}(\Phi)\geq s_{\min}(\Phi) \geq \delta$.
 Furthermore, $||\Phi||_{2}=s_{max}(\Phi) \geq s_{min}(\Phi)=\frac{1}{ || \Phi^{-1} ||_2} \geq \delta$ \footnote{Proposition 5.2.7 and 5.2.8 in \cite{mayer_matrix_analsysis}}, therefore, we  also have $||\Phi^{-1} ||_2 \leq  \frac{1}{\delta}$.
 Notice that the operator $\bar{\mathfrak{A}}$ is single-valued and $\Phi^{-1}$ is also nonsingular, so that for any $x,y\in \bar{\Omega}$,
\begin{equation}
\begin{array}{ll}\nonumber
&|| \Phi^{-1}\bar{\mathfrak{A}}(x)-\Phi^{-1}\bar{\mathfrak{A}}(y) ||^2_{\Phi} \\
&= \langle \Phi\Phi^{-1}(\bar{\mathfrak{A}}(x)-\bar{\mathfrak{A}}(y)),\Phi^{-1}(\bar{\mathfrak{A}}(x)-\bar{\mathfrak{A}}(y)) \rangle \\
&= || \bar{\mathfrak{A}}(x)-\bar{\mathfrak{A}}(y)||^2_{\Phi^{-1}} \leq \frac{1}{\delta} || \bar{\mathfrak{A}}(x)-\bar{\mathfrak{A}}(y)||_2^2
\end{array}
\end{equation}
By the $\beta-$cocoercive property of $\bar{\mathfrak{A}}$ in Lemma \ref{lem_monotone} and the above inequality,
\begin{equation}
\begin{array}{lll}
&& \langle \Phi^{-1}\bar{\mathfrak{A}}(x)-\Phi^{-1}\bar{\mathfrak{A}}(y),  x-y\rangle_{\Phi} = \langle \bar{\mathfrak{A}}(x)-\bar{\mathfrak{A}}(y),  x-y \rangle \\
&&  \geq \beta  ||\bar{\mathfrak{A}}(x)-\bar{\mathfrak{A}}(y) ||_2^2
\geq  {\beta}{\delta}  || \Phi^{-1}\bar{\mathfrak{A}}(x)-\Phi^{-1}\bar{\mathfrak{A}}(y) ||^2_{\Phi}.
\end{array}
\end{equation}
Therefore, the operator $\Phi^{-1}\bar{\mathfrak{A}}$ is $\beta\delta-$cocoercive under the $\Phi-$induced norm $||\cdot||_{\Phi}$.

Moreover, $\beta\delta\Phi^{-1}\bar{\mathfrak{A}}$ is firmly nonexpansive by the definition of cocoercive operator. This implies that there exists a nonexpansive operator $\breve{T}$ such that $\beta\delta\Phi^{-1}\bar{\mathfrak{A}}=\frac{1}{2} \breve{T}+\frac{1}{2} {\rm Id}$. Then
$$T_1= {\rm Id}-\Phi^{-1}\bar{\mathfrak{A}}=(1-\frac{1}{2\beta\delta}) {\rm Id}+ \frac{1}{2\beta\delta} (-\breve{T}) \in \mathcal{A}(\frac{1}{2\beta\delta}) $$
since $1<2\beta\delta$ by the assumption that $\delta > \frac{1}{2\beta}$ and $-\breve{T}$ is also nonexpansive.

(ii). $\Phi$ is symmetric positive definite and nonsingular.
For any $(x,u)\in gra \Phi^{-1}\bar{\mathfrak{B}}$ and $(y,v)\in gra \Phi^{-1}\bar{\mathfrak{B}}$,
$\Phi u \in \Phi \Phi^{-1}\bar{\mathfrak{B}}(x)\in \bar{\mathfrak{B}}(x) $ and
$\Phi v \in \Phi \Phi^{-1}\bar{\mathfrak{B}}(y)\in \bar{\mathfrak{B}}(y) $.
Then
  $\langle x-y,u-v\rangle_{\Phi}= \langle  x-y, \Phi(u-v) \rangle \geq 0, \forall x,y$ since $\bar{\mathfrak{B}}$ is monotone by Lemma \ref{lem_monotone}.   Therefore, $\Phi^{-1}\bar{\mathfrak{B}}$ is monotone under the $\Phi-$matrix induced product $\langle \cdot, \cdot \rangle_{\Phi}$.

Furthermore, take $(y,v) \in \bar{\Omega} \times \mathbf{R}^{n+2mN} $, and
$ \langle x-y, u-v \rangle_{\Phi} \geq 0$, for any other $(x,u) \in  gra(\Phi^{-1}\bar{\mathfrak{B}})$.
For any $(x, \tilde{u}) \in gra \bar{\mathfrak{B}}$, we have $(x,\Phi^{-1}\tilde{u})\in gra(\Phi^{-1}\bar{\mathfrak{B}})$.
$ \langle  x-y, \Phi( \Phi^{-1}\tilde{u}-v) \rangle  \geq 0$, or equivalently,
$ \langle  x-y, \tilde{u}-\Phi v) \rangle  \geq 0$.
Since $\bar{\mathfrak{B}}$ is maximally monotone, then $(y,\Phi v) \in gra \bar{\mathfrak{B}}$. We conclude that $v\in \Phi^{-1}\bar{\mathfrak{B}} (y)$ which implies that $\Phi^{-1}\bar{\mathfrak{B}}$ is maximally monotone.

Therefore, by Lemma \ref{lem_maxiamlly_monotone} $T_2= ({\rm Id}+\Phi^{-1}\bar{\mathfrak{B}})^{-1}$ is firmly nonexpansive under  the $\Phi-$matrix induced norm $|| \cdot||_{\Phi}$.
\hfill $\Box$

\vskip 3mm
Summarizing the above results,  take $0<\beta \leq \min\{\frac{1}{2d^*}, \frac{\eta}{ \theta^2}\}$, $\delta > \frac{1}{2\beta}$, and $\tau_i, \nu_i, \sigma_i$ satisfying \eqref{step_size_choice}, then
$T_1:= {\rm Id}-\Phi^{-1}\bar{\mathfrak{A}} \in  \mathcal{A}(\frac{1}{2\delta\beta})$,
and $T_2:= ({\rm Id}+\Phi^{-1}\bar{\mathfrak{B}})^{-1}\in \mathcal{A}(\frac{1}{2})$ in the Euclidean space with $\Phi-$matrix induced norm $||\cdot||_{\Phi}$. The next result shows the  convergence of  Algorithm \ref{dal_1} based on its compact reformulation \eqref{com_fix} and the properties of $T_1$ and $T_2$.

\begin{thm}
Suppose Assumptions \ref{assum1}-\ref{assum3} hold.
Take $0<\beta \leq \min\{\frac{1}{2d^*}, \frac{\eta}{ \theta^2}\}$ where $d^*$ is the maximal weighted degree of multiplier graph $\mathcal{G}_{\lambda}$, and $\eta,\theta$ are parameters in Assumption \ref{assum2}.
Take $ \delta > \frac{1}{2\beta}$.
The step-sizes $\tau_i, \nu_i, \sigma_i$ in Algorithm \ref{dal_1} are chosen to satisfy \eqref{step_size_choice}.
Then with Algorithm \ref{dal_1}, each player has its local strategy $x_{i,k}$ converging
to its corresponding component in the variational GNE of game \eqref{GM}, and  the  local Lagrangian multipliers $\lambda_{i,k}$ of all the agents converge to the same Lagrangian multiplier corresponding with KKT condition \eqref{kkt_2}, i.e.,
\begin{equation}\label{con_1}
\lim_{k\rightarrow \infty} x_{i,k}\rightarrow x_i^*,  \quad \lim_{k\rightarrow \infty}  \lambda_{i,k} \rightarrow \lambda^*, \forall i=1,\cdots,N.
\end{equation}
\end{thm}
{\bf Proof:}
With Lemma \ref{lem_monotone} and \ref{lem_metric_matrix}, Algorithm \ref{dal_1} can be written in a compact form \eqref{com_fix} according to Lemma \ref{lem_fix_approximation}, i.e., $\varpi_{k+1}=T_2T_1\varpi_{k}$. The convergence analysis will be conducted via the analysis of this iterative computation of fixed points of $T_2 \circ T_1$.

Firstly, by (i) and (ii) of Lemma \ref{lem_alpha_opertor} and the fact that  $T_1,T_2$ are averaged operators due to Lemma \ref{lem_monotone_metric}, $T_1,T_2$ are  also nonexpansive operators under the $\Phi-$matrix induced norm $ || \cdot||_{\Phi}$.
Take any $\varpi^* \in zer(\bar{\mathfrak{A}}+\bar{\mathfrak{B }})$, or equivalently any fixed point of $T_2 \circ T_1$, i.e., $\varpi^* =T_2 T_1 \varpi^*$, and then by Lemma \ref{lem_fix_approximation} and \eqref{com_fix},
\begin{equation}
\begin{array}{lll}\label{equ_thm_5_2}
&&||\varpi_{k+1}- \varpi^*||_{\Phi} = ||  T_2( T_1 \varpi_{k})-  T_2 (T_1 \varpi^*)||_{\Phi} \\
&&\leq  || T_1 \varpi_{k}-   T_1 \varpi^* || \leq  || \varpi_{k}-  \varpi^* ||_{\Phi}
\end{array}
\end{equation}
 Hence the sequence $\{||\varpi_{k}-\varpi^* ||_{\Phi}\}$ is  non-increasing and bounded from below.
 By the monotonic convergence theorem, $\{|| \varpi_{k}-  \varpi^* || \}$ is bounded  and converges for every $\varpi^* \in zer(\mathfrak{A}+\mathfrak{B })$.

By Lemma \ref{lem_monotone_metric}, $T_1\in \mathcal{A}(\frac{1}{2\delta\beta})$ and $T_2\in \mathcal{A}(\frac{1}{2})$. Denote $\xi = \frac{1}{2\delta \beta} \in (0,1)$.  Then with (ii) of Lemma \ref{lem_alpha_opertor} and \eqref{com_fix} we have,
\begin{equation}
\begin{array}{l}\label{equ_thm_57_1}
||\varpi_{k+1}- \varpi^* ||_{\Phi}^2= ||  T_2 ( T_1 \varpi_{k})-  T_2  (T_1 \varpi^*) ||_{\Phi}^2\\
  \leq  ||T_1 \varpi_{k}- T_1 \varpi^*||_{\Phi}^2 \\
  \quad - || (  T_1 \varpi_{k}-  T_1 \varpi^*) - (T_2 T_1 \varpi_{k}-  T_2  T_1 \varpi^*)　||^2_{\Phi} \\
\leq ||\varpi_{k}-  \varpi^*  ||^2_{\Phi}- || (  T_1 \varpi_{k}-  T_1 \varpi^*) - (T_2 T_1 \varpi_{k}-  T_2  T_1 \varpi^*)　||^2_{\Phi} \\
-\frac{1-\xi}{\xi}|| \varpi_{k}-\varpi^*- (T_1 \varpi_{k}- T_1 \varpi^*) ||^2_{\Phi}
\end{array}
\end{equation}
where the first inequality follows by $T_2\in \mathcal{A}(\frac{1}{2})$ and the second inequality follows by $T_1\in \mathcal{A}(\frac{1}{2\delta\beta})$, both utilizing (ii) of Lemma \ref{lem_alpha_opertor}.
Notice that
\begin{equation}\label{equation_norm}
\alpha ||x||^2 + (1-\alpha)||y ||^2 = ||\alpha x +(1-\alpha)y ||^2 + \alpha(1-\alpha) ||x-y||^2.
\end{equation}
For the second and third terms on the right hand side of \eqref{equ_thm_57_1},
\begin{equation}
\begin{array}{ll}\label{equ_thm_57_2}
& || (  T_1 \varpi_{k}-  T_1 \varpi^*)  - (T_2 T_1 \varpi_{k}-  T_2  T_1 \varpi^*)　||^2_{\Phi} \\
&\qquad+\frac{1-\xi}{\xi}||  \varpi_{k}-\varpi^*- (T_1 \varpi_{k}- T_1 \varpi^*) ||^2_{\Phi} \\
& = \frac{1}{\xi} [\xi|| (  T_1 \varpi_{k}-  T_1 \varpi^*) - (T_2 T_1 \varpi_{k}-  T_2  T_1 \varpi^*)　||^2_{\Phi} \\
&\qquad +(1-\xi)|| (T_1 \varpi_{k}- T_1 \varpi^*)-(\varpi_{k}-\varpi^*)||^2_{\Phi}  ]\\
& =\frac{1}{\xi}|| (T_1 \varpi_{k}- T_1 \varpi^*)-\xi (T_2 T_1 \varpi_{k}-  T_2  T_1 \varpi^*) \\
&\qquad\qquad  - (1-\xi)(\varpi_{k}-\varpi^*)  ||^2_{\Phi}\\
&\qquad+ \frac{1}{\xi}\xi(1-\xi)|| (  T_1 \varpi_{k}-  T_1 \varpi^*) - (T_2 T_1 \varpi_{k}-  T_2  T_1 \varpi^*)\\
&\qquad \qquad -(T_1 \varpi_{k}- T_1 \varpi^*)+(\varpi_{k}-\varpi^*) ||_{\Phi}^2\\
&\geq (1-\xi) || (\varpi_{k}-\varpi^*)-(T_2 T_1 \varpi_{k}-  T_2  T_1 \varpi^*)||^2_{\Phi}\\
&= (1-\xi) ||\varpi_{k}-T_2 T_1 \varpi_{k} ||^2_{\Phi}\\
&= (1-\xi) || \varpi_{k}-\varpi_{k+1}||^2_{\Phi}
\end{array}
\end{equation}
where the second equality follows from \eqref{equation_norm} by setting $\alpha=\xi$, $x= (  T_1 \varpi_{k}-  T_1 \varpi^*) - (T_2 T_1 \varpi_{k}-  T_2  T_1 \varpi^*)$ and $y=(T_1 \varpi_{k}- T_1 \varpi^*)-(\varpi_{k}-\varpi^*)$.

Combining \eqref{equ_thm_57_1} and \eqref{equ_thm_57_2} yields $\forall k\geq 0$,
\begin{equation}\label{equ_thm_5}
||\varpi_{k+1}- \varpi^* ||_{\Phi}^2 \leq   ||\varpi_{k}-  \varpi^*  ||^2_{\Phi} - (1-\xi) || \varpi_{k}-\varpi_{k+1}||^2_{\Phi}.
\end{equation}
Using \eqref{equ_thm_5} from $0$ to $k$ and adding all $k+1$ inequalities yields
\begin{equation}\nonumber
||\varpi_{k+1}- \varpi^* ||_{\Phi}^2 \leq   ||\varpi_{0}-  \varpi^*  ||^2_{\Phi} - (1-\xi) \sum_{l=0}^{k}|| \varpi_{l}-\varpi_{l+1}||^2_{\Phi}.
\end{equation}
Taking limit as $k\rightarrow \infty$ we have,
$(1-\xi) \sum_{k=1}^{\infty} || \varpi_{k}-\varpi_{k+1}||^2_{\Phi} \leq ||\varpi_{0}-  \varpi^*  ||^2_{\Phi} $.
Since $1-\xi>0$, it follows that $\sum_{k=1}^{\infty} || \varpi_{k}-\varpi_{k+1}||^2_{\Phi}$ converges and  $\lim_{k\rightarrow \infty} \varpi_{k}-\varpi_{k+1}=0$.

Since $\{|| \varpi_{k}-  \varpi^* || \}$ is bounded   and converges, $\{\varpi_k\}$ is a bounded sequence. There exists a subsequence$\{\varpi_{n_k}\}$  that converges to $\tilde{\varpi}^*$.
Notice that the composition  $T_2\circ T_1$ is (Lipschitz) continuous and single-valued,   because \eqref{com_fix}
is just an equivalent expression of
Algorithm \ref{dal_1}, and obviously the right hand side of Algorithm \ref{dal_1} is continuous.
${\varpi}_{n_{k}+1}=T_2T_1 \varpi_{n_k}$. Since $T_2T_1$ is continuous,
 and $\lim_{n_{k}\rightarrow \infty}\varpi_{n_k}-{\varpi}_{n_{k}+1}=0$, passing to limiting point, we have $\tilde{\varpi}^* =T_2T_1\tilde{\varpi}^* $. Therefore, the limiting point $\tilde{\varpi}^*$  is a fixed point of $T_2T_1$, or equivalently, $\tilde{\varpi}^* \in zer(\bar{\mathfrak{A}}+\bar{\mathfrak{B}})$.

Setting $\varpi^*=\tilde{\varpi}^*$ in \eqref{equ_thm_5_2}, we have  $\{|| \varpi_{k}-  \tilde{\varpi}^* || \}$ is bounded   and converges. Since there exists a subsequence $\{\varpi_{n_k}\}$ that converges to $\tilde{\varpi}^*$, it follows that  $\{|| \varpi_{k}-  \tilde{\varpi}^* || \}$ converges to zero. Therefore, $\lim_{k\rightarrow \infty} \varpi_k \rightarrow \tilde{\varpi}^*$.
By Theorem \ref{thm_zero_is_correct}, this just implies  \eqref{con_1}.
\hfill $\Box$

\section{Distributed  algorithm with inertia}\label{sec_alg_inertial}

In this section, we propose a distributed algorithm with  inertia for variational GNE seeking,
which possibly accelerates  the convergence under some mild additional computation burden.

There are various modifications of Picard fixed point iteration to achieve  the possible acceleration of convergence speed, and most
of them fall into the domains of {\it relaxation algorithm} and {\it inertial algorithm} (Refer to \cite{hendrickx} for reviews and numerical comparisons for  optimization problems). The relaxation algorithm that simply combines the current operator output with previous iterate, leads to the well-known Krasnosel'ski\u{i}-Mann type of fixed point  iteration\footnote{Chapter 5 of \cite{combettes1}}, and
has been utilized in (generalized)  Nash equilibrium computation in \cite{krawczyk1} and \cite{krawczyk2}. Meanwhile, inertial
algorithms  in operator splitting methods have received attention in recent years, such as
\cite{attouch1}, \cite{attouch2}, \cite{intertia1} and \cite{vu}.　
These efforts are partially motivated by the heavy ball method in \cite{polyak} and Nesterov's acceleration algorithm in \cite{nesterov}) for optimization problems and their recent success in machine learning applications (refer to \cite{jordan}).
In particular,  Nesterov's acceleration algorithm is proved to enjoy an  optimal convergence speed with a specific step-size choice.
Thereby, in this work we consider a distributed algorithm with inertia for  variational GNE seeking given as below:

\begin{alg}\label{dal_3}
\quad \\
\noindent\rule{0.49\textwidth}{0.5mm}
\begin{align}
&Acceleration \; phase \nonumber\\
&\tilde{x}_{i,k}       = x_{i,k}       + \alpha(x_{i,k}-x_{i,k-1})\\
& \tilde{z}_{i,k}      = z_{i,k}       + \alpha(z_{i,k}-z_{i,k-1})\\
&\tilde{\lambda}_{i,k} = \lambda_{i,k} + \alpha(\lambda_{i,k}-\lambda_{i,k-1})\\
&Update \; phase \nonumber\\
&x_{i,k+1}   =  P_{\Omega_i}\big{[}\tilde{x}_{i,k}-\tau_i ( \nabla_{x_i}f_i(\tilde{x}_{i,k}, \{\tilde{x}_{j,k}\}_{j\in \mathcal{N}_i^{f}}) -A_i^T \tilde{\lambda}_{i,k})\big{]} \label{ial_1}\\
&z_{i,k+1}   =  \tilde{z}_{i,k}                + \nu_i  \sum_{j\in \mathcal{N}^{\lambda}_i} w_{ij}(\tilde{\lambda}_{i,k}-\tilde{\lambda}_{j,k}) \label{ial_2}\\
&\lambda_{i,k+1} = P_{\mathbf{R}^{m}_{+}}\Big{\{} \tilde{\lambda}_{i,k} - \sigma_i\big{[} 2 A_i x_{i,k+1}-A_i\tilde{x}_{i,k}-b_i + \sum_{j\in \mathcal{N}^{\lambda}_i} w_{ij}(\tilde{\lambda}_{i,k}-\tilde{\lambda}_{j,k}) \nonumber \\
&\quad \qquad +2\sum_{j\in \mathcal{N}^{\lambda}_i} w_{ij}(z_{i,k+1}-z_{j,k+1})-\sum_{j\in \mathcal{N}^{\lambda}_i} w_{ij}(\tilde{z}_{i,k}-\tilde{z}_{j,k})  \big{]} \Big{\}}   \label{ial_3}
\end{align}
\noindent\rule{0.49\textwidth}{0.5mm}
$\alpha >0$ is a fixed step-size in the acceleration phase,  and $\tau_i >0, \nu_i >0,  \sigma_i >0 $ are fixed step-sizes  of player $i$, and $W=[w_{ij}]$ is the weighted adjacency matrix of multiplier graph $\mathcal{G}_{\lambda}$.
\end{alg}

Compared with Algorithm \ref{dal_1}, Algorithm \ref{dal_3} has two phases.
In the acceleration phase, each player  uses the local state information of the last two steps to get
 predictive variables by a simple linear extrapolation.  In the update phase, the players just feed the predictive variables to Algorithm \ref{dal_1} to get the next iterates. Hence compared with Algorithm \ref{dal_1},
  Algorithm \ref{dal_3} has only an additional simple
 local computation burden. Obviously, Algorithm \ref{dal_3} is also totally distributed, and shares all the features of Algorithm \ref{dal_1}. However, there is an additional need to choose a proper  step-size $\alpha$.

In the following two subsections, we will first give some intuitive interpretation of Algorithm \ref{dal_3} from the
viewpoint of a discretization of continuous-time dynamical systems, and then prove its convergence.

\subsection{Interpretations from viewpoints of dynamical systems}

The interpretation of inertial (acceleration) algorithms from a continuous-time dynamical system viewpoint can be found in \cite{polyak} and most recently in \cite{jordan} for optimization problems and in \cite{attouch2} for proximal point algorithms. Here we give a  comparative development of  Algorithm \ref{dal_1} and Algorithm \ref{dal_3} just  for illustrations of the differences behind the algorithms.

Firstly, let us show that Algorithm \ref{dal_1}, or equivalently its compact reformulation \eqref{com_fix} in Lemma \ref{lem_fix_approximation},  can be interpreted as the discretization of the following dynamical system:
\begin{equation}\label{ode_1}
\dot{\varpi } \in -\Phi^{-1}\bar{\mathfrak{A}}(\varpi)-\Phi^{-1}\bar{\mathfrak{B}}(\varpi).
\end{equation}
In fact, for differential inclusion \eqref{ode_1}, we have the following implicit/explicit  discretization with step-size of $h$,
\begin{equation}\label{equ_discretization_1}
\frac{\varpi({kh+h})-\varpi(kh)}{h} \in -\Phi^{-1}\bar{\mathfrak{A}}(\varpi(kh)) - \Phi^{-1}\bar{\mathfrak{B}}(\varpi(kh+h)).
\end{equation}
Denote $\varpi_{k}=\varpi(kh)$ and take $h=1$,  then  \eqref{equ_discretization_1} can be written as
\begin{equation}
\varpi_{k}-\Phi^{-1}\bar{\mathfrak{A}}(\varpi_{k}) \in \varpi_{k+1} + \Phi^{-1}\bar{\mathfrak{B}}(\varpi_{k+1}).
\end{equation}
Therefore, the implicit/explicit discretization of \eqref{ode_1} is
exactly \eqref{compact_operator_2} that leads to  \eqref{com_fix}, or equivalently Algorithm \ref{dal_1}.
Moreover, the explicit discretization corresponds with the forward step,
and the implicit discretization corresponds with the backward step.
 That's the reason why Algorithm \ref{dal_1} is called a forward-backward splitting algorithm.

Adopt similar compact notations as  in Section \ref{sec_algorithm_develpo},  and denote $\mathbf{\tilde{x}}=col(\tilde{x}_1,\cdots,\tilde{x}_N)$,
$\bar{\tilde{\lambda}}= col(\tilde{\lambda}_1,\cdots,\tilde{\lambda}_N)$, and $\bar{\tilde{z}}=col(\tilde{z}_1,\cdots,\tilde{z}_N)$. And further denote $\tilde{\varpi}=col(\mathbf{\tilde{x}},\bar{\tilde{\lambda}},\bar{\tilde{z}}  )$.
Then by similar  arguments as in Section \ref{sec_algorithm_develpo} and using operators $\bar{\mathfrak{A}}$ and $\bar{\mathfrak{B}}$ defined in \eqref{operator_A_B} and $\Phi$ defined in \eqref{metric_matrix},   Algorithm \ref{dal_3} can be written in a compact form (assume $\bar{\mathfrak{B}}$ is maximally monotone),
\begin{align}
&\tilde{\varpi}_{k}  =  \varpi_{k}+\alpha(\varpi_k-\varpi_{k-1})\label{cal_acc_1}\\
&{\varpi}_{k+1}      = ({\rm Id}+\Phi^{-1}\bar{\mathfrak{B}})^{-1}({\rm Id}-\Phi^{-1}\bar{\mathfrak{A}})\tilde{\varpi}_{k} \label{cal_acc_2}
\end{align}
where $\varpi_k$ is defined as in Section \ref{sec_algorithm_develpo}.

We can show that Algorithm \ref{dal_3},  or equivalently  \eqref{cal_acc_1}-\eqref{cal_acc_2},
can be interpreted as the discretization of the following second-order continuous-time dynamical system,
\begin{equation}\label{ode_2}
\ddot{\varpi}+\tilde{\alpha }\dot{\varpi} \in -\Phi^{-1}\bar{\mathfrak{A}}(\varpi)-\Phi^{-1}\bar{\mathfrak{B}}(\varpi).
\end{equation}
In fact, for differential inclusion \eqref{ode_2} consider the following type of
implicit/explicit discretization,
\begin{align}
&\frac{\varpi(kh+h)-2\varpi(kh)+\varpi(kh-h)}{h^2} +  \tilde{\alpha}\frac{\varpi(kh)-\varpi(kh-h)}{h}  \in  \nonumber\\
&\qquad -\Phi^{-1}\bar{\mathfrak{A}}(\tilde{\varpi}(kh)) - \Phi^{-1}\bar{\mathfrak{B}}({\varpi}(kh+h)), \label{equ_discretization_2}
\end{align}
where $\tilde{\varpi}(kh)$ is an  interpolation  point to be determined later.
 Denote $\varpi_{k}=\varpi(kh)$ and take $h=1$, then \eqref{equ_discretization_2} can be written as
\begin{equation}\label{equ_discretization_3}
\varpi_{k}+(1-\tilde{\alpha})(\varpi_{k}-\varpi_{k-1})- \Phi^{-1}\bar{\mathfrak{A}}(\tilde{\varpi}_{k}) \in \varpi_{k+1} + \Phi^{-1}\bar{\mathfrak{B}}(\varpi_{k+1}).
\end{equation}
Denote $\alpha=1-\tilde{\alpha}$ and take $\tilde{\varpi}_{k}= \varpi_{k}+\alpha(\varpi_{k}-\varpi_{k-1}) $, then \eqref{equ_discretization_3} can be  written as
\begin{equation}
\begin{array}{l}\label{equ_discretization_4}
\tilde{\varpi}_{k} = \varpi_{k}+\alpha(\varpi_{k}-\varpi_{k-1}), \\
\tilde{\varpi}_{k}- \Phi^{-1}\bar{\mathfrak{A}}(\tilde{\varpi}_{k}) \in \varpi_{k+1} + \Phi^{-1}\bar{\mathfrak{B}}(\varpi_{k+1}).
\end{array}
\end{equation}
\eqref{equ_discretization_4} leads to  equations \eqref{cal_acc_1}-\eqref{cal_acc_2}, or equivalently Algorithm \ref{dal_3}.

\begin{rem}
Compared with \eqref{ode_1}, \eqref{ode_2} is a second order dynamical system with  an additional inertial term $\alpha \dot{\varpi}$, hence \eqref{ode_2}  enjoys better convergence properties than \eqref{ode_1}. Therefore,
 it is expected that Algorithm \ref{dal_3}, as a discretization of \eqref{ode_2}, would have  better convergence properties than Algorithm \ref{dal_1}.
\end{rem}

\subsection{Convergence analysis}

The following result proves the  convergence of Algorithm \ref{dal_3} by providing sufficient
step-size choices for $\alpha$ as well as $\tau_i, \nu_i, \sigma_i$.  The sufficient choice condition for $\alpha$  can be ensured by solving a simple algebra inequality.
The proof idea of the following result is motivated by  inertial algorithms works for optimization and operator splitting  such as
\cite{attouch1}, \cite{attouch2}, \cite{vu}, \cite{hendrickx}, and especially \cite{intertia1}.
However, since this work considers a noncooperative game setup and adopts a fixed step-size in the distributed algorithm, Theorem \ref{thm_inertial}'s  proof is also provided  for completeness.

\begin{thm}\label{thm_inertial}
Suppose Assumptions \ref{assum1}-\ref{assum3} hold.
Take $0<\beta \leq \min\{\frac{1}{2d^*}, \frac{\eta}{ \theta^2}\}$ where $d^*$ is the maximal weighted degree of multiplier graph $\mathcal{G}_{\lambda}$, and $\eta,\theta$ are parameters in Assumption \ref{assum2}.
Given a sufficient small $0<\epsilon<1$, take $ \delta > \frac{1}{2\beta}$ and  $0<\alpha<1$ in Algorithm \ref{dal_3} such that $2\beta\delta(1-3\alpha-\epsilon) \geq (1-\alpha)^2$.
Suppose that player $i$ chooses  its step-sizes $\tau_i, \nu_i, \sigma_i$ in Algorithm \ref{dal_3} satisfying \eqref{step_size_choice}.
Then with  Algorithm \ref{dal_3}, players' local strategies converge
to the variational GNE of game in \eqref{GM}, and  the  local  multipliers $\lambda_{i,k}$ of all the agents converge to the same multiplier corresponding with KKT condition \eqref{kkt_2}, i.e.,
\begin{equation}\label{con_2}
\lim_{k\rightarrow \infty} x_{i,k}\rightarrow x_i^*,  \quad \lim_{k\rightarrow \infty}  \lambda_{i,k} \rightarrow \lambda^*, \forall i=1,\cdots,N.
\end{equation}
\end{thm}
{\bf Proof:}
By the choice of $\beta$, operator $\bar{\mathfrak{A}}$ is $\beta-$cocoercive
and operator $\bar{\mathfrak{B}}$ is maximally monotone due to Lemma \ref{lem_monotone}.
By the choice of $\delta, \tau_i, \nu_i$ and $ \sigma_i$,  $\Phi-\delta I_{n+2mN}$ is positive semi-definite due to Lemma \ref{lem_metric_matrix}, and $\Phi^{-1}\bar{\mathfrak{A}}$ and $\Phi^{-1}\bar{\mathfrak{B}}$ satisfy  the properties in Lemma \ref{lem_monotone_metric}.
Therefore, Algorithm \ref{dal_3} can be exactly  written in  the compact form of \eqref{cal_acc_1}-\eqref{cal_acc_2} with similar arguments as in  Lemma \ref{lem_fix_approximation}.

Resorting to Theorem \ref{thm_zero_is_correct}, we only need to show that Algorithm \ref{dal_3} converges and its limiting point belongs to $zer(\bar{\mathfrak{A}}+\bar{\mathfrak{B}})$.
In fact, any limiting point of  Algorithm \ref{dal_3} satisfies $\varpi^* \in zer(\bar{\mathfrak{A}} + \bar{\mathfrak{B}})$ as shown next.
Suppose $\lim_{k\rightarrow \infty} \varpi_k \rightarrow \varpi^*$, then $\tilde{\varpi}_k \rightarrow \varpi^* $ and $\varpi^*=({\rm Id}+\Phi^{-1}\bar{\mathfrak{B}})^{-1}({\rm Id}-\Phi^{-1}\bar{\mathfrak{A}})\varpi^*$ using the continuity of the right hand of \eqref{cal_acc_2}. Therefore, $\varpi^* \in zer(\bar{\mathfrak{A}} + \bar{\mathfrak{B}})$ because $\Phi$ is a positive definite matrix.

 The following relationship (similar  with the cosine rule) will be heavily utilized in the convergence analysis.
\begin{equation}\label{lemma_norm_equ_known}
||a-c||_{Q}^2 - ||b-c||_{Q}^2 =2 \langle a-b, a-c \rangle_{Q} - ||a-b ||^2_{Q},
\end{equation}
which  can be verified  by directly expanding with $ ||a+b ||^2_{Q}= ||a ||^2_{Q}+2 \langle a, b\rangle_{Q}+ ||b ||^2_{Q}$.

The proof is divided into three parts:
\begin{itemize}
\item {\it Part 1:} Given any $\varpi^* \in zer(\bar{\mathfrak{A}}+\bar{\mathfrak{B}})$, $||\varpi_{k+1}-\varpi^* ||^2_{\Phi}-||\varpi_k -  \varpi^* ||^2_{\Phi}$ follows a recursive inequality as follows,
\begin{equation}\label{equ_thm_6_13}
  \begin{array}{l}
  ||\varpi_{k+1}-\varpi^* ||^2_{\Phi}-||\varpi_k -  \varpi^* ||^2_{\Phi} \\
  \leq \alpha(||\varpi_{k}-\varpi^* ||^2_{\Phi}-||\varpi_{k-1} -  \varpi^* ||^2_{\Phi})+ (\alpha+\alpha^2) ||\varpi_{k}-\varpi_{k-1} ||^2_{\Phi}
\end{array}
\end{equation}
\item {\it Part 2:} Given any $\varpi^* \in zer(\bar{\mathfrak{A}}+\bar{\mathfrak{B}})$, $\sum_{k}^{\infty} || \varpi_{k+1}-\varpi_{k}||^2_{\Phi}  < \infty$ and $\lim_{k\rightarrow \infty}\varpi_{k+1}-\varpi_{k}=0$.
\item {\it Part 3:} We first show the convergence of $ \{||\varpi_{k}-\varpi^* ||^2_{\Phi}\}$ given any $\varpi^* \in zer(\bar{\mathfrak{A}}+\bar{\mathfrak{B}})$, and then show the convergence of Algorithm \ref{dal_3}.
\end{itemize}

\vskip 3mm
{\it Part 1}:  Given any point $\varpi^* \in zer(\bar{\mathfrak{A}} + \bar{\mathfrak{B}})$, we first prove a
recursive inequality \eqref{equ_thm_6_13} for $||\varpi_{k+1}-\varpi^* ||^2_{\Phi}-||\varpi_k -  \varpi^* ||^2_{\Phi}$.
\vskip 1mm

Using \eqref{lemma_norm_equ_known} to expand the left hand of \eqref{equ_thm_6_13} yields,
\begin{equation}
\begin{array}{ll}\label{equ_thm_6_1}
&||\varpi_k -  \varpi^* ||^2_{\Phi}  - ||\varpi_{k+1}-\varpi^* ||^2_{\Phi}\\
&= 2\langle \varpi_k-\varpi_{k+1},  \varpi_k -  \varpi^*  \rangle_{\Phi} - || \varpi_k-\varpi_{k+1} ||^2_{\Phi} \\
&=  2\langle \varpi_k-\varpi_{k+1},  \varpi_k-\varpi_{k+1} +\varpi_{k+1}-  \varpi^*  \rangle_{\Phi}
- || \varpi_k-\varpi_{k+1} ||^2_{\Phi}\\
&= || \varpi_k-\varpi_{k+1} ||^2_{\Phi}+ 2\langle \varpi_k-\varpi_{k+1},  \varpi_{k+1}-  \varpi^* \rangle_{\Phi} \\
&= || \varpi_k-\varpi_{k+1} ||^2_{\Phi}+ 2\langle \tilde{\varpi}_{k} -\varpi_{k+1},  \varpi_{k+1}-  \varpi^* \rangle_{\Phi} \\
& \qquad -2\alpha \langle {\varpi}_{k} -\varpi_{k-1},  \varpi_{k+1}-  \varpi^* \rangle_{\Phi}
\end{array}
\end{equation}
where the last step is derived by incorporating \eqref{cal_acc_1}.

To tackle the second term on the right hand of \eqref{equ_thm_6_1}, we proceed as follows.
By \eqref{cal_acc_2} we have,
$
\tilde{\varpi}_{k}- \Phi^{-1}\bar{\mathfrak{A}}(\tilde{\varpi}_{k}) \in  \varpi_{k+1} + \Phi^{-1}\bar{\mathfrak{B}}(\varpi_{k+1}), \nonumber
$
or equivalently,
\begin{equation}\label{equ_thm_6_2}
\Phi(\tilde{\varpi}_{k}-\varpi_{k+1})- \bar{\mathfrak{A}}(\tilde{\varpi}_{k}) \in \bar{\mathfrak{B}}(\varpi_{k+1}).
\end{equation}
Because $\varpi^* \in zer(\bar{\mathfrak{A}}+\bar{\mathfrak{B}})$, we also have
\begin{equation}\label{equ_thm_6_3}
-\bar{\mathfrak{A}}(\varpi^*) \in \bar{\mathfrak{B}}(\varpi^*).
\end{equation}
Due to the maximal monotonicity of $\bar{\mathfrak{B}}$ proved in Lemma \ref{lem_monotone},
\begin{equation}\label{equ_thm_6_4}
\langle u-v, \varpi_{k+1}-\varpi^*  \rangle \geq 0, \forall\; u \in \bar{\mathfrak{B}}(\varpi_{k+1}), v\in \bar{\mathfrak{B}}(\varpi^*)
\end{equation}
By incorporating \eqref{equ_thm_6_2} and \eqref{equ_thm_6_3} into \eqref{equ_thm_6_4}, we have
\begin{equation}
\langle\Phi(\tilde{\varpi}_{k}-\varpi_{k+1})-\bar{ \mathfrak{A}}(\tilde{\varpi}_{k})+ \bar{\mathfrak{A}}(\varpi^*), \varpi_{k+1}-\varpi^*  \rangle \geq 0 \nonumber
\end{equation}
or equivalently,
\begin{equation}\label{equ_thm_6_5}
\langle \tilde{\varpi}_{k}-\varpi_{k+1},\varpi_{k+1}-\varpi^*  \rangle_{\Phi} - \langle \bar{\mathfrak{A}}(\tilde{\varpi}_{k})- \bar{\mathfrak{A}}(\varpi^*), \varpi_{k+1}-\varpi^*  \rangle \geq 0.
\end{equation}

Using \eqref{equ_thm_6_5} for the second term on the right hand of  \eqref{equ_thm_6_1} yields
\begin{equation}\label{equ_thm_6_6}
\begin{array}{ll}
&||\varpi_k -  \varpi^* ||^2_{\Phi}  - ||\varpi_{k+1}-\varpi^* ||^2_{\Phi} \\
&\geq || \varpi_k-\varpi_{k+1} ||^2_{\Phi}
+ 2\langle \bar{\mathfrak{A}}(\tilde{\varpi}_{k})-\bar{ \mathfrak{A}}(\varpi^*), \varpi_{k+1}-\varpi^*  \rangle \\
&-2\alpha \langle {\varpi}_{k} -\varpi_{k-1},  \varpi_{k+1}-  \varpi^* \rangle_{\Phi}
\end{array}
\end{equation}

By Lemma \ref{lem_monotone}, $\bar{\mathfrak{A}}$ is $\beta-$cocoercive.
For the second term on the right hand of \eqref{equ_thm_6_6}, we have
 \begin{equation}
 \begin{array}{lll}\label{equ_thm_6_7}
 & \langle \bar{\mathfrak{A}}(\tilde{\varpi}_{k})-\bar{\mathfrak{A}}(\varpi^*), \varpi_{k+1}-\varpi^*  \rangle  \\
 &=
   \langle \bar{\mathfrak{A}}(\tilde{\varpi}_{k})-\bar{\mathfrak{A}}(\varpi^*), \varpi_{k+1}-\tilde{\varpi}_{k}+\tilde{\varpi}_{k}-\varpi^*  \rangle  \\
 &  \geq {\beta} ||\bar{\mathfrak{A}}(\tilde{\varpi}_{k})-\bar{\mathfrak{A}}(\varpi^*)||_2^2 + 2 \langle \sqrt{\beta}(\bar{\mathfrak{A}}(\tilde{\varpi}_{k})-\bar{\mathfrak{A}}(\varpi^*)), \frac{1}{2\sqrt{\beta}}({\varpi}_{k+1}-\tilde{\varpi}_{k})  \rangle \\
 & \geq \beta ||\bar{\mathfrak{A}}(\tilde{\varpi}_{k})-\bar{\mathfrak{A}}(\varpi^*)||_2^2 - \beta || \bar{\mathfrak{A}}(\tilde{\varpi}_{k})-\bar{\mathfrak{A}}(\varpi^*)||_2^2 -\frac{1}{4\beta} ||{\varpi}_{k+1}-\tilde{\varpi}_{k}||_2^2\\
 & = -\frac{1}{4\beta} ||{\varpi}_{k+1}-\tilde{\varpi}_{k}||_2^2
  \end{array}
 \end{equation}
where the first inequality is obtained by the cocoercive property \eqref{cocoercive} and the second inequality is obtained by $2\langle a, b\rangle \geq -||a ||_2^2-||b ||_2^2$.

Combining \eqref{equ_thm_6_6} with \eqref{equ_thm_6_7} we have
 \begin{equation}
 \begin{array}{ll}\nonumber
&||\varpi_k -  \varpi^* ||^2_{\Phi}  - ||\varpi_{k+1}-\varpi^* ||^2_{\Phi}\geq || \varpi_k-\varpi_{k+1} ||^2_{\Phi}
\\
&-\frac{1}{2\beta} ||{\varpi}_{k+1}-\tilde{\varpi}_{k}||_2^2-2\alpha \langle {\varpi}_{k} -\varpi_{k-1},  \varpi_{k+1}-  \varpi^* \rangle_{\Phi},
\end{array}
\end{equation}
or equivalently we have
\begin{equation}\label{equ_thm_6_8}
\begin{array}{ll}
||\varpi_{k+1}-\varpi^* ||^2_{\Phi}-||\varpi_k -  \varpi^* ||^2_{\Phi}\leq  -|| \varpi_k-\varpi_{k+1} ||^2_{\Phi}
\\+\frac{1}{2\beta} ||{\varpi}_{k+1}-\tilde{\varpi}_{k}||_2^2
+2\alpha \langle {\varpi}_{k} -\varpi_{k-1},  \varpi_{k+1}-  \varpi^* \rangle_{\Phi}.
\end{array}
\end{equation}

Utilizing the  equality \eqref{lemma_norm_equ_known} we also have,
\begin{equation}\label{equ_thm_6_9}
\begin{array}{ll}
 &\alpha(||\varpi_{k}-\varpi^* ||^2_{\Phi}-||\varpi_{k-1} -  \varpi^* ||^2_{\Phi})\\
 &=
 \alpha( 2 \langle \varpi_{k}-\varpi_{k-1}, \varpi_{k}-\varpi^{*} \rangle_{\Phi} - ||\varpi_{k}-\varpi_{k-1} ||^2_{\Phi} )
\end{array}
\end{equation}

Combining \eqref{equ_thm_6_8} and \eqref{equ_thm_6_9}, we have,
\begin{equation}
\begin{array}{ll}\label{equ_thm_6_10}
&||\varpi_{k+1}-\varpi^* ||^2_{\Phi}-||\varpi_k -  \varpi^* ||^2_{\Phi} \\
&\quad-\alpha(||\varpi_{k}-\varpi^* ||^2_{\Phi}-||\varpi_{k-1} -  \varpi^* ||^2_{\Phi})\\
&\leq  -|| \varpi_k-\varpi_{k+1} ||^2_{\Phi}
+\frac{1}{2\beta} ||{\varpi}_{k+1}-\tilde{\varpi}_{k}||_2^2  \\
&+
2\alpha \langle {\varpi}_{k} -\varpi_{k-1},  \varpi_{k+1}-  \varpi^* \rangle_{\Phi} \\
& -\alpha( 2 \langle \varpi_{k}-\varpi_{k-1}, \varpi_{k}-\varpi^{*} \rangle_{\Phi} - ||\varpi_{k}-\varpi_{k-1} ||^2_{\Phi} ) \\
& \leq  -|| \varpi_k-\varpi_{k+1} ||^2_{\Phi}
+\frac{1}{2\beta} ||{\varpi}_{k+1}-\tilde{\varpi}_{k}||_2^2  \\
&+
2\alpha \langle {\varpi}_{k} -\varpi_{k-1},  \varpi_{k+1}-  \varpi_k\rangle_{\Phi} +\alpha ||\varpi_{k}-\varpi_{k-1} ||^2_{\Phi}
\end{array}
\end{equation}

Next using \eqref{cal_acc_1} and  \eqref{lemma_norm_equ_known} for the first and third terms on the right hand of \eqref{equ_thm_6_10} yields
\begin{equation}
  \begin{array}{lll}\label{equ_thm_6_11}
 & \;&2\alpha \langle {\varpi}_{k} -\varpi_{k-1},  \varpi_{k+1}-  \varpi_k\rangle_{\Phi}-|| \varpi_k-\varpi_{k+1} ||^2_{\Phi} \\
  \; &= &2\alpha \langle \frac{1}{\alpha }(\tilde{\varpi}_k-\varpi_k),  \varpi_{k+1}-  \varpi_k\rangle_{\Phi}-|| \varpi_k-\varpi_{k+1} ||^2_{\Phi}  \\
  \;&=& 2 \langle  \tilde{\varpi}_k-\varpi_k,  \varpi_{k+1}-  \varpi_k\rangle_{\Phi}-|| \varpi_k-\varpi_{k+1} ||^2_{\Phi}\\
    \;&=&  ||\varpi_{k}-\tilde{\varpi}_{k}　||^2_{\Phi}- ||\varpi_{k+1}-\tilde{\varpi}_{k} ||_{\Phi}\\
    \;&=& \alpha^2 || \varpi_{k}-{\varpi}_{k-1}||^2_{\Phi}- ||\varpi_{k+1}-\tilde{\varpi}_{k} ||^2_{\Phi}
  \end{array}
\end{equation}

For the second term on the right hand of \eqref{equ_thm_6_10}, $\frac{1}{2\beta} ||{\varpi}_{k+1}-\tilde{\varpi}_{k}||_2^2 \leq \frac{1}{2\beta \delta} ||{\varpi}_{k+1}-\tilde{\varpi}_{k}||^2_{\Phi}$,  since $\Phi > \delta I_{n+2mN}$.
Incorporating this and  \eqref{equ_thm_6_11} into \eqref{equ_thm_6_10}
\begin{equation}
\begin{array}{ll}\label{equ_thm_6_12}
& ||\varpi_{k+1}-\varpi^* ||^2_{\Phi}-||\varpi_k -  \varpi^* ||^2_{\Phi} \\
&\quad -
 \alpha(||\varpi_{k}-\varpi^* ||^2_{\Phi}-||\varpi_{k-1} -  \varpi^* ||^2_{\Phi})\\
& \leq
- ||\varpi_{k+1}-\tilde{\varpi}_{k} ||^2_{\Phi}+\frac{1}{2\beta} ||{\varpi}_{k+1}-\tilde{\varpi}_{k}||_2^2  + (\alpha+\alpha^2) ||\varpi_{k}-\varpi_{k-1} ||_{\Phi} \\
& \leq -(1-\frac{1}{2\delta\beta})  ||\varpi_{k+1}-\tilde{\varpi}_{k} ||^2_{\Phi} + (\alpha+\alpha^2) ||\varpi_{k}-\varpi_{k-1} ||^2_{\Phi}
\end{array}
\end{equation}
Since $\delta > \frac{1}{2\beta}$, we derive \eqref{equ_thm_6_13}.

\vskip 3mm
{\it Part 2}: In this step, we will prove $\sum_{k}^{\infty} || \varpi_{k+1}-\varpi_{k}||^2_{\Phi}  < \infty$ and $\lim_{k\rightarrow \infty}\varpi_{k+1}-\varpi_{k}=0$.

\vskip 1mm
Denote $S=\Phi- \frac{1}{2\beta}I_{n+2mN}$. Then $S$ is symmetric and positive definite  since $\Phi\geq \delta I_{n+2mN}$ and $\delta > \frac{1}{2\beta}$.
The first inequality of \eqref{equ_thm_6_12} can also be written as:
\begin{equation}
\begin{array}{l} \label{equ_thm_6_14}
||\varpi_{k+1}-\varpi^* ||^2_{\Phi}-||\varpi_k -  \varpi^* ||^2_{\Phi}-
 \alpha(||\varpi_{k}-\varpi^* ||^2_{\Phi}-||\varpi_{k-1} -  \varpi^* ||^2_{\Phi})\\
\leq - ||\varpi_{k+1}-\tilde{\varpi}_{k} ||^2_{S}   + (\alpha+\alpha^2) ||\varpi_{k}-\varpi_{k-1} ||^2_{\Phi}.
\end{array}
\end{equation}

For the first term on the right hand of \eqref{equ_thm_6_14},
\begin{equation}
\begin{array}{ll}\label{equ_thm_6_15}
 & - ||\varpi_{k+1}-\tilde{\varpi}_{k} ||^2_{S}\\
 &= 2 \langle \varpi_k-\varpi_{k+1}, \varpi_k-\tilde{\varpi}_k \rangle_{S} - ||\varpi_k-\varpi_{k+1} ||^2_{S} -||\varpi_{k}-\tilde{\varpi}_k||_{S}^2 \\
 &= -2 \langle \varpi_k-\varpi_{k+1}, \alpha(\varpi_{k}-\varpi_{k-1}) \rangle_{S} \\
 &- ||\varpi_k-\varpi_{k+1} ||^2_{S} -\alpha^2||\varpi_{k}-\varpi_{k-1}||_{S}^2 \\
 & \leq \alpha  ||\varpi_k-\varpi_{k+1} ||^2_S + \alpha  ||\varpi_{k}-\varpi_{k-1} ||^2_S \\
 &- ||\varpi_k-\varpi_{k+1} ||^2_{S} -\alpha^2||\varpi_{k}-\varpi_{k-1}||_{S}^2 \\
 & =  (\alpha -1) ||\varpi_k-\varpi_{k+1} ||^2_S + (\alpha-\alpha^2)  ||\varpi_{k}-\varpi_{k-1} ||^2_S.
 \end{array}
\end{equation}
where the first equality follows from \eqref{lemma_norm_equ_known}, the second equality follows from \eqref{cal_acc_1}, and the third inequality follows from $-2 \langle x,y \rangle \leq  || x||^2 + ||y ||^2$.

Denote $ Q= 2\Phi- \frac{1-\alpha}{2\beta} I_{n+2mN}$. Then $Q$ is also symmetric and  positive definite, since $\alpha<1$ and $\Phi \geq \delta I_{n+2mN}$.
 Combining \eqref{equ_thm_6_14} with \eqref{equ_thm_6_15},
  \begin{equation}
\begin{array}{lll}\label{equ_thm_6_16}
& ||\varpi_{k+1}-\varpi^* ||^2_{\Phi}-||\varpi_k -  \varpi^* ||^2_{\Phi} \\
&\qquad-
 \alpha(||\varpi_{k}-\varpi^* ||^2_{\Phi}-||\varpi_{k-1} -  \varpi^* ||^2_{\Phi})\\
& \leq (\alpha -1) ||\varpi_k-\varpi_{k+1} ||^2_S + (\alpha-\alpha^2)  ||\varpi_{k}-\varpi_{k-1} ||^2_S   \\
& \qquad + (\alpha+\alpha^2) ||\varpi_{k}-\varpi_{k-1} ||^2_{\Phi} \\
&=(\alpha -1) ||\varpi_k-\varpi_{k+1} ||^2_S\\
 &\quad + \alpha \langle [(1-\alpha)(\Phi- \frac{1}{2\beta}I)+(1+\alpha)\Phi] \varpi_{k}-\varpi_{k-1},\varpi_{k}-\varpi_{k-1} \rangle \\
& = (\alpha -1) ||\varpi_k-\varpi_{k+1} ||^2_S + \alpha ||\varpi_{k}-\varpi_{k-1} ||^2_Q.
\end{array}
\end{equation}

Denote $\mu_k=  ||\varpi_{k}-\varpi^* ||^2_{\Phi} - \alpha|| \varpi_{k-1}-\varpi^*||^2_{\Phi} + \alpha ||\varpi_{k}-\varpi_{k-1} ||^2_Q $, then
\begin{equation}
\begin{array}{ll}\label{equ_thm_6_17}
&\mu_{k+1}-\mu_{k} =||\varpi_{k+1}-\varpi^* ||^2_{\Phi} - \alpha|| \varpi_{k}-\varpi^*||^2_{\Phi} + \alpha ||\varpi_{k+1}-\varpi_{k} ||^2_Q\\
&\quad - ||\varpi_{k}-\varpi^* ||^2_{\Phi} + \alpha|| \varpi_{k-1}-\varpi^*||^2_{\Phi} - \alpha ||\varpi_{k}-\varpi_{k-1} ||^2_Q\\
&\; = ||\varpi_{k+1}-\varpi^* ||^2_{\Phi}- ||\varpi_{k}-\varpi^* ||^2_{\Phi}-\alpha (|| \varpi_{k}-\varpi^*||^2_{\Phi} \\
&\quad-|| \varpi_{k-1}-\varpi^*||^2_{\Phi})
+  \alpha ||\varpi_{k+1}-\varpi_{k} ||^2_Q-\alpha ||\varpi_{k}-\varpi_{k-1} ||^2_Q\\
&\;\leq (\alpha -1) ||\varpi_k-\varpi_{k+1} ||^2_S + \alpha ||\varpi_{k}-\varpi_{k-1} ||^2_Q \\
 &\quad +\alpha ||\varpi_{k+1}-\varpi_{k} ||^2_Q -\alpha ||\varpi_{k}-\varpi_{k-1} ||^2_Q\\
&\; = (\alpha -1) ||\varpi_k-\varpi_{k+1} ||^2_S  +\alpha ||\varpi_{k+1}-\varpi_{k} ||^2_Q \\
&\; = \langle [(\alpha-1)(\Phi-\frac{1}{2\beta}I)+\alpha(2\Phi-\frac{1-\alpha}{2\beta}I)] \varpi_{k+1}-\varpi_{k} , \\
& \qquad\varpi_{k+1}-\varpi_{k} \rangle
 \end{array}
\end{equation}
where the third inequality follows by \eqref{equ_thm_6_16}.

Given a sufficient small $0<\epsilon<1$, choose $0<\alpha<1$  and $\delta > \frac{1}{2\beta}$ such that $2\beta\delta(1-3\alpha-\epsilon) \geq (1-\alpha)^2$, then
\begin{equation}
\begin{array}{ll}\nonumber
&-(\alpha-1)(\Phi-\frac{1}{2\beta}I_{n+2mN})-\alpha(2\Phi-\frac{1-\alpha}{2\beta}I_{n+2mN}) \\
&= (1-3\alpha)\Phi- \frac{(1-\alpha)^2}{2\beta} I_{n+2mN}\geq  \epsilon \Phi
\end{array}
\end{equation}

Therefore,  \eqref{equ_thm_6_17} yields
$\mu_{k+1}-\mu_{k} \leq -\epsilon || \varpi_{k+1}-\varpi_{k}||^2_{\Phi}$. Therefore,
$\mu_{k+1}\leq\mu_{k} \leq \mu_{1}$.
By the definition of $\mu_k$,
$ ||\varpi_{k}-\varpi^* ||^2_{\Phi} - \alpha|| \varpi_{k-1}-\varpi^*||^2_{\Phi} \leq ||\varpi_{k}-\varpi^* ||^2_{\Phi} - \alpha|| \varpi_{k-1}-\varpi^*||^2_{\Phi} + \alpha ||\varpi_{k}-\varpi_{k-1} ||^2_Q \leq \mu_{1} $.
Therefore, $||\varpi_{k}-\varpi^* ||^2_{\Phi} \leq \alpha|| \varpi_{k-1}-\varpi^*||^2_{\Phi}+ \mu_1$.
By Lemma 1 at Page 44 of \cite{polyak},
$||\varpi_{k}-\varpi^* ||^2_{\Phi} \leq \alpha^{k} ( || \varpi_{1}-\varpi^*||^2_{\Phi}-\frac{\mu_1}{1-a}) +\frac{ \mu_1}{1-\alpha}$, and  $\{\varpi_{k}\}$ is bounded sequence.

We also have $\mu_{k+1}-\mu_{1} \leq  -\epsilon \sum_{i=1}^{k} || \varpi_{i+1}-\varpi_{i}||^2_{\Phi}$.
Then $
\epsilon \sum_{i=1}^{k} || \varpi_{i+1}-\varpi_{i}||^2_{\Phi} \leq \mu_{1}-\mu_{k+1}
\leq \mu_{1}+\alpha||\varpi_{k}-\varpi^* ||^2_{\Phi}
\leq \mu_{1}+ \alpha^{k+1}( || \varpi_{1}-\varpi^*||^2_{\Phi}-\frac{\mu_1}{1-a}) +\frac{ \alpha\mu_1}{1-\alpha} $.
Let $k$ goes to infinity, we have,
\begin{equation}
\sum_{k}^{\infty} || \varpi_{k+1}-\varpi_{k}||^2_{\Phi}  < \infty
\end{equation}
Therefore, $\lim_{k\rightarrow \infty}\varpi_{k+1}-\varpi_{k}=0$, and $   \sum_{k=1}^{\infty}(\alpha+\alpha^2) ||\varpi_{k}-\varpi_{k-1} ||^2_{\Phi}  <\infty $.

\vskip 3mm
{\it Part 3}: In this part, we first show the convergence of $ \{||\varpi_{k}-\varpi^* ||^2_{\Phi}\}$ given any $\varpi^* \in zer(\bar{\mathfrak{A}}+\bar{\mathfrak{B}})$, and then show the convergence of
Algorithm \ref{dal_3}.
\vskip 1mm

 Denote $\phi_k=\max\{0,||\varpi_{k}-\varpi^* ||^2_{\Phi}-||\varpi_{k-1} -  \varpi^* ||^2_{\Phi}\}$ and $\psi_k=(\alpha+\alpha^2) ||\varpi_{k}-\varpi_{k-1} ||^2_{\Phi}$, and recall \eqref{equ_thm_6_13},  we  have
 $\phi_{k+1} \leq \alpha \phi_{k} + \psi_k.$ Apply this relationship recursively,
 \begin{equation}\label{equ_thm_6_18}
 \phi_{k+1} \leq \alpha^k \phi_{1}+ \sum_{i=0}^{k-1} \alpha^{i}\psi_{k-i}.
 \end{equation}
 Summing \eqref{equ_thm_6_18} from $k=1$ to $k=J$,
 \begin{equation}
 \begin{array}{lll}
 \sum_{k=1}^J \phi_{k+1} & \leq &  \phi_1 \frac{1-\alpha^J}{1-\alpha} + \sum_{k=1}^J\sum_{i=0}^{k-1} \alpha^{i}\psi_{k-i} \\
 \qquad                  & \leq &  \phi_1 \frac{1-\alpha^J}{1-\alpha} + \sum_{i=0}^{J-1}\alpha^{i} \sum_{k=i+1}^{J}\psi_{k-i}
 \end{array}
 \end{equation}
Let $J\rightarrow \infty$, then since $0<\alpha<1$,
 \begin{equation}
 \begin{array}{lll}
 \sum_{k=1}^{\infty}\phi_{k} &\leq & \frac{\phi_1}{1-\alpha} + \sum_{i=0}^{\infty}\alpha^{i}\sum_{k=i+1}^{\infty}\psi_{k-i} \\
                             &\leq & \frac{\phi_1}{1-\alpha} + \sum_{i=0}^{\infty}\alpha^{i}\sum_{t=1}^{\infty}\psi_{t}\\
                             &\leq & \frac{\phi_1}{1-\alpha} + \frac{1}{1-\alpha}\sum_{t=1}^{\infty}\psi_{t}
 \end{array}
 \end{equation}
Noticing that $\sum_{t=1}^{\infty}\psi_{t}=\sum_{k=1}^{\infty}(\alpha+\alpha^2) ||\varpi_{k}-\varpi_{k-1} ||^2_{\Phi} <\infty$,
$ \sum_{k=1}^{\infty}\phi_{k} <\infty$, and hence the sequence $\{ \sum_{i=1}^k\phi_{i}\}$, being a nonnegative and non-decreasing sequence,  converges and is bounded.

Consider another  sequence $\{ ||\varpi_{k}-\varpi^* ||^2_{\Phi}-\sum_{i=1}^k\phi_{i} \}$.
Since $||\varpi_{k}-\varpi^* ||^2_{\Phi}$ is nonnegative and $\{ \sum_{i=1}^k\phi_{i}\}$ is bounded,  $\{ ||\varpi_{k}-\varpi^* ||^2_{\Phi}-\sum_{i=1}^k\phi_{i} \}$ is bounded from below. Furthermore, $\{ ||\varpi_{k}-\varpi^* ||^2_{\Phi}-\sum_{i=1}^k\phi_{i} \}$ is a non-increasing sequence. In fact,
 \begin{equation}
 \begin{array}{lll}\nonumber
 &\;&||\varpi_{k+1}-\varpi^* ||^2_{\Phi}-\sum_{i=1}^{k+1}\phi_{i} \\
 &=& ||\varpi_{k+1}-\varpi^* ||^2_{\Phi}- \phi_{k+1}-\sum_{i=1}^{k}\phi_{i}\\
  &\leq&  ||\varpi_{k+1}-\varpi^* ||^2_{\Phi}-||\varpi_{k+1}-\varpi^* ||^2_{\Phi}+||\varpi_{k} -  \varpi^* ||^2_{\Phi} -\sum_{i=1}^{k}\phi_{i} \\
  &=& ||\varpi_{k}-\varpi^* ||^2_{\Phi}-\sum_{i=1}^{k}\phi_{i}
\end{array}
\end{equation}
where the second inequality follows from the definition of $\phi_k$,  $\phi_{k+1} \geq ||\varpi_{k+1}-\varpi^* ||^2_{\Phi}-||\varpi_{k} -  \varpi^* ||^2_{\Phi}$.
As  $\{ ||\varpi_{k}-\varpi^* ||^2_{\Phi}-\sum_{i=1}^k\phi_{i} \}$ is a non-increasing sequence and bounded from below,
$\{ ||\varpi_{k}-\varpi^* ||^2_{\Phi}-\sum_{i=1}^k\phi_{i} \}$ converges.

Therefore, $ \{||\varpi_{k}-\varpi^* ||^2_{\Phi}\} $, being the sum of two convergent sequences $\{ ||\varpi_{k}-\varpi^* ||^2_{\Phi}-\sum_{i=1}^k\phi_{i} \} $ and $\{\sum_{i=1}^k\phi_{i}\}$, also converges.

\vskip 5mm
We are ready to show the convergence of Algorithm \ref{dal_3} using the results in {\it Part 1} and {\it Part 2}.
Since $\{\varpi_k\}$ is a bounded sequence,  it has a convergent subsequence $\{\varpi_{n_k}\}$ that converges to $\tilde{\varpi}^*$.　
Because $\lim_{k\rightarrow \infty}\varpi_{k+1}-\varpi_{k}=0$ by {\it Part 2},  we have
$\lim_{k \rightarrow \infty} \varpi_{n_{k}-1}-\varpi_{n_k}=0$ and $\lim_{k\rightarrow \infty} \varpi_{n_{k}+1}-\varpi_{n_k}=0$.
Pass to limiting point of $\{\varpi_{n_k}\}$, then we have $\tilde{\varpi}^*=T_1T_2\tilde{\varpi}^*$ because the righthand side of \eqref{cal_acc_1}-\eqref{cal_acc_2} is continuous. Hence, $\tilde{\varpi}^*\in zer(\bar{\mathfrak{A}}+\bar{\mathfrak{B}})$. Taking $\varpi^*=\tilde{\varpi}^*$ in \eqref{equ_thm_6_13} of {\it Part 1}, we also have $\{||\varpi_{k}-\tilde{\varpi}^* ||^2_{\Phi}\} $ converges by {\it Part 3}. Because a subsequence  $\{||\varpi_{n_k}-\tilde{\varpi}^* ||^2_{\Phi}\}$ converges to zero,  the whole sequence $\{||\varpi_{k}-\tilde{\varpi}^* ||^2_{\Phi}\} $  converges to zero. Therefore, the whole sequence of $\{\varpi_k\}$ converges to
$\tilde{\varpi}^*\in zer(\bar{\mathfrak{A}}+\bar{\mathfrak{B}})$. Resorting to Theorem \ref{thm_zero_is_correct} gives the desired result.
\hfill $\Box$

\begin{rem}
A sufficient and simple  choice of parameters to ensure the conditions in Theorem \eqref{thm_inertial} is  $\delta = \frac{1}{\beta}$, $\epsilon=\alpha$ and $0<\alpha < \sqrt{10}-3$.
In fact,  $\delta > \frac{1}{2\beta}$ implies that $2\beta\delta$ could be be any real number $\varrho > 1$.
If we take $\epsilon=\varsigma \alpha, \varsigma >0$,  then the quadratic inequality becomes
$\alpha^2- (2-3\varrho-\varsigma \varrho)\alpha+ 1-\varrho <0$. Since $1-\varrho<0$, $\alpha^2- (2-3\varrho-\varsigma \varrho)\alpha+ 1-\varrho$ takes value  strictly less than zero when $\alpha$ takes $0$.  By the continuity of quadratic equation, there always exists $0<\alpha<1$ that ensures the above quadratic inequality given any $\varrho>1$  and $\varsigma>0$.
\end{rem}

\section{Network Cournot game and simulation studies}\label{sec_cournot_simulaton}

There are various practical problems that can be well modeled by the game in \eqref{GM},
such as the river basin pollution game in \cite{krawczyk2},   the power market competition in \cite{krawczyk1}, plug-in electric vehicles charging management in \cite{lygeros2}, and communication network congestion game in \cite{shanbhag4}.
All above examples  can be regarded as the type of the  network Cournot game described below, which is a generalization of the network Cournot competition in \cite{cournotgame} by introducing additional market capacity constraints or equivalently globally shared coupling affine constraints. This type of network Cournot game with affine coupling constraints also appeared in the numerical studies of \cite{sayed}.

\subsection{Network Cournot game}\label{subsection_cournot}

Suppose that there are $N$ companies (players) with labels $F_1,\cdots,F_N$ and $m$ markets with labels $M_1,\cdots,M_m$.
Company $F_i$ decides its strategy to participate in the competition in $n_i$ markets by producing and delivering $ x_i \in \mathbf{R}^{n_i}$ amounts of products to the markets it connects with.
The production limitation of company $F_i$ is $x_i \in \Omega_i \subset \mathbf{R}^{n_i}$.
Company $F_i$ has a local matrix $A_i\in \mathbf{R}^{m\times n_i}$ that specifies which market it will participate in.
The $j$-th column of $A_i$, that is $[A_i]_{:j}$, has only one element being $1$ and all other elements being $0$, and
$[A_i]_{:j}$ has its $k$-th element being $1$ if and only if player $F_i$ delivers $[x_i]_j$ amount of production to the market $M_k$.
Therefore,  matrices $A_1,\cdots,A_N$ can be used to specify a  bipartite graph that represents the connections between the companies and the markets.
Denote  $n=\sum_{i=1}^N n_i$, $ \mathbf{x}= col(x_1,\cdots,x_N) \in \mathbf{R}^n$,  and   $A=[A_1,\cdots,A_N] \in \mathbf{R}^{m\times n}$.
Then $A\mathbf{x}\in \mathbf{R}^{m}= \sum_{i=1}^N A_i x_i$ is just the total product supply to all the markets given the action profile $\mathbf{x}$ of all the companies.
Market $M_j$ has  a maximal capacity of $r_j >0$, therefore,  it should be satisfied that $A\mathbf{x} \leq r$ where $r=col(r_1,\cdots,r_m) \in \mathbf{R}^m$.
Suppose that $P: \mathbf{R}^m \rightarrow \mathbf{R}^m$ is a price vector function that maps the total supply of each market to the corresponding market's price. Each company has also a local production cost function $c_i(x_i): \Omega_i \rightarrow \mathbf{R}$. Then the local objective function of company (player) $F_i$ is $f_i(x_i,\mathbf{x}_{-i}): c_i(x_i)-P^T(A\mathbf{x})A_ix_i$.

Overall, in this network Cournot game, each company needs to solve the following optimization problem given the other companies' profile $\mathbf{x}_{-i}$,
\begin{equation}
\begin{array}{ll}\label{network_cournot_game}
\min_{x_i \in \Omega_i } & \quad c_i(x_i)-P^T(A\mathbf{x})A_ix_i \\
s.t. &  A_ix_i \leq r- \sum_{j=1,j \neq i}^N A_jx_j
\end{array}
\end{equation}

Obviously, the above network Cournot game in \eqref{network_cournot_game} is a particular problem of game in \eqref{GM}.
Some practical decision problems in engineering networks can be well described by the network Cournot game in \eqref{network_cournot_game}, such as the rate control game in communication network (\cite{shanbhag4}) and the demand response game in smart grids (\cite{hu}).

\begin{exa}[Rate control game]
Consider a group of  source-destination pairs (nodes) in a communication network, that is  $\{S_1,...,S_N\}$, to decide their  data rates in a non-cooperative setting. The data  is transferred through a group of communication links (channels), that is $\{L_1,...,L_m\}$, and each link $L_j$ has a maximal data rate capacity of $c_j >0$.
Assume that an additional layer has decided the
routine table for each source-destination pair  $S_i$, which is encoded by $A_i \in \mathbf{R}^{m\times n_i}$.
Each column of $A_i$ has only one element being $1$ and all the other elements being zero, and  the $k$-th element of column $j$ is $1$ if $S_i$ utilizes the link $L_k$ and transfers data rate $[x_i]_j$ on link $L_k$.
 The local decision variable of $S_i$ is the data rate on each link  that it utilizes, denoted by $x_i \in \mathbf{R}^{n_i}$.
$x_i$  also has a local feasibility  constraint $x_i\in \Omega_i$. Denote $n=\sum_{i=1}^N n_i$,
$\mathbf{x}=col(x_1,\cdots,x_N)\in \mathbf{R}^{n}$, $A=[A_1,\cdots,A_N]\in \mathbf{R}^{m\times n}$ and $c=col(c_1,\cdots,c_m)\in \mathbf{R}^m$.
The total data rate on each link should be less than the capacity of that link:
$A\mathbf{x}\leq c. $
Given the data rate profile of all the nodes $\mathbf{x}$, the payoff function of $S_i$, $J_i(x_i,\mathbf{x}_{-i}): \mathbf{R}^n \rightarrow \mathbf{R}$, takes the form as
$J_i(x_i,\mathbf{x}_{-i})= -u_i(x_i) + \mathcal{D}^T(A\mathbf{x})A_ix_i$, where
$u_i(x_i): \Omega_i \rightarrow \mathbf{R}$ is the utility of source $S_i$, and $\mathcal{D}: R^{m}\rightarrow R^{m}$ is a delay function that maps the total data  rate on each link to the unit delay of that link.
 Thereby,   the data rate control game can be well described by the network Cournot game  in \eqref{network_cournot_game}.
\end{exa}

\begin{exa}[Demand response game]
Given a distribution network in power grids, suppose that there are $T$ time periods, and each period has a desirable minimal total load shedding $d_i>0$.
Suppose that there are $N$ load managers (energy management units or players) in the network, and each load manager  $i$ can decide a local vector $x_i\in \mathbf{R}^{t_i}$ as its local  load shedding vector in some specific time periods. Each load manager $i$ also has a local matrix $A_i \in \mathbf{R}^{T\times t_i}$ that specifies which time period player
$i$ will participate. For $j-$th column of $A_i$, it has one element being $1$ while all  other elements being zero.
The $k-$th element of the $j-$th column of $A_i$ is $1$ if load manager $i$ decides to decrease its load by $[x_i]_j$ at time $k$.
Denote $\mathbf{x}=col(x_1,\cdots, x_N)$, and $A=[A_1,\cdots,A_N] \in \mathbf{R}^{T\times \sum_{i=1}^N t_i}$, $d=col(d_1,\cdots, d_T)$.
Naturally, it is required  that the total load shedding of all the load managers should meet the minimal value, $A\mathbf{x }\geq d$.
Each player has a local feasible constraint $x_i \in \Omega_i \subset \mathbf{R}^{t_i}$, and a cost function $c_i(x_i): \Omega_i\rightarrow \mathbf{R}$ due to local load shedding. $P: \mathbf{R}^T \rightarrow \mathbf{R}^{T}$ is the payment price vector  function that maps total load shedding of each period to the payment price vector,  therefore, $P^T(Ax)A_ix_i$ is the payment awards of player $i$ for its load shedding.
The disutility function of player $i$ is $J_i(x_i,\mathbf{x}_{-i}) = c_i(x_i) - P^T(A\mathbf{x})A_ix_i $ given all the players' action profile $\mathbf{x}$.
All in all,  the demand response management  game is well described by the network Cournot game model in \eqref{network_cournot_game}.
\end{exa}

Moreover, the Assumptions \ref{assum1} and \ref{assum2} for Algorithm \ref{dal_1} and \ref{dal_3} can easily be satisfied for many practical cost functions and price functions.
For example,  take  company $F_i$'s  production cost function to be  a strongly convex  function with Lipschitz continuous gradients (A quadratic function $c_i(x_i)=x_i^T Q_i x_i+b_i^T x_i$ with
$Q_i\in \mathbf{R}^{n_i \times n_i}$ being a symmetric and positive definite matrix and $b_i\in \mathbf{R}^{n_i}$ is one possible choice).
The price of market $M_j$ is taken as the linear  function of the total supplying $p_j(\mathbf{x})= \bar{P}_j -  d_j [A\mathbf{x}]_{j}$ (known as a linear inverse demand function in economics) with $\bar{P}_j >0,   d_j>0$.
Denote $P=col(p_1,\cdots, p_m): \mathbf{R}^n \rightarrow \mathbf{R}^m$, $\bar{P}=col(\bar{P}_1,\cdots,\bar{P}_m)\in \mathbf{R}^m$,
$D=\diag\{d_1,\cdots,d_m\} \in \mathbf{R}^{m\times m}$. Then $P=\bar{P}-DA\mathbf{x}$ is the vector price function.
The payments of company $F_i$ by selling
product $x_i$ to the markets that it  connects with  is just $P^TA_ix_i$. Therefore, the objective function of company $F_i$ is,
\begin{equation}
\begin{array}{lll}\label{network_cournot_game_function}
J_i(x_i,\mathbf{x}_{-i}) &=&c_i(x_i)- (\bar{P}-DA\mathbf{x})^T A_ix_i \\
\;                       &=&c_i(x_i)-\bar{P}^TA_ix_i+ (\sum_{j=1}^N DA_jx_j)^T  A_ix_i\\
\;                       &=& c_i(x_i) + x_i^T A_i^T D A_i x_i  -A^T_i\bar{P}^T x_i  \\
\;                       & \quad & \qquad +  \sum_{j=1, j\neq i}^N x_j^T A_j^T DA_i x_i
\end{array}
\end{equation}
Denote $F(\mathbf{x})=col( \nabla_{x_1} J_1(x_1,\mathbf{x}_{-1}), \cdots, \nabla_{x_N} J_N(x_N, \mathbf{x}_{-N}) )$, $\nabla c(\mathbf{x})= col(\nabla c_1(x_1), \nabla c_2(x_2),\cdots, \nabla c_N(x_N))$ and $Q\in \mathbf{R}^{n\times n}$,
and $$Q=\left(
                \begin{array}{cccc}
                  2A_1^T D A_1 & A_1^T DA_2         & \cdots & A_1^T D A_N \\
                  A_2^T DA_1        & 2A_2^T D A_2  & \cdots & A_2^TDA_N \\
                  \cdots            & \cdots             & \cdots & \cdots  \\
                  A_N^TDA_1         & A_N^T D A_2        & \cdots &2A_N^T D A_N  \\
                \end{array}
              \right),$$
then
\begin{equation}
F(\mathbf{x})= \nabla c(\mathbf{x}) + Q\mathbf{x}-A^T\bar{P}.
\end{equation}

Notice that $Q$  can be written as
\begin{equation}
Q= \diag\{A^T_1DA_1, \cdots, A_N^T DA_N \}+S^T S
\end{equation}
where $S$ is a  block matrix defined as  $S=[\sqrt{D}A_1,\cdots, \sqrt{D}A_N] \in \mathbf{R}^{m\times n} $ and $\sqrt{D} =\diag\{\sqrt{d}_1,\cdots, \sqrt{d}_m\}$. Therefore, $Q$ is  positive semi-definite matrix.
Hence,  the Jacobian
matric of $F(\mathbf{x})$, $JF(\mathbf{x})= diag\{ \nabla^2 c_1(x_1), \nabla^2 c_2(x_2),\cdots, \nabla^2 c_N(x_N)\}+Q$ is postive definite since the cost functions  $c_i(x_i)$ are strongly convex.
Therefore,  $F(x)$ is strongly monotone\footnote{Proposition 2.3.2 of \cite{pang2}} and Lipschitz continuous, and Assumptions \ref{assum1} and \ref{assum2} are satisfied.

\begin{rem}
Notice that in the network Cournot game of \eqref{network_cournot_game}, each player $i$'s local objective function
$c_i(x_i)-P^T(A\mathbf{x})A_ix_i$ only depends on the decisions of the players that participate the same markets as player $i$.
Denote $\mathbb{Z}[a,b]$ as the set of integers from $a$ to $b$.
Mathematically, $j\in \mathcal{N}_i^f $ if and only if $\exists k\in\mathbb{ Z}[1,m]$, $p\in \mathbb{Z}[1,n_i]$ and $q\in \mathbb{Z}[1,n_j]$  such that $[A_i]_{kp}=1$ and $[A_j]_{kq}=1$. Since $A=[A_1,\cdots,A_N]$ is usually a sparse matrix, the interference graph $\mathcal{G}_f$ of network Cournot game also has sparse edge connections.
\end{rem}

\subsection{Simulation studies}

\begin{figure}
\begin{center}
\begin{tikzpicture}[->,>=stealth',shorten >=0.3pt,auto,node distance=2.1cm,thick,
  rect node/.style={rectangle, ball color={rgb:red,0;green,20;yellow,0},font=\sffamily,inner sep=1pt,outer sep=0pt,minimum size=12pt},
  wave/.style={decorate,decoration={snake,post length=0.1mm,amplitude=0.5mm,segment length=3mm},thick},
  main node/.style={shape=circle, ball color=green!1,text=black,inner sep=1pt,outer sep=0pt,minimum size=12pt},scale=0.85]

%\node[cloud, fill=blue!20, cloud puffs=16, cloud puff arc=50, minimum width=9cm, minimum height=3cm, aspect=3] at (3,0.7) {};
%\node[cloud, fill=blue!50, cloud puffs=16, cloud puff arc=60, minimum width=6cm, minimum height=2.6cm, aspect=3] at (3.4,-1.6) {};

% Plant nodes
  \foreach \place/\i in {{(-3.1,3.1)/1},{(-3.2,2.2)/2},
  {(-3.3,-0.3)/3},
  {(-3.4,-1.4)/4},
  {(-2.5,3.5)/5},
  {(-2.1,0.5)/6},
  {(-0.6,1.2)/7},
  {(-1.3,-0.6)/8},
  {(-1.4,-1.4)/9},
  {(0.1,3.5)/10},
  {(0.6,2.1)/11},
  {(1.8,0.2)/12},
  {(1.3,-1.3)/13},
  {(2.4,3.4)/14},
  {(1.5,1.5)/15},
  {(2.1,-1.1)/16},
  {(3.2,3.2)/17},
  {(3.3,2.3)/18},
  {(3.4,1.4)/19},
  {(3.5,0.5)/20}}
    \node[main node] (a\i) at \place {};

      \node at (-3.1,3.1){\rm \color{black}{$F_1$}};
      \node at (-3.2,2.2){\rm \color{black}{$F_2$}};
      \node at (-3.3,-0.3){\rm \color{black}{$F_3$}};
      \node at (-3.4,-1.4){\rm \color{black}{$F_4$}};
      \node at (-2.5,3.5){\rm \color{black}{$F_5$}};
      \node at (-2.1,0.5){\rm \color{black}{$F_6$}};
      \node at (-0.6,1.2){\rm \color{black}{$F_7$}};
      \node at (-1.3,-0.6){\rm \color{black}{$F_8$}};
      \node at (-1.4,-1.4){\rm \color{black}{$F_9$}};
      \node at (0.1,3.5){\rm \color{black}{$F_{10}$}};
      \node at (0.6,2.1){\rm \color{black}{$F_{11}$}};
      \node at (1.8,0.2){\rm \color{black}{$F_{12}$}};
      \node at (1.3,-1.3){\rm \color{black}{$F_{13}$}};
      \node at(2.4,3.4) {\rm \color{black}{$F_{14}$}};
      \node at (1.5,1.5){\rm \color{black}{$F_{15}$}};
      \node at (2.1,-1.1){\rm \color{black}{$F_{16}$}};
       \node at (3.2,3.2){\rm \color{black}{$F_{17}$}};
      \node at (3.3,2.3){\rm \color{black}{$F_{18}$}};
      \node at(3.4,1.4) {\rm \color{black}{$F_{19}$}};
      \node at (3.5,0.5){\rm \color{black}{$F_{20}$}};

  \foreach \place/\x in {{(-1.2,2.8)/1},{(-3,1)/2},{(-2.2,-1)/3},
    {(0,0)/4}, {(1,0)/5}, {(2,2)/6}, {(2.5,1)/7}}
  \node[rect node] (b\x) at \place {};

      \node at (-1.2,2.8){\rm \color{red}{$M_1$}};
      \node at (-3,1){\rm \color{red}{$M_2$}};
      \node at (-2.2,-1){\rm \color{red}{$M_3$}};
      \node at (0,0){\rm \color{red}{$M_4$}};
      \node at (1,0){\rm \color{red}{$M_5$}};
      \node at (2,2){\rm \color{red}{$M_6$}};
      \node at (2.5,1){\rm \color{red}{$M_7$}};

         \path[->,blue,thick]               (a1) edge (b1);

         \path[->,blue,thick]               (a2) edge (b1);
         \path[->,blue,thick]               (a2) edge (b2);

         \path[->,blue,thick]               (a3) edge (b2);

         \path[->,blue,thick]               (a4) edge (b3);

         \path[->,blue,thick]               (a5) edge (b1);

         \path[->,blue,thick]               (a6) edge (b1);
         \path[->,blue,thick]               (a6) edge (b2);
         \path[->,blue,thick]               (a6) edge (b3);
         \path[->,blue,thick]               (a6) edge (b4);

         \path[->,blue,thick]               (a7) edge (b4);

         \path[->,blue,thick]               (a8) edge (b3);
         \path[->,blue,thick]               (a8) edge (b4);

         \path[->,blue,thick]               (a9) edge (b3);

         \path[->,blue,thick]               (a10) edge (b1);
         \path[->,blue,thick]               (a10) edge (b4);
         \path[->,blue,thick]               (a10) edge (b6);

         \path[->,blue,thick]               (a11) edge (b4);
         \path[->,blue,thick]               (a11) edge (b5);

         \path[->,blue,thick]               (a12) edge (b5);

         \path[->,blue,thick]               (a13) edge (b5);

         \path[->,blue,thick]               (a14) edge (b6);

         \path[->,blue,thick]               (a15) edge (b5);
         \path[->,blue,thick]               (a15) edge (b6);
         \path[->,blue,thick]               (a15) edge (b7);

         \path[->,blue,thick]               (a16) edge (b5);
         \path[->,blue,thick]               (a16) edge (b7);

         \path[->,blue,thick]               (a15) edge (b7);

         \path[->,blue,thick]               (a17) edge (b6);
         \path[->,blue,thick]               (a17) edge (b7);

         \path[->,blue,thick]               (a18) edge (b7);
         \path[->,blue,thick]               (a19) edge (b7);
         \path[->,blue,thick]               (a20) edge (b7);
\end{tikzpicture}
\end{center}
\caption{Network cournot game: An edge from $F_i$ to $M_j$ on this graph implies that  $F_i$ participates in the competition in $M_j$.}\label{fig_network_cournot_game}
\end{figure}
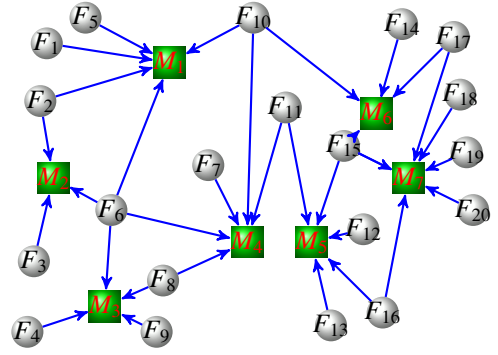

\begin{figure}
\begin{center}
\begin{tikzpicture}[->,>=stealth',shorten >=0.3pt,auto,node distance=2.1cm,thick,
  rect node/.style={rectangle,ball color=blue!10,font=\sffamily,inner sep=1pt,outer sep=0pt,minimum size=12pt},
  wave/.style={decorate,decoration={snake,post length=0.1mm,amplitude=0.5mm,segment length=3mm},thick},
  main node/.style={shape=circle,ball color=green!1,text=black,inner sep=1pt,outer sep=0pt,minimum size=12pt},scale=0.85]

%\node[cloud, fill=blue!20, cloud puffs=16, cloud puff arc=50, minimum width=9cm, minimum height=3cm, aspect=3] at (3,0.7) {};
%\node[cloud, fill=blue!50, cloud puffs=16, cloud puff arc=60, minimum width=6cm, minimum height=2.6cm, aspect=3] at (3.4,-1.6) {};

% Plant nodes
  \foreach \place/\i in {{(-3.1,3.1)/1},{(-3.2,2.2)/2},
  {(-3.3,-0.3)/3},
  {(-3.4,-1.4)/4},
  {(-1.5,3.8)/5},
  {(-2.1,0.5)/6},
  {(-0.4,1.1)/7},
  {(-1.3,-0.6)/8},
  {(-1.4,-1.4)/9},
  {(0.1,3.5)/10},
  {(0.3,-1.1)/11},
  {(1.8,0.2)/12},
  {(1.3,-1.3)/13},
  {(2.4,3.9)/14},
  {(1.5,1.5)/15},
  {(2.8,-0.8)/16},
  {(2.5,3)/17},
  {(3.5,2.9)/18},
  {(4.2,1.4)/19},
  {(3.8,0.1)/20}}
    \node[main node] (a\i) at \place {};

      \node at (-3.1,3.1){\rm \color{black}{$F_1$}};
      \node at (-3.2,2.2){\rm \color{black}{$F_2$}};
      \node at (-3.3,-0.3){\rm \color{black}{$F_3$}};
      \node at (-3.4,-1.4){\rm \color{black}{$F_4$}};
      \node at (-1.5,3.8){\rm \color{black}{$F_5$}};
      \node at (-2.1,0.5){\rm \color{black}{$F_6$}};
      \node at (-0.4,1.1){\rm \color{black}{$F_7$}};
      \node at (-1.3,-0.6){\rm \color{black}{$F_8$}};
      \node at (-1.4,-1.4){\rm \color{black}{$F_9$}};
      \node at (0.1,3.5){\rm \color{black}{$F_{10}$}};
      \node at (0.3,-1.1){\rm \color{black}{$F_{11}$}};
      \node at (1.8,0.2){\rm \color{black}{$F_{12}$}};
      \node at (1.3,-1.3){\rm \color{black}{$F_{13}$}};
      \node at(2.4,3.9) {\rm \color{black}{$F_{14}$}};
      \node at (1.5,1.5){\rm \color{black}{$F_{15}$}};
      \node at (2.8,-0.8){\rm \color{black}{$F_{16}$}};
       \node at (2.5,3){\rm \color{black}{$F_{17}$}};
      \node at (3.5,2.9){\rm \color{black}{$F_{18}$}};
      \node at(4.2,1.4) {\rm \color{black}{$F_{19}$}};
      \node at (3.8,0.1){\rm \color{black}{$F_{20}$}};

         \path[<->,blue,thick]               (a1) edge (a2);
         \path[<->,blue,thick]               (a1) edge (a5);
         \path[<->,blue,thick]               (a1) edge (a6);
         \path[<->,blue,thick]               (a1) edge (a10);

         \path[<->,blue,thick]               (a2) edge (a3);
         \path[<->,blue,thick]               (a2) edge (a5);
         \path[<->,blue,thick]               (a2) edge (a6);
         \path[<->,blue,thick]               (a2) edge (a10);

         \path[<->,blue,thick]               (a3) edge (a2);
          \path[<->,blue,thick]              (a3) edge (a6);

          \path[<->,blue,thick]               (a4) edge (a8);
          \path[<->,blue,thick]              (a4) edge (a6);

        \path[<->,blue,thick]               (a5) edge (a1);
         \path[<->,blue,thick]               (a5) edge (a2);
         \path[<->,blue,thick]               (a5) edge (a6);
         \path[<->,blue,thick]               (a5) edge (a10);

          \path[<->,blue,thick]               (a6) edge (a4);
           \path[<->,blue,thick]               (a6) edge (a8);
           \path[<->,blue,thick]               (a6) edge (a9);

           \path[<->,blue,thick]               (a6) edge (a8);
           \path[<->,blue,thick]               (a6) edge (a7);
           \path[<->,blue,thick]               (a6) edge (a10);
           \path[<->,blue,thick]               (a6) edge (a11);

           \path[<->,blue,thick]               (a7) edge (a8);
           \path[<->,blue,thick]               (a7) edge (a6);
           \path[<->,blue,thick]               (a7) edge (a10);
           \path[<->,blue,thick]               (a7) edge (a11);

           \path[<->,blue,thick]               (a8) edge (a9);
            \path[<->,blue,thick]              (a8) edge (a10);
           \path[<->,blue,thick]               (a8) edge (a11);

             \path[<->,blue,thick]               (a9) edge (a4);

           \path[<->,blue,thick]               (a10) edge (a11);
           \path[<->,blue,thick]               (a10) edge (a15);
           \path[<->,blue,thick]               (a10) edge (a14);
           \path[<->,blue,thick]               (a10) edge (a17);

           \path[<->,blue,thick]               (a11) edge (a15);
           \path[<->,blue,thick]               (a11) edge (a12);
            \path[<->,blue,thick]               (a11) edge (a13);
           \path[<->,blue,thick]               (a11) edge (a16);

           \path[<->,blue,thick]               (a12) edge (a11);
           \path[<->,blue,thick]               (a12) edge (a15);
           \path[<->,blue,thick]               (a12) edge (a13);
           \path[<->,blue,thick]               (a12) edge (a16);

           \path[<->,blue,thick]               (a13) edge (a11);
           \path[<->,blue,thick]               (a13) edge (a15);
           \path[<->,blue,thick]               (a13) edge (a12);
           \path[<->,blue,thick]               (a13) edge (a16);

              \path[<->,blue,thick]             (a14) edge (a10);
           \path[<->,blue,thick]             (a14) edge (a15);
           \path[<->,blue,thick]             (a14) edge (a17);

           \path[<->,blue,thick]             (a15) edge (a16);
           \path[<->,blue,thick]             (a15) edge (a17);
            \path[<->,blue,thick]             (a15) edge (a18);
             \path[<->,blue,thick]             (a15) edge (a19);
              \path[<->,blue,thick]             (a15) edge (a20);

            \path[<->,blue,thick]               (a10) edge (a15);
            \path[<->,blue,thick]               (a10) edge (a14);
            \path[<->,blue,thick]               (a10) edge (a17);

            \path[<->,blue,thick]               (a16) edge (a17);
            \path[<->,blue,thick]               (a16) edge (a18);
            \path[<->,blue,thick]               (a16) edge (a19);
            \path[<->,blue,thick]               (a16) edge (a20);

            \path[<->,blue,thick]               (a17) edge (a18);
            \path[<->,blue,thick]               (a17) edge (a19);
            \path[<->,blue,thick]               (a17) edge (a20);

            \path[<->,blue,thick]               (a18) edge (a19);
            \path[<->,blue,thick]               (a18) edge (a20);

            \path[<->,blue,thick]               (a19) edge (a20);
\end{tikzpicture}
\end{center}
\caption{Interference graph $\mathcal{G}_f$: An undirected edge from $F_i$ to $F_j$ on this graph means that the objective function of $F_i$ depends on $x_j$, and vice versa. Company $F_i$ is assumed to be able to observe the decisions of its neighbors on this graph .}\label{fig_inference_game}
\end{figure}
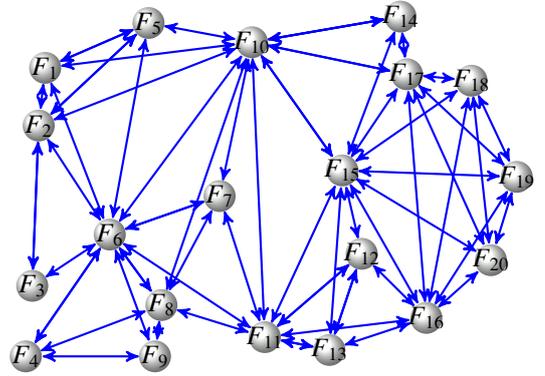

\begin{figure}
\begin{center}
\begin{tikzpicture}[->,>=stealth',shorten >=0.3pt,auto,node distance=2.1cm,thick,
  rect node/.style={rectangle,ball color=blue!10,font=\sffamily,inner sep=1pt,outer sep=0pt,minimum size=12pt},
  wave/.style={decorate,decoration={snake,post length=0.1mm,amplitude=0.5mm,segment length=3mm},thick},
  main node/.style={shape=circle,ball color=green!1,text=black,inner sep=1pt,outer sep=0pt,minimum size=12pt},scale=0.85]

%\node[cloud, fill=blue!20, cloud puffs=16, cloud puff arc=50, minimum width=9cm, minimum height=3cm, aspect=3] at (3,0.7) {};
%\node[cloud, fill=blue!50, cloud puffs=16, cloud puff arc=60, minimum width=6cm, minimum height=2.6cm, aspect=3] at (3.4,-1.6) {};

% Plant nodes
  \foreach \place/\i in {{(-3,2)/1},{(-2,2)/2},
  {(-1,2)/3},
  {(0,2)/4},
  {(1,2)/5},
  {(2,2)/6},
  {(3,2)/7},
  {(3,1)/8},
  {(3,0)/9},
  {(3,-1)/10},
  {(3,-2)/11},
  {(2,-2)/12},
  {(1,-2)/13},
  {(0,-2)/14},
  {(-1,-2)/15},
  {(-2,-2)/16},
  {(-3,-2)/17},
  {(-3,-1)/18},
  {(-3,-0)/19},
  {(-3,1)/20}}
    \node[main node] (a\i) at \place {};

      \node at (-3,2){\rm \color{black}{$F_1$}};
      \node at (-2,2){\rm \color{black}{$F_2$}};
      \node at (-1,2){\rm \color{black}{$F_3$}};
      \node at (0,2){\rm \color{black}{$F_4$}};
      \node at (1,2){\rm \color{black}{$F_5$}};
      \node at (2,2){\rm \color{black}{$F_6$}};
      \node at (3,2){\rm \color{black}{$F_7$}};
      \node at (3,1){\rm \color{black}{$F_8$}};
      \node at (3,0){\rm \color{black}{$F_9$}};
      \node at (3,-1){\rm \color{black}{$F_{10}$}};
      \node at (3,-2){\rm \color{black}{$F_{11}$}};
      \node at (2,-2){\rm \color{black}{$F_{12}$}};
      \node at (1,-2){\rm \color{black}{$F_{13}$}};
      \node at (0,-2) {\rm \color{black}{$F_{14}$}};
      \node at (-1,-2){\rm \color{black}{$F_{15}$}};
      \node at (-2,-2){\rm \color{black}{$F_{16}$}};
      \node at (-3,-2){\rm \color{black}{$F_{17}$}};
      \node at (-3,-1){\rm \color{black}{$F_{18}$}};
      \node at(-3,0) {\rm \color{black}{$F_{19}$}};
      \node at (-3,1){\rm \color{black}{$F_{20}$}};

        \path[<->,blue,thick]               (a1) edge (a2);
         \path[<->,blue,thick]               (a2) edge (a3);
         \path[<->,blue,thick]               (a3) edge (a4);
         \path[<->,blue,thick]               (a4) edge (a5);

         \path[<->,blue,thick]               (a5) edge (a6);
         \path[<->,blue,thick]               (a6) edge (a7);
         \path[<->,blue,thick]               (a7) edge (a8);
         \path[<->,blue,thick]               (a8) edge (a9);
           \path[<->,blue,thick]               (a9) edge (a10);
         \path[<->,blue,thick]               (a10) edge (a11);
             \path[<->,blue,thick]             (a11) edge (a12);
            \path[<->,blue,thick]               (a12) edge (a13);
            \path[<->,blue,thick]               (a13) edge (a14);
            \path[<->,blue,thick]               (a14) edge (a15);

            \path[<->,blue,thick]               (a15) edge (a16);
            \path[<->,blue,thick]               (a16) edge (a17);
            \path[<->,blue,thick]               (a17) edge (a18);

            \path[<->,blue,thick]               (a18) edge (a19);
            \path[<->,blue,thick]               (a19) edge (a20);
             \path[<->,blue,thick]               (a1) edge (a20);

              \path[<->,blue,thick]               (a2) edge (a15);
             \path[<->,blue,thick]               (a6) edge (a13);

\end{tikzpicture}
\end{center}
\caption{Multiplier graph $\mathcal{G}_{\lambda}$: Company $F_i$ and $F_j$ are able to exchange their local $\{\lambda_i,z_i\}$ and $\{\lambda_j,z_j\}$ if there exists an edge between them on this graph. }\label{fig_multiplier_graph}
\end{figure}
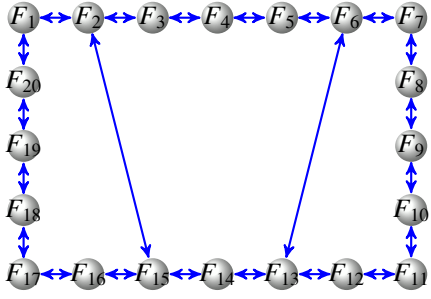

In the studies, we adopt a similar simulation setting as \cite{sayed} without considering the stochastic factors.
Consider $20$ companies and $7$ markets, and the connection relationship between the companies and the markets is depicted in
Figure \ref{fig_network_cournot_game}. If there is an edge from $F_i$ to $M_j$ in Figure \ref{fig_network_cournot_game}, then company $F_i$ participates the competition in market $M_j$ by producing and delivering products to market $M_j$.
Each company $F_i$ has a local constraint as $ \mathbf{0}< x_i < \Theta_i $ where each component of $\Theta_i$ is randomly drawn from $(10,25)$. Each market $M_j$ has a maximal capacity  of $r_j >0, j=1,\cdots, 7 $ and $r_j$ is  randomly  drawn from $(20,80)$. The local objective function is taken as \eqref{network_cournot_game_function}. The local cost function of company $i$ is $c_i(x_i)=\pi_i(\sum_{j=1}^{n_i}[x_{i}]_j)^2+b_i^T x_i$ which is a strongly convex function with Lipschitz continuous gradient. Here $\pi_i$ is randomly drawn from $(1,8)$, and each component of $b_i$ is randomly drawn from $(1,4)$.
The price function is taken as the linear function $ P=\bar{P}-DA\mathbf{x}$, and  $\bar{P}_j$ and $d_j$ are  randomly drawn from $(250,500)$ and $(1,5)$, respectively.

\begin{figure}
  \centering
  \includegraphics[width=3.5in]{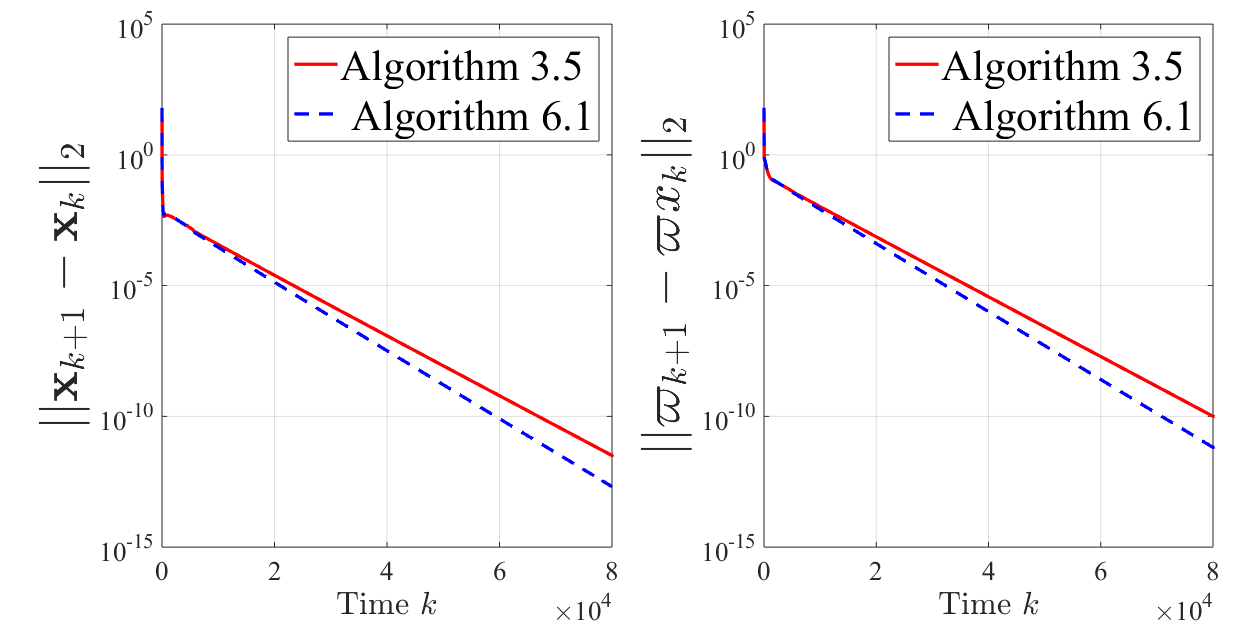}\\
  \caption{ The trajectories of $||\mathbf{x}_{k+1}-\mathbf{x}_k ||_2$ and $||\varpi_{k+1}-\varpi_{k} ||_2$ generated with Algorithm \ref{dal_1} and \ref{dal_3}: This shows the convergence of both algorithms, and the superior convergence speed of the algorithm with inertia \ref{dal_3}.\label{fig_sim_1}  }
\end{figure}

\begin{figure}
  \centering
  \includegraphics[width=3.5in]{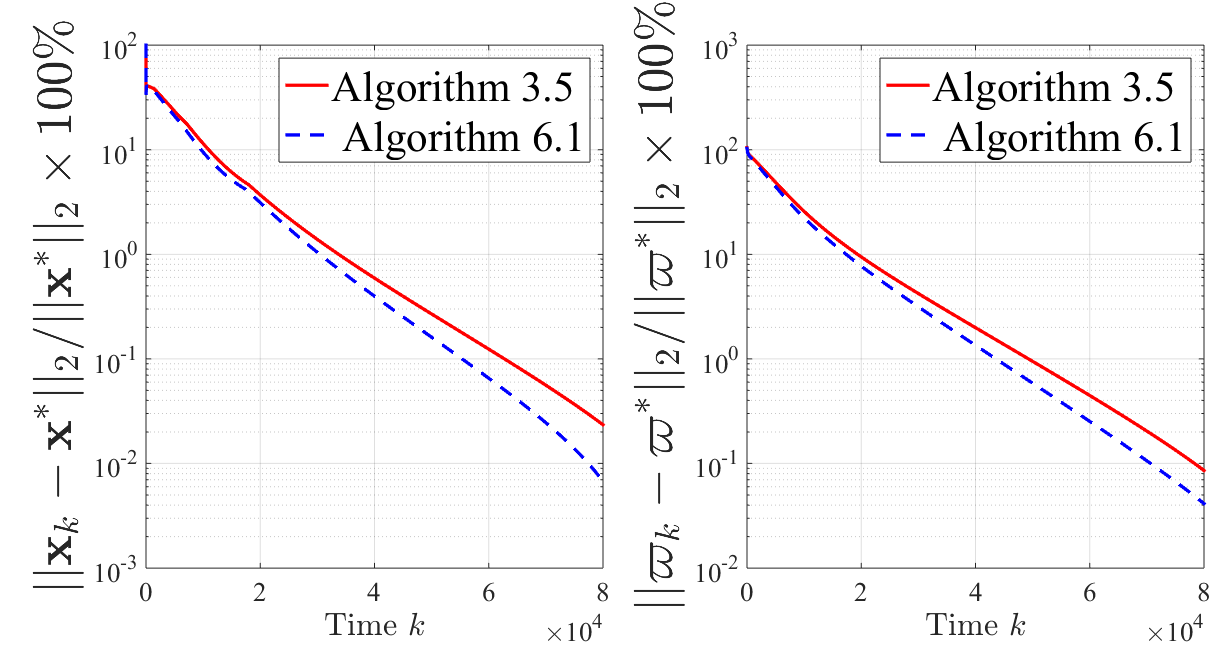}\\
  \caption{ The trajectories of $\frac{||\mathbf{x}_{k}-\mathbf{x}^* ||_2}{||x^*||_2}\times 100\%$ and $\frac{||\varpi_{k}-\varpi^* ||_2}{||\varpi^* ||_2}\times 100\%$ generated with Algorithm \ref{dal_1} and \ref{dal_3}: This also shows the convergence of both algorithms, and the superior convergence speed of Algorithm \ref{dal_3}. }\label{fig_sim_6}
\end{figure}

\begin{figure}
  \centering
  \includegraphics[width=3.5in]{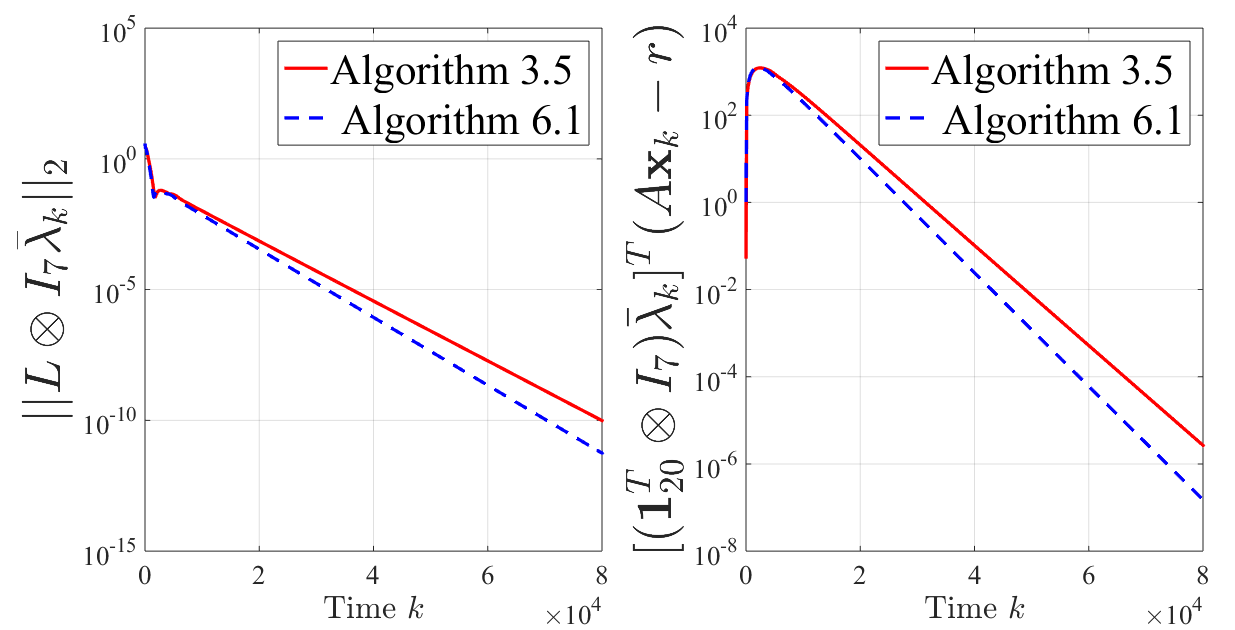}\\
  \caption{The trajectory of $||(L\otimes I_7)\bar{\lambda}_{k} ||_2$ and $[(\mathbf{1}^T_{20}\otimes I_7)\bar{\lambda}_k]^T (A\mathbf{x}_k-r)$: The left sub-figure shows that $\bar{\lambda}_{k}$ converges to the null space of $L\otimes I_7$, hence the local multipliers of all companies reach consensus asymptotically. $(\mathbf{1}^T_{20}\otimes I_7)\bar{\lambda}_k$ is the averaged vector $\sum_{i=1}^{20} \lambda_{i,k}$, therefore, the right sub-figure shows that the complementary condition $\langle \lambda^*, Ax^*-r \rangle=0$  is  asymptotically satisfied.  } \label{fig_sim_2}
\end{figure}

\begin{figure}
  \centering
  \includegraphics[width=3.5in]{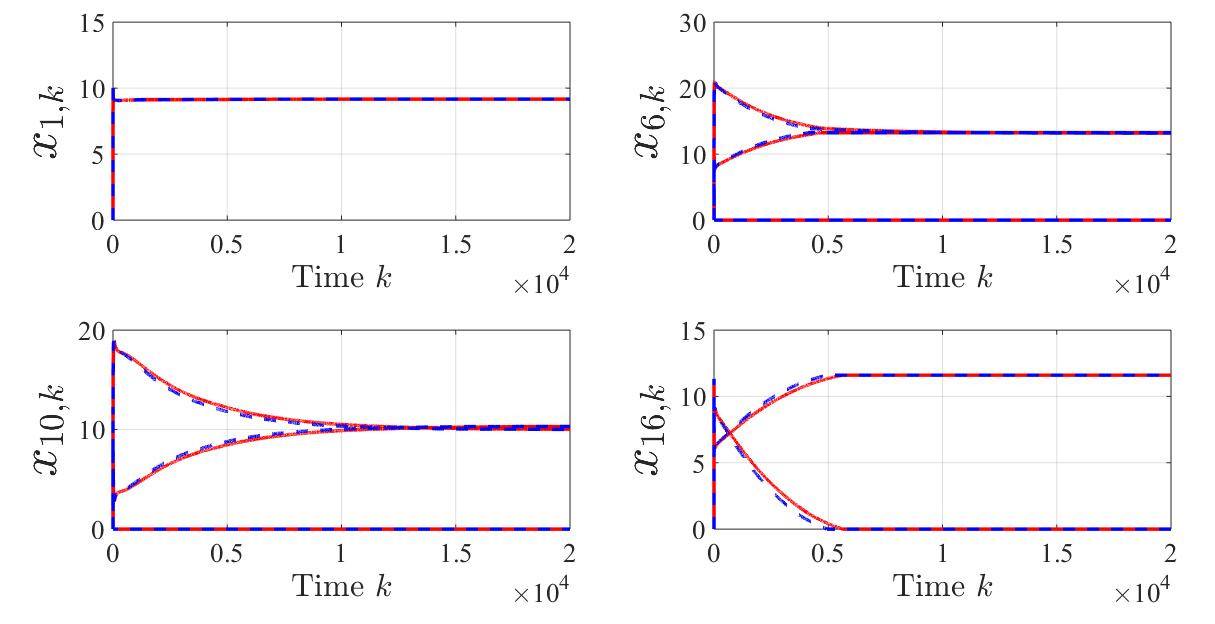}\\
  \caption{The trajectories of local decisions $x_{i,k}$ of $F_1$, $F_6$, $F_{10}$ and $F_{16}$.
   The solid red lines are generated with Algorithm \ref{dal_1} and dashed lines are generated with Algorithm \ref{dal_3}}\label{fig_sim_3}
\end{figure}

\begin{figure}
  \centering
  \includegraphics[width=3.5in]{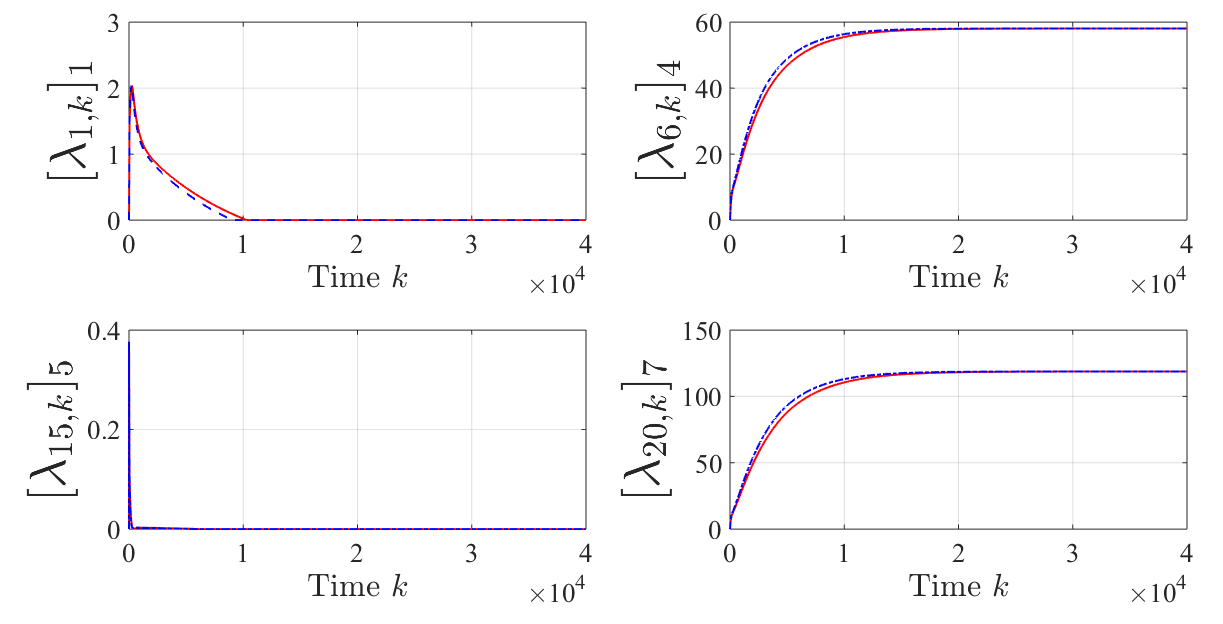}\\
  \caption{The trajectories of some components of local multiplier $[\lambda_{i,k}]_j$ where solid red lines are generated with Algorithm \ref{dal_1} and dashed lines are generated with Algorithm \ref{dal_3}. $[\lambda_{i,k}]_j$ is the $j-$th component of local multiplier $\lambda_i$ at iteration time $k$.  Here we select local multipliers of $F_1$, $F_6$, $F_{10}$ and $F_{16}$} \label{fig_sim_4}
\end{figure}

With Figure \ref{fig_network_cournot_game} and the definition of objective function in \eqref{network_cournot_game_function}, the interference graph  $\mathcal{G}_{f}$  can be easily obtained and is depicted in Figure \ref{fig_inference_game}. Meanwhile, we adopte
the multiplier graph $\mathcal{G}_{\lambda}$ shown in Figure \ref{fig_multiplier_graph}. The weighted adjacency matrix $W=[w_{ij}]$ of multiplier graph $\mathcal{G}_{\lambda}$ has all its nonzero elements to be $1$.

With Figure \ref{fig_network_cournot_game} and the definition of objective function in \eqref{network_cournot_game_function}, the interference graph  $\mathcal{G}_{f}$  can be easily obtained and is depicted in Figure \ref{fig_inference_game}. Meanwhile, we adopte
the multiplier graph $\mathcal{G}_{\lambda}$ shown in Figure \ref{fig_multiplier_graph}. The weighted adjacency matrix $W=[w_{ij}]$ of multiplier graph $\mathcal{G}_{\lambda}$ has all its nonzero elements to be $1$.

Set the step-sizes in Algorithm \ref{dal_1} as $\tau_i=0.03, \nu_i=0.2, \sigma=0.02$ for all companies,
and for  Algorithm \ref{dal_3} set $\alpha=0.12$ while other step-sizes are the same as Algorithm \ref{dal_1}.  The initial starting points $x_{i,0}, \lambda_{i,0}$ and $z_{i,0}$ of both algorithms are set to be zeros.

The trajectories of  selected  algorithm  performance indexes, including $||\mathbf{x}_{k+1}-\mathbf{x}_{k} ||_2$, $||\varpi_{k+1}-\varpi_{k} ||_2$, $\frac{||\mathbf{x}_k-\mathbf{x}^* ||}{|| \mathbf{x}^*||}\times 100\%$,
 $\frac{||\varpi_k-\varpi^* ||}{||\varpi^* ||}\times 100 \%$, $||L\otimes I_7 \bar{\lambda}_k　||_2$ and $[(\mathbf{1}^T_{20}\otimes I_7)\bar{\lambda}_k]^T (A\mathbf{x}_k-r)$, are shown in Figures \ref{fig_sim_1}, \ref{fig_sim_6} and  \ref{fig_sim_2}. The trajectories of the local decisions $x_{i,k}$ of some companies are shown in Figure \ref{fig_sim_3}, and the trajectories of $[\lambda_{i,k}]_j$, which stands for the $j$-th  component of local Lagrangian multiplier $\lambda_{i,k}$, are shown in Figure \ref{fig_sim_4}.
%The trajectories of $k ||\varpi_{k+1}-\varpi_k  ||_2$ and $k^2||\varpi_{k+1}-\varpi_{k}　||_2$ are also given in Figure \ref{fig_sim_5} for a numerical illustration of convergence speed.

\section{Conclusions}\label{sec_concluding}

In this paper, we proposed a  primal-dual distributed algorithm based on operator splitting methods for iterative computation of a variational GNE in noncooperative games with globally shared affine coupling constraints. The algorithm is motivated  by the forward-backward operator splitting method for finding zeros of a sum of monotone operators.  Each player only needs to knows its local information, especially a block of the affine coupling constraints. The proposed algorithm is proved to converge with fixed step-sizes under some mild assumptions by exploiting properties of  composition of averaged operators. Furthermore, a distributed algorithm with inertia is also proposed and analyzed for possible acceleration of convergence speed.
Numerical simulation  studies for a network Cournot game  demonstrate the efficiency of the proposed algorithms and the superior convergence speed of the inertial algorithm.

 Many challenging and exciting topics are still open for distributed NE/GNE seeking. Here we only list some problems with  probable solution hints.
 Finding all the generalized Nash equilibria has its only interests, and this could be  partially solved  by combining the design in this paper with the {\it parameterized variational inequality method} in \cite{pseng}.
 The algorithm  requires that each player is able to observe all its neighbors' decisions through the interference graph $\mathcal{G}_f$. This assumption could be relaxed by adopting the local consensus dynamics in \cite{pavel4}, and then it could only be  required that the players were able to  observe parts of its neighbors' decisions through a {\it maximal triangle-free spanning subgraph} of $\mathcal{G}_f$.
 The methodology of this  paper could be extended for stochastic GNE seeking with noisy  gradient observations and noisy information sharing by resorting to the  {\it stochastic forward-backward splitting algorithm} in \cite{vu}.
 The strong monotonicity  assumption on the pseudo-gradient might be relaxed to monotonicity  assumption by utilizing the {\it forward-backward-forward splitting method} in \cite{combettes2}.

\end{document}